\documentclass[a4paper,11pt]{article}
\usepackage{epsf,epsfig}
\usepackage{graphics,color}
\usepackage{amsmath}
\usepackage{amssymb}
\usepackage{cite}
\usepackage{verbatim}
\usepackage{float}
\usepackage{graphicx}
\usepackage{amsthm}
\usepackage{textcomp}
\usepackage{subfig}
\usepackage{hyperref}
\usepackage{yhmath}
\usepackage{mathrsfs}
\newif\ifdraft
\draftfalse 

\numberwithin{equation}{section}
\newtheorem{theorem}{Theorem}[section]
\newtheorem{Proposition}[theorem]{Proposition}
\newtheorem{Corollary}[theorem]{Corollary}
\newtheorem{lemma}[theorem]{Lemma}

\newtheorem{Remark}[theorem]{Remark}
\newtheorem{Definition}[theorem]{Definition}

\newcommand{\area}{\mathcal A}

\newcommand{\connessionebordo}{\omega}
\newcommand{\connessionepoloi}{\gamma_{z_i}}
\newcommand{\curl}{\,{\rm Curl\ }}
\newcommand{\Det}{\,{\rm Det\ }}
\newcommand{\dist}{d}
\renewcommand{\div}{{\rm div}}
\newcommand{\Div}{\,{\rm Div\ }}
\newcommand{\dom}{{\rm Dom}}
\newcommand{\dOm}{\partial \Om}
\newcommand{\dovOm}{d_\Om}
\newcommand{\eps}{\epsilon}
\newcommand{\effeell}{T_{r}(\connessionebordo_\ell)}

\newcommand{\findistr}{T}
\newcommand{\diagonal}{{\rm Diag}_\Om}
\newcommand{\diagonals}{{\rm Diag}_\Om^\star}
\newcommand{\grad}{\nabla}

\newcommand{\indice}{k}
\newcommand{\indicefinalesomma}{N}
\newcommand{\indicefinalesommalemma}{m}
\newcommand{\indicefinalesommateorema}{N}
\newcommand{\indicefinalesommaproposition}{N}
\newcommand{\indicegenerico}{j}
\newcommand{\ingrassatoconnessionepoloi}{T_{2r}(\gamma_{z_i})}

\newcommand{\infindistr}{\Lambda}
\newcommand{\jump}[1]{\text{{\rm \textlbrackdbl}}{#1}\text{{\rm \textrbrackdbl}}}
\newcommand{\NN}{\mathbb N}

\newcommand{\Om}{\Omega}
\newcommand{\outofballs}{G_\eps}
\newcommand{\piccoloraggio}{{\small R}}
\newcommand{\punto}{c}
\newcommand{\R}{\mathbb{R}}
\newcommand{\res}      {\mathop{\hbox{\vrule height 7pt width .5pt depth
                        0pt\vrule height .5pt width 6pt depth 0pt}}\nolimits}
\newcommand{\rilarea}{\overline{\mathcal A}}
\newcommand{\rilrilarea} {\overline {\overline{\mathcal A}}}
\newcommand{\Suno}{{\mathbb S}^1}
\newcommand{\supp}{\,{\rm supp \ }}
\newcommand{\testfunctions}{\mathcal C^\infty_c}
\newcommand{\zerocurrfin}{\mathcal R_f}

\title{Upper bounds for the relaxed area of $\Suno$-valued
Sobolev maps 
and its countably subadditive interior envelope}

\author{Giovanni Bellettini\footnote{
		Dipartimento di Ingegneria dell'Informazione e Scienze Matematiche, Universit\`a di Siena, 53100 Siena, Italy,
		and International Centre for Theoretical Physics ICTP,
		Mathematics Section, 34151 Trieste, Italy.
		E-mail: bellettini@diism.unisi.it
	}
	\and
	Riccardo Scala\footnote{ 
		Dipartimento di Ingegneria dell'Informazione e Scienze Matematiche, Universit\`a di Siena, 53100 Siena, Italy.
		E-mail: riccardo.scala@unisi.it}
	\and
Giuseppe Scianna \footnote{	Dipartimento di Ingegneria dell'Informazione e Scienze Matematiche, Universit\`a di Siena, 53100 Siena, Italy.
	E-mail: giuseppe.scianna@unisi.it
}
}

\begin{document}

\maketitle
\begin{abstract}
Given a bounded open connected Lipschitz set $\Om \subset \R^2$,  
we show that the relaxed Cartesian
area functional $\rilarea(u, \Om)$
of a map $u\in W^{1,1}(\Om;\Suno)$ is finite,
and provide a useful upper bound for its value. Using
this estimate, 
we prove a modified version of a De Giorgi conjecture \cite{DG} adapted to $W^{1,1}(\Om; \Suno)$,
on the largest countably subadditive set function
$\rilrilarea(u, \cdot)$ smaller than or equal to $\rilarea(u, \cdot)$.

\end{abstract}


\vspace{2mm}

\noindent {\bf AMS (MOS) 2020 Subject Clas\-si\-fi\-ca\-tion:}  49J45, 49Q05, 49Q15, 28A75.

\section{Introduction}
Let  $\Omega \subset \R^2$ be a bounded open set.
For a given $v\in C^1(\Om;\R^2)$ we indicate by 
\begin{align}
	\area(v,\Om):=\int_\Om\sqrt{1+|\nabla v|^2+|Jv|^2}dx
\end{align}
the classical $2$-dimensional area of the graph $G_v=\{(x,y)\in \Om\times \R^2:y=v(x)\}$ of $v$,
where $Jv=\frac{\partial v_1}{\partial x_1}\frac{\partial v_2}{\partial x_2}-\frac{\partial v_2}{\partial x_1}\frac{\partial v_1}{\partial x_2}$
denotes the Jacobian determinant of $v$. For any $u\in L^1(\Om;\R^2)$
we consider the $L^1$-relaxed area of the graph of $u$, namely
\begin{align}\label{def_arearel}
\rilarea
(u,\Om):=\inf\{\liminf_{k\rightarrow +\infty}\area(v_k;\Om),\;v_k\in C^1(\Om;\R^2),\;v_k\rightarrow u\text{ in }L^1(\Om;\R^2)\}.
\end{align}
It is well known that, when $v$ is scalar valued,
 the study of the relaxed area 
is crucial in the study of the Cartesian Plateau problem \cite{Giu:84}.
The characterization of the domain $\dom(\rilarea(\cdot, \Omega))$ of 
$\rilarea(\cdot, \Omega)$, 
and the computation of its corresponding 
values,
seem at the moment out of reach, due to the presence of highly nonlocal
phenomena. 
More specifically, for a given\footnote{This is true for general maps, 
unless some specific cases which trivialize the functional 
(see \cite{ADM} for details).} 
$u\in L^1(\Om;\R^2)$, the set function 
$\Om\supseteq A\rightarrow \rilarea(u,A)$ turns out to be not subadditive 
when restricted to open sets. In particular
$\rilarea(u, \cdot)$ 
is not a measure, and 
thus it
cannot be represented in integral form; for this reason, 
only a few partial results are available 
(see e.g. \cite{BEPS,BePaTe:15,BePaTe:16}). In 
these references it is shown that nonlocality is due to at least 
two reasons: one is the presence of singularities
in the map $u$; the other one is the possible interaction of such 
singularities with $\partial \Om$. In both cases, it appears that,
in general, 
interesting and rather involved Plateau-type problems must be solved,
in order to get the exact value of $\rilarea(u, \Om)$ (see
the discussion below on the maps $u_V$ and $u_T$). So, the computation
of $\rilarea(u, \Om)$ is, in general, quite difficult; on the other
hand, looking for upper bounds that do not take into account
the above mentioned Plateau problems in full generality, seems
realistic\footnote{Notice that,
	if one replaces the $L^1$ convergence 
	in \eqref{def_arearel} with stronger topologies, some 
	sharp estimates can be given (see
	for instance \cite{BCS,BCS2,Mu} where the strict convergence in $BV$ has been investigated).
}

In this paper we are
 concerned with maps in 
$$
	W^{1,1}(\Om;\Suno):=\{u\in W^{1,1}(\Om;\R^2):|u|=1\text{ a.e. in }\Om\},
$$
where $\Suno=\{(x_1,x_2)\in\R^2:x_1^2+x_2^2=1\}$.
Given a distribution $\infindistr
\in \mathcal D'(\Om)$, let us introduce the quantity
\begin{equation}\label{eq:def_M}
\Vert \infindistr\Vert_{{\rm flat},\alpha}:=\sup
\Big\{
\langle \infindistr,\varphi\rangle: 
\varphi\in {\rm Lip}_0(\Om), 
\|\varphi\|_{L^\infty(\Om)}\leq 1, \ 		
\alpha \|\nabla \varphi\|_{L^\infty(\Om)}\leq 1
\Big\},
\end{equation}
where 
\begin{equation}\label{eq:alpha}
\alpha:= \frac{\vert B_1\vert}{\vert \partial B_1\vert} = \frac{1}{2}
\end{equation}
and ${\rm Lip}_0(\Om)$ are the Lipschitz functions on $\Om$ vanishing 
on $\partial \Om$. 
Our first result (Section \ref{sec:main_results}) reads as follows: 
\begin{theorem}\label{teo:upper_bound_W11Suno} 
	Let $\Om\subset\R^2$ be a  bounded open set with Lipschitz boundary 
and $u\in W^{1,1}(\Om;\Suno)$. Then 
		\begin{align}\label{eq:teo1}
		\rilarea(u,\Om)\leq \int_{\Om}\sqrt{1+|\nabla u|^2}dx+
\Vert \Det(\nabla u)\Vert_{{\rm flat},\alpha} < +\infty.
	\end{align}
In particular
$$W^{1,1}(\Om;\Suno)\subset \dom(\rilarea(\cdot,\Om)).$$
\end{theorem}

Estimate \eqref{eq:teo1} in general is not sharp. 
Indeed, consider the map $u_V(x):=\frac{x}{|x|}$ 
defined on $B_r(0)\setminus \{0\}$ for $r>0$, whose distributional
Jacobian determinant is $\Det(\nabla u_v)=\pi\delta_0$. Theorem \ref{teo:upper_bound_W11Suno} implies 
that 
\begin{align}\label{eq:vortice}
\rilarea(u_V,B_r(0))\leq \int_{B_r(0)}\sqrt{1+|\nabla u|^2}dx+\min\{2\pi r,\pi\}.
\end{align}
On the other hand, according to \cite[Theorem 1.1]{BES}, one has
\begin{align}
\rilarea(u_V,B_r(0))=\int_{B_r(0)}\sqrt{1+|\nabla u|^2}dx+F(r),
\end{align}
where the singular contribution $F(r)\in (0,\pi]$ has
 the meaning of the area of a minimal surface solving 
a suitable non-parametric Plateau problem with partial free boundary. Specifically, $F(r)$ coincides with half of the area of a sort of 
catenoid $S\subset\R^3=\R^2\times \R$ with boundary $(\Suno\times\{0\})\cup(\Suno\times \{2r\})
$ and constrained to contain the segment $\{0\}\times [0,2r]$. In particular it can be seen that there exists a number $\overline r\in (0,\frac12)$ such that for $r\geq\overline r$ this catenoid reduces to 
two disks,
 and then $F(r)=\pi$, whereas for $r<\overline r$ there exists a non-trivial catenoid whose area is strictly smaller than the lateral area of the 
solid portion of the (smallest) cylinder containing it, namely $F(r)< 2\pi r$. This shows that for $r\geq \overline r$ estimate in \eqref{eq:vortice} is an equality, and that for $r<\overline r$ is not sharp. 
We emphasize that a more precise estimate than 
\eqref{eq:teo1}, and hopefully the sharp value of the left-hand side, 
seems quite difficult to obtain. On the one hand 
we expect that, when the singularities of a map $u$ are far
from each other, \eqref{eq:teo1} becomes sharp\footnote{For
instance, under the further assumption that the $1$-current
$S_{\min}$ given by Lemma \ref{lem:structure_of_minimizers_of_the_combinatorial_problem} vanishes.}. However, in the opposite case, a characterization as in \eqref{eq:vortice} needs some strong improvements of the techniques used in \cite{BES}. 
Indeed, in \cite{BES} the 
rotational invariance of the domain and of the map $u_V$ itself are strongly 
exploited to prove the lower bound, which is based on a cilindrical Steiner-type symmetrization for integral currents.
A similar technique has been employed in \cite{S} (see also \cite{BP}), 
yielding the  value of $\rilarea(u_T,B_r(0))$, where $u_T$ is the symmetric
triple junction map, a piecewise constant map taking three
 values in $\Suno$, each value on a $120^o$ sector. 
Also in that case the symmetries of the source and the target spaces allow to use such techniques. 
Without these symmetries, at the moment
little can be said about the exact expression of $\rilarea(\cdot,\Om)$. 
   
So, the nonlocality of the
$L^1$-relaxed functional
 seems not removable.
Thus, 
following De Giorgi \cite{DG}, it seems interesting
to consider a further ``relaxation'', this time looking at
the functional $\rilarea(u, \cdot)$, i.e., looking at it  
as a function of the open set: For every $V\subseteq\Om$, we set
\begin{equation}\label{eq:subadditive_interior_envelope}
	\rilrilarea(u,V):=\inf\left\{
\sum_{k=1}^{+\infty}
\rilarea(u,A_k):A_k\subseteq\Om\text{ open }, \bigcup_{k=1}^{+\infty}
 A_k\supseteq V\right\}.
\end{equation}
Actually, 
notice that, for all $u\in L^1(\Om;\R^2)$, $\rilrilarea(u,\cdot)$ is the trace of a regular Borel measure restricted to open sets (in $\Om$).

The estimate provided by Theorem \ref{teo:upper_bound_W11Suno} allows us to prove (Proposition
\ref{prop_2rel}) that 
$$\rilrilarea(u,U)=\int_U\sqrt{1+|\nabla u|^2}dx,\qquad \forall
u\in W^{1,1}(\Om;\Suno), 
\ 
\text{ for every open set }U\subseteq\Om.
$$
Using this, we are able  to show
our next
main result (Corollary \ref{cor:De_Giorgi}):
\begin{theorem}\label{teo:doppio_rilassato_introduction}
Let $\Om\subset\R^2$ be a bounded open set with Lipschitz boundary and
$u\in W^{1,1}(\Om;\Suno)$. Then 
\begin{equation*}
\rilrilarea(u,\Om)=
\inf
\Big\{
\rilarea(u,\Om\setminus C): 
C\subset\Om,
\mathcal H^0(C)< +\infty
\Big\}.
\end{equation*}
\end{theorem}
This theorem positively answers to an adaptation of a De Giorgi conjecture 
\cite[Conjecture 3]{DG}, 
provided one restricts the analysis to the space $W^{1,1}(\Om;\Suno)$.

Before concluding the introduction, it is worth recalling that
in several works (see \cite{BMP,BMP2,BBM} and references therein; a general survey can also be found in  \cite{BMbook}) the authors studied the analogue of our relaxation 
problem, with the area functional 
replaced by the total variation, for $W^{1,1}$ maps defined 
on a closed simply connected surface taking values in $\Suno$. They
were able to characterize the corresponding relaxed functional,
and showed that the singular contribution is given by 
$$
L(\Lambda)
:=
\sup_{\substack{\varphi\in 
{\rm Lip}_0(\Om), 
{{\rm lip}(\varphi)}\leq 1}}\langle \Lambda,\varphi\rangle
= 
\inf
\{
|S|_\Om:S\in\mathcal {\mathcal D}_1(\Om), \findistr=\partial S
\},
$$
which has the geometric meaning of the (geodesic)  length 
of a minimal connection between the poles of $\Lambda$.
The case considered in the present paper seems much more involved,
due to the presence of the minimal surfaces briefly discussed above.

The plan of the paper is the following: In Section \ref{sec:notation_and_preliminaries} we fix the setting and notation needed in the sequel. In Section \ref{sec:a_minimization_problem_for_atomic_distributions}
we investigate the minimization problem dual to \eqref{eq:def_M} (see \eqref{Mdef_loc} below), and we prove some regularity result for the minimizing currents (see also Remark \ref{finalremark} in the Appendix). In Section \ref{subsec:distributional_Jacobian_determinant} we collect some results on the distributional Jacobian for Sobolev maps taking
values in the circle. We briefly add some details to extend the well-known results for simply connected domains \cite{BMP,BMP2} to the case of non-simply connected domains, for the reader convenience. Notice that however many of these results were stated in the aforementioned references and also summarized in \cite{BMbook}. In Section \ref{sec:density_results_in_W11} we 
prove a density result for circle valued Sobolev maps, see Proposition \ref{prop:density_in_flat_norm} which needs some preparatory lemmata. Finally in Section \ref{sec:main_results} we prove Theorem \ref{teo:upper_bound_W11Suno}, whereas in the last section we investigate the countably subadditive interior envelope of the relaxed area functional and we prove Theorem \ref{teo:doppio_rilassato_introduction}. 

To conclude, we mention that
it would be interesting to extend Theorem 
\ref{teo:upper_bound_W11Suno} to maps $u \in BV(\Om; \Suno)$. 
 We leave this effort for future investigations; we only mention that in \cite{SS} some estimates can be given for specific piecewise constant maps.

\section{Notation and preliminaries}
\label{sec:notation_and_preliminaries}
In what follows,
$\Omega\subset \R^2$ 
is a fixed connected (but not necessarily
simply connected) bounded open set 
with Lipschitz boundary. 
We denote by $d(\cdot, \dOm)$ the distance from $\dOm$, 
and following \cite[pag. 96]{BMP}, 
by $\dovOm:\overline\Om\times\overline\Om\rightarrow [0,+\infty)$ the function
$$
\dovOm(x,y):=\min\big\{|x-y|,
d(x,\partial\Om)+d(y,\partial\Om)\big\}.
$$ 
Hence, if $\dovOm(x,y)=|x-y|$, then the closed segment $\overline{xy}$ 
joining $x$ and $y$ is contained in $\overline \Om$.

Given a vector $V=(V_1,V_2)\in \R^2$, we set 
$V^\perp:=(-V_2,V_1)$ its $\pi/2$-counterclockwise rotation.
If $V = \grad u$, then $\grad^\perp u$ stands for $(\grad u)^\perp$.
The distributional divergence of a vector field $V = (V_1, V_2)\in L^{1}(\Om;\R^2)$ 
is the distribution
\begin{equation}\label{eq:distributional_divergence}
	\langle\Div V,\varphi\rangle :=-\int_\Om 
V\cdot \nabla\varphi dx
\qquad \forall \varphi\in \testfunctions(\Om).
\end{equation}
If $V$ is sufficiently smooth,
${\rm Div} V$ equals
the pointwise divergence $\div V=
\frac{\partial V_1}{\partial x_1}+\frac{\partial V_2}{\partial x_2}$.

The distributional curl of $V \in L^1(\Om; \R^2)$ is 
the distribution 
\begin{equation}\label{eq:distributional_curl}
\langle\curl V,\varphi\rangle :=\int_\Om V\cdot \nabla^\perp \varphi dx
\qquad \forall \varphi\in 
\testfunctions(\Om),
\end{equation}
where $\nabla^\perp \varphi=
(-\frac{\partial\varphi}{\partial x_2},\frac{\partial \varphi}{\partial x_1})$.
If $V$ is sufficiently smooth then 
$\curl V=\div(V^\perp)=-\frac{\partial V_2}{\partial x_1}+\frac{\partial V_1}{\partial x_2}$.

The symbols $BV(A)$ (resp. $SBV(A)$) denotes the space of functions of bounded
variation (resp. special functions of bounded variation) in the open set 
$A \subseteq \R^2$; 
if $u\in BV(A)$, $\grad u$ stands for the absolutely 
continuous part of the gradient measure $Du$. Further,  $\jump{u}$ stands
for the difference $u^+ - u^-$ of the two traces
of $u$ on its jump set $J_u$, provided 
a unit normal vector field
to $J_u$ is assigned. We denote by
$BV(A; \R^2)$ the space of functions
of bounded variation in $A$ taking values in $\R^2$; if $u \in BV(A; \R^2)$,
$\vert \grad u\vert$ stands for the Frobenius norm of $\grad u$;
see \cite{AmFuPa:00}.

\begin{Definition}[\textbf{Dipole map}]\label{def:dipole_map}
Let $p,n\in \R^2$ be distinct, 
 and consider two polar coordinate systems $(\rho_p,\theta_p)$ and $(\rho_n,\theta_n)$ centered at $p$ and $n$, 
respectively,
chosen\footnote{The orientation of these
systems is always counterclockwise.}
so that both 
$\theta_p$ and $\theta_n$ have a jump of size $2\pi$ on 
$\ell_n\subset \ell$, where $\ell$ is the line containing $\overline{pn}$ and $\ell_n$ is the halfline with endpoint $n$ and not containing $p$. 
We let $w_{p,n}\in 
BV_{{\rm loc}}(\R^2)$ be the 
dipole map, defined as
\begin{align}\label{dip_map}
w_{p,n}:=\theta_p-\theta_n.
\end{align} 
\end{Definition}

Thus $w_{p,n}$ does not jump on $\ell_n$, while it jumps (of $2\pi)$ 
on the relative interior of $\overline{pn}$.
Notice that 
\begin{align}\label{gradw}
|\nabla w_{p,n}(x)|\leq |\nabla \theta_p(x)|+|\nabla\theta_n(x)|\leq C\Big(\frac{1}{\vert x-
p\vert}+\frac{1}{\vert x-n\vert}\Big)
\qquad \forall x\in \R^2 \setminus \ell.
\end{align}
For any open set $A \subset \R^2$, 
 $\mathcal D_0(A)$ and 
$\mathcal D_1(A)$ denote the zero-dimensional and $1$-dimensional
currents in $A$, respectively. The symbol $\vert \cdot\vert_A$ stands for
the mass of a current $\Lambda$ in $A$, while ${\rm supp}$ denotes the 
support of $\Lambda$ \cite{Federer}.

\subsection{Lipschitz maps; the flat norm}\label{subsec:Lipschitz_maps_and_flat_norms}
For any bounded open set $A \subseteq \R^2$,  we let 
${\rm Lip}_0(A)$ 
be the space of Lipschitz functions on $A$
vanishing on $\partial A$, endowed 
with the norm
\begin{align}\label{lip_norm}
	\|\varphi\|_{{\rm Lip}_0(A)}:=
\max
	\Big\{\|\varphi\|_{L^\infty(A)},
{\rm lip}(\varphi, A) 
\Big\},
\end{align}
where
\begin{equation}\label{eq:lip}
{\rm lip}(\varphi, A) 
:= \sup_{\substack{x,y\in A\\x\neq y}}\frac{|\varphi(x)-\varphi(y)|}{|x-y|}.
\end{equation}
In the Banach space ${{\rm Lip}}_0(A)$,
 the norms ${\rm lip}(\cdot, A)$ and
	$\|\cdot\|_{{\rm Lip}_0(A)}$
are equivalent.

In what follows it is also convenient to introduce
the equivalent norm
\begin{equation}\label{eq:Lip_alpha}
\|\varphi\|_{{\rm Lip}_0(A),\alpha}=\max
\Big\{
\|\varphi\|_{L^\infty(A)},\alpha {\rm lip}(\varphi,A)
\Big\}
\qquad\forall \varphi \in {\rm Lip}_0(A),
\end{equation}
with $\alpha$ as in \eqref{eq:alpha}.
For all these norms we 
drop the symbol $A$ in the sequel, when $A=\Om$.

We denote by 
${\rm Lip}_0(A)'$ 
the dual space of 
${\rm Lip}_0(A)$ (endowed with one of these norms). 
The (equivalent to each other) 
dual norms to \eqref{lip_norm}, \eqref{eq:Lip_alpha} on 
${\rm Lip}_0(A)'$ 
are respectively: 
\begin{equation}\label{Mdef_loc}
\begin{aligned}
\|\Lambda\|_{{\rm flat},A}
:=& \sup_{\substack{\varphi\in {\rm Lip}_0(A)
\\ \|\varphi\|_{{\rm Lip}_0(A)}\leq1}}\langle \Lambda,\varphi\rangle, \qquad\qquad
\|\Lambda\|_{{\rm flat},\alpha,A}
:= &
\sup_{\substack{\varphi\in {\rm Lip}_0(A)
\\
			\| \varphi\|_{{\rm Lip}_0(A),\alpha}\leq 1}}\langle \Lambda,\varphi\rangle,
\end{aligned}
\end{equation}
for all $\Lambda \in
		{\rm Lip}_0(A)'$,
see \eqref{eq:def_M}. 
Again, for these dual norms we 
usually drop the symbol $A$ when $A=\Om$.
The reason of the notation $\Vert \cdot\Vert_{{\rm flat}}$ is explained
by formula \eqref{flat=flat} below.

\subsection{The classes $X(\Om)$ and $X_f(\Om)$}
\label{subsec:the_classes}
Let
$((x_i,y_i))_{i \in \NN}\subset\overline\Om\times\overline\Om$ be a sequence
of 
pairs of points for which $\sum_{i=1}^{+\infty} \dovOm(x_i,y_i)<+\infty$.
We shall always
 suppose that $x_i \neq y_i$, while we do not exclude that $x_i = x_j$
and/or $y_h = y_k$ for some $i\neq j$, $h \neq k$. Namely,
$((x_i,y_i))_{i \in \NN}\subset\overline\Om\times\overline\Om \setminus 
\diagonal$, where $\diagonal$ is the diagonal of 
$\overline \Om \times \overline \Om$.

The 
measures 
$$
\Lambda_n:=\sum_{i=1}^n(\delta_{x_i}-\delta_{y_i}), \qquad
n \in \mathbb N,
$$
converge in ${{\rm Lip}}_0(\Om)'$ to 
$\sum_{i=1}^{+\infty}
 (\delta_{x_i}-\delta_{y_i})$.
Indeed, for any $\varphi\in \text{{\rm Lip}}_0(\Om)$ with $\|\varphi\|_{{\rm Lip}}\leq 1$ and any $n \in \mathbb N$, setting
$I_n := \{i \geq n+1: 
\dovOm
(x_i,y_i)=
\vert x_i - y_i\vert\}$, 
$B_n := \{i \geq n+1 : \dovOm
(x_i,y_i) = d(x_i, \dOm) + 
d(y_i, \dOm)\}$, we have
$$
\begin{aligned}
& \left\vert
\langle \sum_{i=n+1}^{+\infty} \delta_{x_i} -
\delta_{y_i},
\varphi \rangle
\right\vert 
= 
\left|
\sum_{i=n+1}^{+\infty} 
(\varphi(x_i)-\varphi(y_i))\right|
\leq
\left|
\sum_{i\in I_n}
(\varphi(x_i)-\varphi(y_i))\right|
+
\sum_{i\in B_n}
(\vert \varphi(x_i)\vert + \vert \varphi(y_i)\vert)
\\
\leq & 
\sum_{i\in I_n}
\vert x_i- y_i\vert + 
\sum_{i\in B_n} (d(x_i, \dOm) + d(y_i, \dOm))
= \sum_{i=n+1}^{+\infty} \dovOm(x_i,y_i)\rightarrow 0
\end{aligned}
$$ 
as $n\rightarrow +\infty$, 
where we use $\Vert \grad \varphi \Vert_\infty\leq 1$ and 
$\varphi=0$ on $\dOm$.

\begin{Remark}[\textbf{Non uniqueness of the representation}]\rm
The representation $\Lambda=\sum_{i=1}^{+\infty} 
(\delta_{x_i}-\delta_{y_i})$ is not unique; clearly two sequences
$((x_i,y_i))_{i\in \mathbb N} 
\subset \overline \Om \times \overline \Om
\setminus \diagonal$  and 
$((\widehat x_i,\widehat y_i))_{i\in \mathbb N} \subset \overline 
\Om \times \overline \Om \setminus \diagonal$
with 
$\sum_{i\in \mathbb N}\dovOm(x_i, y_i) < +\infty$, 
$\sum_{i\in \mathbb N}\dovOm(\widehat x_i, \widehat y_i) < +\infty$, 
define the same linear functional on $ \text{{\rm Lip}}_0(\Om)$ if 
\begin{align}
\langle 
\sum_{i=1}^{+\infty} (\delta_{x_i}-\delta_{y_i}),
\varphi \rangle
=\langle
\sum_{i=1}^{+\infty} (\delta_{\widehat x_i}-\delta_{\widehat y_i}),
\varphi \rangle \qquad \forall 
\varphi \in  {\rm Lip}_0(\Om).
\end{align}
\end{Remark}
We emphasize that the hypothesis 
$((x_i,y_i))_{i \in \NN}\subset\overline\Om\times\overline\Om\setminus \diagonal$ 
(instead that
$((x_i,y_i))_{i \in \NN}\subset \Om\times \Om \setminus \diagonal$)
is done for convenience, and it
may happen
that for some $i\in \mathbb N$,
either $x_i\in \partial \Om$ or $y_i\in \partial\Om$ (or both). 
Of course, if $x_i\in \partial \Om$ then $\delta_{x_i}=0$ in $\text{\rm Lip}_0(\Om)'$;
 the presence of $x_i$ affects the representation of
$\Lambda$, but not its action on ${\rm Lip}_0(\Om)$. 
Nevertheless, we can always assume that for all $i\in \mathbb N$, at least one among $x_i$ and $y_i$ belongs to $\Om$. 
To indicate such a property, we briefly write 
$$
(x_i, y_i) \in
\overline \Om \times \overline \Om \setminus \diagonals.
$$
Preferred
representations will be discussed at
the beginning of Section \ref{sec:a_minimization_problem_for_atomic_distributions}.

\begin{Definition}\label{def:X_f_Omega}
We set\footnote{We take the union with $\{0\}$
since later $\Lambda$ will be typically the Jacobian
determinant of a suitable map, and we want to include
the case in which the map is constant.}
\begin{equation}\label{X}
\begin{aligned}
	X(\Omega):=&
\Big\{\Lambda\in \text{{\rm Lip}}_0(\Om)':
\exists ((x_i,y_i))_{i \in \mathbb N} \subset \overline\Om\times\overline\Om
\setminus \diagonals,
\\
& \ \ \ \sum_{i=1}^{+\infty} 
\dovOm
(x_i,y_i)<+\infty, 
 \Lambda=
\sum_{i=1}^{+\infty}
 (\delta_{x_i}-\delta_{y_i})
\Big\} 
\cup \{0\}.
\end{aligned}
\end{equation}
\end{Definition}

We have seen
that

\begin{align}\label{serie_Lambda}
 \forall \Lambda\in X(\Om)
\qquad	 \ \ \langle \Lambda,\varphi\rangle=\sum_{i=1}^{+\infty} (\varphi(x_i)-\varphi(y_i)) 
\qquad\forall \varphi\in \text{{\rm Lip}}_0(\Om),
\end{align}
the series in \eqref{serie_Lambda} being convergent. 

\begin{Definition} We set
\begin{equation}\label{Xf}
\begin{aligned}
	X_f(\Omega):= & 
\Big\{
\findistr
\in \text{{\rm Lip}}_0(\Om)':
\exists m\in \mathbb N,\;  \exists (x_i,y_i)
\in \overline\Om\times\overline\Om \setminus \diagonals {\rm ~for~} i=1,\dots,m,
\\
& \ \ \findistr
=\sum_{i=1}^m( \delta_{x_i}-\delta_{y_i})
\Big\}\cup \{0\}.
\end{aligned}
\end{equation}
\end{Definition}
Every $\findistr\in X_f(\Om)$ is a Radon measure and can be 
identified with an integral $0$-current in $\mathcal D_0(\Om)$.
\begin{Remark}\label{sup_achieved}\rm 
If
 $\Lambda\in X(\Om)$  then, adapting the arguments in \cite[Proposition 18]{Ponce},  
it easily follows that the suprema in \eqref{Mdef_loc} are achieved (taking into account that we have Lipschitz maps which are null on $\partial\Om$).
\end{Remark}

\subsection{The classes 
 $\zerocurrfin$ and $\mathcal S$}
In the sequel we need to 
 consider the following classes of rectifiable currents in $\R^2$:
\begin{equation}\label{currents}
\begin{aligned}
&\zerocurrfin:=\Big\{R\in \mathcal D_0(\R^2):R=\sum_{i=1}^n
\sigma_i\delta_{z_i}\text{ for some }n\geq0,\;z_i\in \R^2,\sigma_i\in\{-1,+1\}\Big\},
\\
	&{\mathcal S}:=
	\Big\{S\in \mathcal D_1(\R^2):S=\sum_{k=1}^{+\infty}\jump{\overline{x_ky_k}} \text{ for some sequence }((x_k,y_k))_k\subset \R^2, \;
\sum_{k=1}^{+\infty} \vert y_k
- x_k\vert < +\infty
\Big\},
\end{aligned}
\end{equation}
and denote by $\zerocurrfin(A)$
and $\mathcal S(A)$ the classes in \eqref{currents} when the currents are 
restricted to an open set $A\subset\R^2$. 

By \cite[page 367]{Federer} and
Lemma \ref{fed} in the appendix
\begin{align}\label{flat=flat}
\|\Lambda\|_{\text{{\rm flat}}}=
\inf\{|R|_\Om+|S|_\Om: (R, S) \in \mathcal D_0(\Om) \times \mathcal D_1(\Om),\;\Lambda=R+\partial S\}
\qquad
\forall \Lambda\in X(\Om),
\end{align}
and similarly
\begin{align}\label{def_MT}
\|\Lambda\|_{\text{{\rm flat}},\alpha}
=\inf
\Big\{|R|_\Om+ \alpha^{-1}|S|_\Om:
(R,S)\in \mathcal D_0(\Om) \times  \mathcal D_1(\Om),\;\Lambda=
R+\partial S\Big\}
\qquad
\forall \Lambda\in X(\Om),
\end{align}
where we recall that $\alpha$ is defined in \eqref{eq:alpha}.
We shall prove 
that 
the infimum in \eqref{def_MT} is
attained and that, 
if $\Lambda \in X_f(\Om)$,
 minimizers $R_{{\rm min}}$ and $S_{{\rm min}}$ satisfy $R_{{\rm min}}\in \zerocurrfin$,  $S_{{\rm min}}\in {\mathcal S}$
(similar properties hold for \eqref{flat=flat}).

\section{A minimization problem for atomic distributions}
\label{sec:a_minimization_problem_for_atomic_distributions}
Our aim in this section is to show that, for all $\Lambda\in X(\Om)$, the infimum
on the right hand side of 
\eqref{def_MT} is a minimum, and to analyze the regularity of its 
minimizers 
(Proposition \ref{prop_inf_attained}); this will be done supposing first
that, in place of $\Lambda$, we consider $\findistr \in X_f(\Om)$.

\subsection{Properties $({\rm P})$ and (${\rm P_f}$)}\label{subsec:properties_P_and_P'}
 Given a distribution $\Lambda=\sum_{i=1}^{+\infty}
(\delta_{x_i}-\delta_{y_i})\in X(\Om)$, we can modify 
the set of points $x_i, y_i$ 
in the following way.
Take $i \in \mathbb N$; 
\begin{itemize}
\item[]
if  $\dovOm(x_i,y_i)=|x_i-y_i|$ we introduce two  (coinciding) points $\widehat x_i=\widehat y_i$ at the center of the segment $\overline{x_iy_i}$;
\item[] 
if $\dovOm(x_i,y_i)=d(x_i,\partial\Om)+d(y_i,\partial\Om)$ we 
choose two points $\widehat x_i, \widehat y_i\in\partial\Om$ so that 
$$d(x_i,\partial\Om)=|x_i-\widehat y_i|,\qquad d(y_i,\partial\Om)=|\widehat x_i-y_i|.$$
\end{itemize}
In this way
$$\sum_{i=1}^{+\infty} (|{\widehat x_i}-{y_i}|+|{x_i}-{\widehat y_i}|)=\sum_{i=1}^{+\infty} \dovOm({x_i},{y_i})<+\infty,
$$
and we can write 
$$\Lambda=\sum_{i=1}^{+\infty}(\delta_{x_i}-\delta_{y_i})=\sum_{i=1}^{+\infty}(\delta_{\widehat x_i}-\delta_{y_i})+\sum_{i=1}^{+\infty}(\delta_{x_i}-\delta_{\widehat y_i}) 
\qquad {\rm in}~ \mathcal D'(\Om).
$$
In particular, 
we may assume, after relabelling and renaming the points, that:
\begin{itemize}
	\item[$({\rm P})$] There are sequences $((x_i,y_i))\subset
\overline\Om \times \overline \Om \setminus \diagonal$ such that 
	\begin{align}\label{representation_T}
		\Lambda=\sum_{i=1}^{+\infty}(\delta_{x_i}-\delta_{y_i})\qquad\text{ and } \qquad \sum_{i=1}^{+\infty}|{x_i}-{y_i}|=\sum_{i=1}^{+\infty} \dovOm(x_i,y_i)<+\infty.
	\end{align}
\end{itemize}   

Using that $\Lambda=\sum_{i=1}^{+\infty}(\delta_{x_i}-\delta_{y_i})$ 
admits a representation as in \eqref{representation_T}, in \eqref{def_MT} we can 
choose
as a competitor the pair $(R,S)$, with
$R=0$,  $S=\sum_{i=1}^{+\infty}\jump{\overline{y_ix_i}}$, 
and we obtain
\begin{equation}\label{eq:flat_bounded_by_sum}
\Vert \Lambda\Vert_{{\rm flat},\alpha}
\leq 2
\sum_{i=1}^{+\infty}|{x_i}-{y_i}|.
\end{equation}
Recall that there can be repetitions
among the $x_i$'s, as well 
as among the $y_i$'s.
Now, let $\findistr
=\sum_{i=1}^n(\delta_{x_i}-\delta_{y_i})\in X_f(\Om)$.
After relabelling (and keeping the same symbols, for simplicity),
$\findistr$ 
admits the
representation
\begin{align}\label{T_finite}
	\findistr
=\sum_{k\in J^+}\delta_{x_k}-\sum_{k\in J^-}\delta_{y_k},\qquad x_k,y_k\in \Om, \
\quad x_k \neq y_k,
\end{align}
in $\mathcal D'(\Omega)$, where $J^+$ and $J^-$ are finite 
(possibly empty) subsets of $\mathbb N$ and 
$J^+\cap J^-=\varnothing$.
It is convenient to add some atoms to $\findistr$ 
as follows: for any $k\in J^+$ we consider a point 
$\hat y_k\in \partial\Om$ 
so that 
$|x_k-\hat y_k|=\dist(x_k,\partial\Om)$, and similarly 
 for any $k\in J^-$ 
we consider a point 
$\hat x_k\in \partial\Om$,
so that $|\hat x_k-y_k|=\dist(y_k,\partial\Om)$. 
In this way,
again without changing the notation and calling
once more $\hat y_k$ by $y_k$ and $\hat x_k$ by $x_k$
for simplicity, 
setting $I=J^+\cup J^-$, 
 we can always write $\findistr$ as 
\begin{align}
\findistr
=\sum_{k\in I}(\delta_{x_k}-\delta_{y_k}),\label{T_finite2}
\end{align}
with the following additional property:
\begin{itemize}
	\item[(${\rm P_f}$)] for every $k\in I$ one and only one 
among  $x_k$ and  $y_k$ belongs to $\Om$, 
$x_i\neq y_j$
for any $x_i,y_j\in\{x_k,y_h:x_k\in \Om,y_h\in \Om\}$, and 
$\vert x_k - y_k\vert = \dovOm(x_k, y_k)$.
\end{itemize}
This implies that in
\eqref{T_finite2} there are no cancellations in $\Om$.
Recall that there can be repetitions
among the $x_i$'s, as well 
as among the $y_i$'s.

\subsection{Analysis of the minimum problem \eqref{def_MT}}
\label{subsec:analysis_of_the_minimum_problem}
Let $T\in X_f(\Om)$ 
be represented as in \eqref{T_finite2} and satisfying (${\rm P_f}$). We consider a 
disjoint partition 
$\{I_P, I_D\}$ 
of $I$ (i.e., $I=I_P\cup I_D$, 
$I_P\cap I_D=\varnothing$, where we allow 
$I_P$ or $I_D$ to be empty) 
and, 
provided $I_D\neq\varnothing$, an injective map $\tau:I_D\rightarrow I$. Along with this, we define the currents 
\begin{equation}\label{eq:R_tau_S_tau}
\begin{cases}
R_\tau:= 
\displaystyle \sum_{k\in I_P}\delta_{x_k}-
\sum_{j\in I\setminus \tau(I_D)}\delta_{y_j},\quad
S_\tau:=\sum_{k\in I_D}\jump{\overline{y_{\tau(k)}x_k}} & {\rm if} ~ 
I_D\neq \varnothing,
\\
\\
R_\tau :=T, \quad  S_\tau:=0 
& 
{\rm if}~ I_D= \varnothing
\end{cases}
\end{equation}
(clearly
$R_\tau\in \zerocurrfin$ and $S_\tau$, being a finite sum, belongs to ${\mathcal S}$).
Namely, we split the set 
$I$ as the union of 
$\tau(I_D)$ and 
$I \setminus \tau(I_D)$; a point
labelled by an index $h=\tau(k) \in \tau(I_D)$ is coupled with
$x_k$, while a point labelled by an index  
$k \in I \setminus \tau(I_D)$ is uncoupled.

Notice that 
$$
\begin{aligned}
R_\tau+\partial S_\tau\qquad 
=& 
\sum_{k\in I_P}\delta_{x_k}
+
\sum_{k\in I_D}\delta_{x_k}
-\left(
\sum_{j\in I \setminus \tau(I_D)}
\delta_{y_j}
+ 
\sum_{k\in I_D}
\delta_{ y_{\tau(k)} }
\right)
\\
=&
\sum_{k\in I}\delta_{x_k}
-
\sum_{k\in I}\delta_{y_k}
=T 
\end{aligned}
\qquad  \qquad
\text{in }\mathcal D_0(\Om).
$$
\begin{lemma}\label{lemma32}
For any $T\in X_f(\Om)$
we have
\begin{equation}\label{eq:n_alpha_hatN_alpha}
\begin{aligned}
& \min\big\{
|R_\tau|_\Om+\alpha^{-1} |S_\tau|_\Om : 
(R_\tau, S_\tau) {\rm ~as ~in }~ \eqref{eq:R_tau_S_tau} 
\big\}
\\
=& 
\min\big\{
|R|_\Om+\alpha^{-1}|S|_\Om:(R,S) \in \zerocurrfin \times  {\mathcal S},
\findistr=R+\partial S \text{ in }\mathcal D_0(\R^2)
\big\},
\end{aligned}
\end{equation}
where on the left-hand side 
the minimum\footnote{
The existence of a minimizer is guaranteed since the number of competitors is finite.}
 is taken over all
disjoint  partitions $\{I_D, I_P\}$ of $I$ 
and all injective maps $\tau:I_D\rightarrow I$, as above. 
In particular, a minimizer of 
the left-hand side is also a minimizer of the right-hand side.
\end{lemma}
\begin{proof}
On the one hand, the inequality 
$\geq$  trivially holds in \eqref{eq:n_alpha_hatN_alpha}.
On the other hand 
also the converse inequality holds, since every 
competitor $(R,S)\in \zerocurrfin
\times  {\mathcal S}$ 
for the right-hand side, can be modified, not increasing its energy, 
into a competitor for the minimum
problem on the left-hand side. More specifically,
let $(R,S)\in \zerocurrfin
\times  {\mathcal S}$ be such that $R+\partial S=T$ in $\mathcal D_0(\R^2)$,
with $T$ represented as in \eqref{T_finite2} and satisfies ($P_f$); 
in particular $\partial S = T - R$ is a 
finite sum of Dirac deltas. By Federer decomposition theorem 
for $1$-currents \cite[4.2.25]{Federer}, we can write 
$$
S=\sum_{i=1}^{+\infty} S_i,\qquad \text{ in }\mathcal D_1(\R^2),
$$
with $S_i\in \mathcal S$ for all $i \in \mathbb N$, and either $\partial S_i=0$ 
(so $S_i$ is a loop) or $\partial S_i=\delta_{z_i}-\delta_{w_i}$
 for some $z_i\neq w_i$, $z_i, w_i\in \{x_k,y_k:k\in I\}$. 
If $\partial S_i =0$ we set $\widehat S_i:=0$, i.e.,
we remove the loop. If
 $\partial S_i=\delta_{z_i}-\delta_{w_i}$ and 
  $\supp(S_i)\cap (\R^2\setminus \Om)= \varnothing$, we 
set $\widehat S_i:=\jump{z_iw_i}$ (the segment $\overline{z_iw_i}$ is
not necessarily included in $\Om$).
If  $\partial S_i=\delta_{z_i}-\delta_{w_i}$ and
$\supp(S_i)\cap (\R^2\setminus \Om)\neq  \varnothing$ then, using $({\rm P_f})$, we set $\widehat S_i:=\jump{z_i\widehat z_i}+\jump{\widehat w_iw_i}$
where  $\widehat z_i\in \{y_k:k\in I\}$ 
is a point on $\partial \Om$ 
such that $d(z_i,\partial \Om)=|z_i-\widehat z_i|$ and 
similarly, $\widehat w_i\in \{x_k:k\in I\}$ 
is a point on $\partial \Om$  such that 
$d(w_i,\partial \Om)=
|w_i-\widehat w_i|$. Finally, if some $S_i=\jump{z_iw_i}$ is such that both $z_i$ and $w_i$ belong to $\partial\Om$, we remove $S_i$, whereas if only one of them belongs to $\partial \Om$, say $w_i\in \partial\Om$, we replace $S_i$ by $\widehat S_i:=\jump{\widehat z_i\widehat w_i}$ where, again, $\widehat w_i\in \{x_k:k\in I\}$ 
is a point on $\partial \Om$  such that 
$d(w_i,\partial \Om)=
|w_i-\widehat w_i|$. 
 
Then $|\widehat S_i|_\Om\leq|S_i|_\Om$ for all $i \in \mathbb N$, and moreover the 
support of  
$$\widehat S:=\sum_{i=1}^{+\infty}\widehat S_i$$
consists of finitely many segments (possibly with repetitions)
joining some point in $ \{x_k:k\in I\}$ to some point in $\{y_k:k\in I\}$. 
Furthermore, $\partial \widehat S=\partial S$. From this remark one can 
easily define two sets $I_P, I_D \subseteq I$ 
of indices 
and an injective map $\tau: I_D \to I$ so that 
$\widehat S=S_\tau$, $R=R_\tau$, and it is checked that 
$$|R_\tau|_\Om+|S_\tau|_\Om\leq |R|_\Om+|S|_\Om.$$
This concludes the proof.
\end{proof} 

\begin{Remark}\rm
As a consequence of the previous arguments, the minimum in the right-hand side can be taken among currents in supported on $\overline\Om$.
\end{Remark}

The following crucial fact is a result of regularity theory for 
minimal currents;
since we were not able to find a specific reference,
for the reader convenience we propose a direct proof,
independent of regularity theory.

\begin{Proposition}
\label{lem_min_finita}
	Let 
 $\findistr=\sum_{i=1}^\indicefinalesomma(\delta_{x_i}-\delta_{y_i})\in X_f(\Om)$. Then the infimum in \eqref{def_MT},
with $T$ in place of $\Lambda$,
is attained and there are minimizers 
$(R_{{\min}}, S_{{\min}}) \in 
\zerocurrfin \times 
 {\mathcal S}$.
\end{Proposition}
\begin{proof}
The minimum problem on the right-hand side of \eqref{eq:n_alpha_hatN_alpha}
is
attained, as a consequence of Lemma \ref{lemma32}, and is trivially
larger than or 
equal to 
$\Vert \findistr\Vert_{{\rm flat},\alpha}$, see 
\eqref{def_MT}. We claim that actually equality holds, which will imply the thesis.
To prove this, recalling \eqref{Mdef_loc},
it is sufficient to show that 
\begin{align}
\min\big\{
|R|_\Om+\alpha^{-1}|S|_\Om:(R,S) \in \zerocurrfin \times  {\mathcal S},
\findistr=R+\partial S
\big\}
\leq \sup_{\substack{\varphi\in {\rm Lip}_0(\Om)
\\
		\|\varphi\|_{{\rm Lip}_0,\alpha}\leq 1}}\langle T,\varphi\rangle,
\end{align}
and this readily follows from Proposition \ref{minimal_lip} in the Appendix.
\end{proof}

Now we prove that, for a general $\Lambda\in X(\Om)$, the	
infimum on the right-hand side of \eqref{def_MT} 
can be obtained infimizing just on pairs $(R,S)\in 
\zerocurrfin
\times\mathcal S$. 
\begin{Corollary}[\textbf{$\Vert \cdot\Vert_{{\rm flat},\alpha}$ 
as an infimum over $\zerocurrfin \times \mathcal S$}]\label{cor3.3}
We have 
\begin{align}\label{claim_int}
\Vert \Lambda\Vert_{{\rm flat}, \alpha}
=
\inf\big\{
|R|_\Om + 
\alpha^{-1} |S|_\Om: 
(R,S) \in \zerocurrfin \times  {\mathcal S},
\Lambda=
R+\partial S\text{ in }\mathcal D_0(\R^2)
\big\} \qquad
\forall \Lambda\in X(\Om).
\end{align}
\end{Corollary}

\begin{proof}
Given $\eps>0$, it is sufficient to show that there exist
$R_\eps\in \zerocurrfin$, 
$S_\eps\in {\mathcal S}$, such that $\Lambda=
R_\eps+\partial S_\eps$ and
$$
|R_\eps|_\Om + 
2
 |S_\eps|_\Om
\leq \Vert \Lambda\Vert_{{\rm flat}, \alpha}+\eps.
$$
Assuming $\Lambda=\sum_{i=1}^{+\infty}(\delta_{x_i}-\delta_{y_i})\in 
X(\Om)$ is a representation satisfying $({\rm P})$, 
select $N_\eps\in \mathbb N$ so that 
\begin{equation}\label{eq:sum_eps_3}
\sum_{i=N_\eps+1}^{+\infty}|x_i-y_i|<\frac\eps3.
\end{equation}
 Thus, for
$$
\Lambda_\eps:=\sum_{i=N_\eps+1}^{+\infty}(\delta_{x_i}-\delta_{y_i})
\in {{\rm Lip}}_0(\Om)', \qquad
T_\eps:=\sum_{i=1}^{N_\eps}(\delta_{x_i}-\delta_{y_i})
\in 	X_f(\Omega), 
$$
we have 
\begin{equation}\label{eq:eps_3_eps_3}
\Vert \Lambda_\eps\Vert_{{\rm flat}, \alpha}
\leq \frac\eps3, \qquad 
\Vert T_\eps\Vert_{{\rm flat}, \alpha}
\leq 
\Vert \Lambda\Vert_{{\rm flat}, \alpha}
+
\Vert \Lambda_\eps\Vert_{{\rm flat}, \alpha}
\leq 
\Vert \Lambda\Vert_{{\rm flat}, \alpha}
+\frac\eps3.
\end{equation}
 By Proposition 
\ref{lem_min_finita}, 
there are integral currents $\widehat R_\eps\in 
\zerocurrfin$ 
and $\widehat S_\eps\in {\mathcal S}$ with $T_\eps
=\widehat R_\eps+\partial 
\widehat S_\eps$ in $\mathcal D_0(\mathbb \R^2)$ 
such that 
\begin{equation}\label{eq:T_flat_R}
\Vert T_\eps\Vert_{{\rm flat}, \alpha}
=|\widehat R_\eps|_\Om+ 2|\widehat S_\eps|_\Om.
\end{equation}
Setting $R_\eps:=\widehat R_\eps$ and $S_\eps:=\widehat S_\eps+\sum_{i=N_\eps+1}^{+\infty}\jump{\overline{y_ix_i}}$ one 
sees that 
$S_\eps\in  {\mathcal S}$, $\Lambda=  R_\eps + \partial S_\eps$, and
using \eqref{eq:sum_eps_3}, \eqref{eq:T_flat_R},
and 
\eqref{eq:eps_3_eps_3},
$$|R_\eps|_\Om+2
|S_\eps|_\Om\leq |\widehat R_\eps|_\Om+ 
2|\widehat S_\eps|_\Om+\frac{2\eps}{3}=
\Vert T_\eps\Vert_{{\rm flat}, \alpha}
+\frac{2\eps}{3}\leq 
\Vert \Lambda\Vert_{{\rm flat}, \alpha}
+\eps.
$$
\end{proof}

\begin{Proposition}[\textbf{Existence of minimizers 
defining $\Vert \cdot\Vert_{{\rm flat},\alpha}$}]\label{prop_inf_attained}
	Let
 $\Lambda\in X(\Om)$. 
Then the infimum in \eqref{def_MT} is attained and there are minimizers $R_{{\min}}\in \mathcal D_0(\Om)$ and $S_{{\min}}\in \mathcal D_1(\Om)$ which are integer multiplicity currents.
\end{Proposition}
\begin{proof}
Represent
$\Lambda=\sum_{i=1}^{+\infty}
(\delta_{x_i}-\delta_{y_i})$ as in $({\rm P})$, with 
$\sum_{i=1}^{+\infty} \dovOm({x_i},{y_i})<+\infty$. 
By \eqref{claim_int}, we can find a sequence $((R
_\indice,S_\indice)) 
\subset \zerocurrfin\times \mathcal S$ 
(in particular, of integer 
multiplicity currents) with 
$\Lambda=R_\indice+\partial S_\indice$ 
in $\mathcal D_0(\Om)$ for any $\indice \in \mathbb N$,
and such that 
	$$\lim_{k\rightarrow +\infty}\left(|R_\indice|_\Om+
2 |S_\indice|_{\Om}\right)=
\Vert \Lambda\Vert_{{\rm flat},\alpha}.$$ 
	By compactness \cite[Theorem 7.5.2]{Krantz_Parks}, up to a (not relabelled)
subsequence, 
we know that $R_\indice\rightharpoonup R_{{\rm min}}\in \mathcal D_0(\Om)$ 
and $S_\indice\rightharpoonup S_{{\rm min}}\in \mathcal D_1(\Om)$ weakly as currents,
and  we have to prove that $R_{{\rm min}}$ and $S_{{\rm min}}$ can be chosen with integer multiplicity.
	Suppose $R_\indice=\sum_{i=1}^{m_\indice}\sigma_i\delta_{z_i}$ 
for some $m_\indice \in \mathbb N$,
with $z_i\in \Om$ and $\sigma_i\in \{-1,+1\}$; we 
may  assume that there are no cancellations in the previous expression. We introduce $w_i\in \partial \Om$ so that $|z_i-w_i|=d(z_i,\partial\Om)$, and write $R_\indice=\sum_{i=1}^{m_\indice}\sigma_i(\delta_{z_i}-\delta_{w_i})$. 
In this way 
$R_\indice=\partial \Sigma_\indice$ with 
 $\Sigma_\indice=\sum_{i=1}^{m_\indice}\sigma_i\jump{\overline{w_iz_i}}
\in{\mathcal S}$. 
	We have
$$|R_\indice|_\Om\leq 
\Vert \Lambda\Vert_{{\rm flat},\alpha}
+1$$
for $\indice$ large enough, and after passing to a not-relabelled subsequence, 
we have 
$R_\indice\rightharpoonup R_{{\rm min }}\in \zerocurrfin$ weakly in $\mathcal D_0(\Om)$. Moreover, since $|R_k|_\Om=m_k$, we find that $(m_k)$ is a bounded sequence, and thus the mass of $\Sigma_k$ satisfies
$$|\Sigma_k|_\Om\leq m_k\text{{\rm diam}}(\Om),$$
and 
is uniformly bounded in $k$.
Since $\partial \Sigma_\indice
=R_\indice$ in $\mathcal D_0(\Om)$, 
also $\Sigma_\indice\rightharpoonup \Sigma$ weakly in $\mathcal D_1(\Om)$, with $\Sigma$ an integral current.
	
	Now we know that $R_\indice+\partial S_\indice=\Lambda$ in $\mathcal D_0(\Om)$. Writing $\Lambda=\partial T$, with $T=\sum_{i=1}^{+\infty}\jump{\overline{x_iy_i}}$ an integer multiplicity current, we see that 
	$$\partial S_\indice=\partial T-\partial \Sigma_\indice\qquad \text{for $\indice$ large enough},$$
	and then $S_\indice+\Sigma_\indice+T\in \mathcal D_1(\Om)$ is an integral current without boundary. By compactness, we can assume that 
the sequence $(S_\indice + \Sigma_\indice + T)$ weakly 
converges in $\mathcal D_1(\Om)$ to an 
integral current $Q$ without boundary. On the other hand, 
since $S_\indice\rightharpoonup S_{{\rm min}}$ weakly in $\mathcal D_1(\Om)$, we 
conclude that $S_{{\rm min}}+\Sigma+T=Q$ is an integral current. In particular, $S_{{\rm min}}=Q-\Sigma-T$ is an integer multiplicity current.
\end{proof}

\subsection{Properties of minimizers}

Here we prove a useful lemma which summarizes some 
properties of the minimizing partition $\{I_P, I_D\}$ and of the minimizing 
map $\tau$ on the left-hand side of \eqref{eq:n_alpha_hatN_alpha}.

\begin{lemma}[\textbf{Structure of minimizers of the combinatorial problem}]
\label{lem:structure_of_minimizers_of_the_combinatorial_problem}
	Let $T\in X_f(\Om)$ 
	be of the form \eqref{T_finite2} satisfying (${\rm P_f}$).
	Then there exist a disjoint 
	partition $\{I_P, I_D\}$ 
	of $I$ and an injective map $\tau:I_D\rightarrow I$,
	minimizing 
the left-hand side of \eqref{eq:n_alpha_hatN_alpha} 
for which, 
	setting
	\begin{align}\label{defRSprimo_bis}
		&R_{{\rm min}}:=\sum_{k\in I_P}\delta_{x_k}-\sum_{j\in I\setminus \tau(I_D)}\delta_{y_j} \in
		\zerocurrfin,
		\qquad\qquad S_{{\rm min}}:=\sum_{k\in I_D}\jump{\overline{y_{\tau(k)}x_k}}
		\in \mathcal S,
	\end{align}
so that $T = R_{{\rm min}}+\partial S_{{\rm min}}$, 
the following properties hold:
	\begin{itemize}
		\item[(a)] 
		for all $k\in I_P$ and $j\in I\setminus \tau(I_D)$ 
		for which $x_k\in \Om$ and $y_j\in \Om$, we have
		\begin{align}\label{min_2}
			|x_k-y_j|\geq 1, \qquad
			d(x_k,\partial\Omega)\geq\frac12,
			\qquad d(y_j,\partial\Omega)\geq\frac12.
		\end{align}
		Moreover, if $k\in I_D$ is such that either $x_k\in \Om$ 
and $y_{\tau(k)}\in \partial\Om$, or $x_k\in \partial \Om$ and $y_{\tau(k)}\in \Om$, then $\tau(k)=k$; 
		\item[(b)] for all $k\in I_D$ the (relative) interior of the segment $\overline{y_{\tau(k)}x_k}$ is contained in $\Om$, and
		\begin{equation}\label{min_3}		
			\begin{aligned}
				|x_k-y_{\tau(k)}|\leq& 
\min\{1,d(x_k, \partial
\Om) + d(y_{\tau (k)}, \partial \Om)\},
				\\
				|x_k-y_{\tau(k)}|\leq& 
				\frac{1}{2}+ \min\{d(x_k,\partial\Om), d(y_{\tau(k)},\partial\Om)\};
			\end{aligned}
		\end{equation}
		\item[(c)] 
		if  $x_k\in \supp({S_{{\rm min}}}) \cap \Om$
for some $k\in I_P$, then $x_k=x_h$ for some $h\in I_D$;
		\item[(d)]
		if $y_j\in \supp({S_{{\rm min}}}) \cap \Om$ 
for some $j\in I\setminus\tau(I_D)$, 
then $y_j=y_{\tau(k)}$ 
for some $k\in I_D$;
		\item[(e)]  if $k,h\in  I_D$, 
$k \neq h$,
 and 
		$\overline{y_{\tau(k)}x_k}\cap \overline{y_{\tau(h)}x_h} = 
		\{r\}$, then either $r=y_{\tau(k)}=y_{\tau(h)}$ or $r=x_k=x_h$;
		
		\item[(f)] if 
$k,h\in  I_D$, $k \neq h$, and 
		$\overline{y_{\tau(k)}x_k}\cap \overline{y_{\tau(h)}x_h}$ contains more than one point, then either  $\overline{y_{\tau(k)}x_k}\cap \overline{y_{\tau(h)}x_h}=\overline{y_{\tau(k)}x_h}$ or $\overline{y_{\tau(k)}x_k}\cap \overline{y_{\tau(h)}x_h}=\overline{y_{\tau(h)}x_k}$;
		\item[(g)]
		if 
		the points in \eqref{defRSprimo_bis} contained in $\Om$ 
		are distinct and  three by three not collinear,
		then 
		the segments  $\overline{y_{\tau(k)}x_k}\cap \Om$, 
		$k\in  I_D$, are disjoint;
			\item[(h)] $|S_{{\rm min}}|_\Om=\sum_{k\in I_D}|x_k-y_{\tau(k)}|$;
in particular $\supp({S_{{\rm min}}})=\cup_{k\in  I_D}\overline{y_{\tau(k)}x_k}$.
	\end{itemize}
\end{lemma}
In words, (c) says that if $x_k$ intersects ${\rm supp}(S_{{\rm min}})$, 
the intersection happens
in one extremum of the intervals composing $S_{{\rm min}}$, and similarly for
$y_j$ in (d). 
(e) says that if two intervals of $S_{{\rm min}}$ intersect
at one point, this point must be an extremum of both. 
(f) says that if two intervals of $S_{{\rm min}}$ intersect
at more than one point, then they cannot be contained one inside the other.

\begin{proof}
	
	\item[(a)] Let us prove the first inequality\footnote{We
shall prove a stronger statement, namely the
validity of \eqref{min_2}
 for any minimizing $I_P$, $I_D$, $\tau$.}
 in \eqref{min_2}.
	Suppose, to the contrary, 
	that there exist $k\in I_P$ and $j\in I\setminus \tau(I_D)$ 
	such that $x_k\in \Om$ , $y_j\in \Om$ and
	$|x_k-y_j|<1$. 
	Define the injective map $\varphi:
	I_D\cup\{k\}
	\rightarrow I$ as follows: $\varphi =\tau$ on $I_D$,
	$\varphi(k):=j$. 
	Then 
$$
R_\varphi =
\sum_{h\in I_P \setminus \{k\}}
\delta_{x_h}
- \sum_{\iota \in I \setminus \varphi(I_D \cup \{k\})}
\delta_{y_\iota} = R_{\min} - \delta_{x_k}- \delta_{y_j},
$$

$$
S_\varphi =\sum_{h \in I_D \cup \{k\}} \jump{y_{\varphi(h)}
 x_h} = S_{\min} + \jump{y_j x_k}.
$$
Thus

	\begin{align*}
		&|R_{\varphi}|_\Om
		+ 2|S_{\varphi}|_\Om
\leq
		|R_{{\rm min}}|_\Om
		+ 2|S_{{\rm min}}|_\Om-|\delta_{x_k}|-|\delta_{y_j}|
		+2 |x_k-y_j|
		\\
		=& 
		|R_{{\rm min}}|_\Om
		+2|S_{{\rm min}}|_\Om
		+2(|x_k-y_j|-1)
		<
		|R_{{\rm min}}|_\Om
		+ 2|S_{{\rm min}}|_\Om,
	\end{align*}
	contradicting the minimality of $(R_{{\rm min}}, S_{{\rm min}})$. 
	
	Let us now show the last statement in (a).
Let $D^+:=\{k\in I_D:x_k\in \Om,\; y_{\tau (k)}\in \partial \Om\}$ 
and $D^-:=\{k\in I_D:x_k\in \partial \Om,\;y_{\tau(k)}\in \Om\}$. 
For all $k\in D^+\cup D^-$ define $\varphi(k):=k$, 
whereas $\varphi(k):=\tau(k)$ for all $k\in I_D\setminus(D^+\cup D^-)$. 
It is easily checked that $\varphi$ is injective,
and 
	$$\sum_{k\in I_D}|x_k-y_{\varphi(k)}|=
\sum_{k\in D^+\cup D^-}|x_k-y_{k}|+\sum_{k\in I_D\setminus (D^+\cup D^-)}|x_k-y_{\tau(k)}|\leq \sum_{k\in I_D}|x_k-y_{\tau(k)}|,$$
the inequality being true since, for $k\in D^+$ (and similarly for $D^-$), by ($P_f$), $y_k$ is a closest point on $\partial \Om$ to $x_k$. 
	In particular, replacing $\tau$ with $\varphi$ we get a minimizing configuration satisfying the last statement in (a). In words,
by assumption $x_k \in \Om$ implies
$y_{\tau(k)} \in \partial \Om$, and 
$d(x_k, \dOm)=  \vert x_k -
y_{\tau(k)}\vert$;
we have shown that there are minimizers for which
$d(x_k, \dOm)=  \vert x_k -
y_k\vert$, 
so we are ``connecting'' $x_k$ with $y_k$.

To conclude the proof of (a), we need to show the second
and third inequalities in \eqref{min_2}. 		
Let $k \in I_P$, and suppose by contradiction that 
	$d(x_k, \partial\Omega)= \vert x_k- y_k\vert <\dfrac12$. We then
extend $\tau$ on $I_D \cup \{k\}$ 
using $\varphi:=\tau$ on $I_D$ and  $\varphi(k):=k$. 
Notice that this extension is well-defined, since $y_k\in \partial\Om$ and 
the last statement of (a) is satisfied by $\tau$. Also in this case 
the new partition with $\varphi$ has smaller energy 
than the original one with $\tau$, since 
		\[
	1=|\delta_{x_k}|_\Om>2|x_k-y_{\varphi(k)}|, 
	\]
	and this is enough to prove that 
	\[
|R_{{\rm min}}|_\Om
+	2|S_{{\rm min}}|_\Om
>
|R_{\varphi}|_\Om +2|S_{\varphi}|_\Om,
	\]
	contradicting the minimality. 
 In a similar manner we prove the third inequality in  \eqref{min_2}.
\item[(b)]
	Let us start to prove that 
	$\overline{y_{\tau(k)}x_k}\subset\overline\Om$
	for all $k\in I_D$. 
	Suppose, to the contrary, there exists 
	$k\in I_D$ for which
	$\overline{y_{\tau(k)}x_k}\cap (\R^2 \setminus
	\overline\Om)\not=\emptyset$, so that
	$|\jump{y_{\tau(k)}x_k}|_\Om<
	|{y_{\tau(k)}-x_k}|$, and moreover it must hold $|\jump{y_{\tau(k)}x_k}|_\Om\geq \dovOm
(x_k,y_{\tau(k)}) $. 
In such a case $\tau(k)\neq k$, and so we set $\varphi(k):=k$ and $\varphi(j):=j$ for $j=\tau(k)$; moreover we set $\varphi:=\tau$ 
on the other indices. Owing to the last assertion in (a), $\varphi$ is well-defined, and since
	$$|x_k-y_k|+|x_j-y_j|\leq |x_k-y_j|,\qquad j=\tau(k),$$
	it easily follows that the new partition and $\varphi$ minimize \eqref{eq:n_alpha_hatN_alpha}. 
	This concludes 
the proof of the first assertion in (b). Let us prove the first expression in \eqref{min_3}: If $|x_k-y_{\tau(k)}|> \dovOm
(x_k,y_{\tau(k)})$ we modify the partition and $\tau$ as before, getting a contradiction with the minimality. If 
 $|x_k-y_{\tau(k)}|>1$, we erase $k$ from $I_D$, and we find out that the new partition with $\tau$ replaced by its restriction on $I_D\setminus \{k\}$ realizes a smaller contribution, contradicting the minimality.
	
	The last inequality in \eqref{min_3} follows 
from the following argument: We may assume without loss of generality that $x_k,y_{\tau(k)}\in \Om$ and that, by the first assertion in (b), the segment joining them has interior in $\Om$; by (a), we can also suppose that $j:=\tau(k)\notin I_D$. Hence 
we can delete $k$ from $I_D$ and add $j$ to it, defining $\varphi(j):=j$ 
and $\varphi:=\tau$ elsewhere. In such a case, by minimality assumption we 
obtain 
	$$1+d(y_j,\partial \Om)=|\delta_{x_k}|+2|\jump{y_{\varphi(j)}x_j}|_\Om\geq 2|\jump{y_{\tau(k)}x_k}|_\Om=2|x_k-y_{\tau(k)}|.$$  
	\item[(c)] 
	Suppose there exists 
	$k\in I_P$ 
with $x_k\in \Om$,  $x_k\in\supp({S_{{\rm min}}})\setminus\bigcup_{h\in I_D}x_h$; then necessarily 
$x_k$ belongs to the relative interior of some segment 
$\overline{y_{\tau(j)}x_j}$, with $j \in I_D$ and
$j\neq k$, so that 
$\vert x_j - y_{\tau(j)} \vert
= \vert x_j - x_k \vert+ 
\vert x_k - y_{\tau(j)} \vert$.
 Set $\widetilde{I}_D=I_D\cup\{k\}\setminus\{j\}$ and 
let $\varphi:\tilde{I}_D\rightarrow I$ be the injective map such that  
$\varphi(i):=\tau(i)$ if $i\neq k$, $\varphi(k):=\tau(j)$. Now 
	\[
	2|x_j-y_{\tau(j)}|+|\delta_{x_k}|=2|x_k-y_{\tau(j)}|+2|x_j-x_k|+|\delta_{x_j}|>2|x_k-y_{\varphi(k)}|+|\delta_{x_j}|, 
	\]
implying that $|x_j-y_{\tau(j)}|>|x_k-y_{\varphi(k)}|$. Since  $R_{{\rm min}}$ and $R_{\varphi}$ have the same mass in $\Om$, the previous inequality readily gives
	\[
|R_{{\rm min}}|_\Om
+	2|S_{{\rm min}}|_\Om
>|R_{\varphi}|_\Om
+2|S_{\varphi}|_\Om,
	\]
	contradicting the minimality of $(R_{{\rm min}}, S_{{\rm min}})$.
	In a similar manner we prove (d).

	\item[(e)] Suppose to the contrary that 
$r$ belongs to 
	$\overline{y_{\tau(k)}x_k}\setminus\{y_{\tau(k)}, x_k\}$.  
Set $\varphi: I_D\rightarrow I$, 
	$\varphi(j):=\tau(j)$ if $j\not=k, h$, $\varphi(k):=\tau(h), \varphi(h)=\tau(k)$. We have 
\begin{equation}\label{eq:summing}	
\begin{aligned}
		&|x_k-y_{\varphi(k)}|=|x_k-y_{\tau(h)}|\leq |r-y_{\tau(h)}|+|x_k-r|,
\\
		&|x_h-y_{\varphi(h)}|=|x_h-y_{\tau(k)}|\leq |r-y_{\tau(k)}|+|x_h-r|,
	\end{aligned}
\end{equation}
	where at least one of these inequalities holds strictly, 
because 
the points
$y_{\tau(k)}, x_k, y_{\tau(h)}, x_h$ are not collinear by construction. 
Summing the inequalities in \eqref{eq:summing}  we get
	\[
	|x_k-y_{\varphi(k)}|+|x_h-y_{\varphi(h)}|<|x_k-y_{\tau(k)}|+|x_h-y_{\tau(h)}|, 
	\]
	and this is enough to deduce that 
	\[
|R_{{\rm min}}|_\Om
+	2|S_{{\rm min}}|_\Om
>
|R_{\varphi}|_\Om
+2|S_{\varphi}|_\Om,
	\]
	contradicting once again the minimality.

\item[(f)] If 	$\overline{y_{\tau(k)}x_k}\cap \overline{y_{\tau(h)}x_h}$ contains more than one point, it must contain a segment. In particular we have to 
exclude the two cases:  	$\overline{y_{\tau(k)}x_k}\cap \overline{y_{\tau(h)}x_h}=\overline{x_kx_h}$ and $\overline{y_{\tau(k)}x_k}\cap \overline{y_{\tau(h)}x_h}=\overline{y_{\tau(k)}y_{\tau(h)}}$; let us discuss the former (the latter being similar). In such a case it is 
sufficient to set $\varphi(k)=h$, $\varphi(h)=k$, and $\varphi=\tau$ otherwise, and check that the new map $\varphi$ associated with the same partition provides
\[
|R_{{\rm min}}|_\Om
+2|S_{{\rm min}}|_\Om
>
|R_{\varphi}|_\Om+
2|S_{\varphi}|_\Om,
\] 
contradicting the hypothesis.

\item[(g)] Follows from (e) and (f).
\item[(h)] Follows from the last assertion in (a) and the first in (b).
\end{proof}

\section{Distributional Jacobian;  maps with values in $\mathbb S^1$}
\label{subsec:distributional_Jacobian_determinant}
If $u = (u_1, u_2)\in W^{1,1}(\Om;\R^2)\cap L^\infty(\Om;\R^2)$ 
its distributional Jacobian determinant is 
the distribution $\Det(\nabla u)\in \mathcal D'(\Om)$ defined by
\begin{equation}\label{eq:distributional_Jacobian}
\langle \Det(\nabla u),\varphi\rangle :=\int_{\Om}
\lambda_u
\cdot \nabla \varphi dx\qquad \forall \varphi\in 
\testfunctions(\Om),
\end{equation}
where
\begin{align}
	\lambda_u
:=\frac12(u_1\nabla^\perp u_2-u_2\nabla^\perp
u_1)=\frac12\Big(-u_1\frac{\partial u_2}{\partial x_2}+u_2\frac{\partial u_1}{\partial x_2},u_1\frac{\partial u_2}{\partial x_1}-u_2\frac{\partial u_1}{\partial x_1}\Big)\in L^1(\Om;\R^2),
\end{align}
hence (cfr. \eqref{eq:distributional_divergence}) 
$$
\Det(\nabla u)=-\Div\lambda_u\in \mathcal D'(\Om).
$$
Moreover, since $\lambda_u\in L^1(\Om;\R^2)$, 
equality \eqref{eq:distributional_Jacobian} extends to $\varphi\in \text{{\rm Lip}}_0(\Om)$, so that 
$$
\Det(\nabla u)\in \text{{\rm Lip}}_0(\Om)'.
$$
It follows from the definition that the distributional
Jacobian enjoys some well-known compactness properties. For instance, 
let $u\in  W^{1,1}(\Om;\R^2)\cap L^\infty(\Om;\R^2)$. Let 
$(v_k)\subset  W^{1,1}(\Om;\R^2)\cap L^\infty(\Om;\R^2)$ be a bounded
sequence 
in $L^\infty(\Om, \R^2)$, and suppose 
$v_k\rightarrow u$  in $W^{1,1}(\Om;\R^2)$. Then
	\begin{align}\label{conv_det}
	\Det(\nabla v_k)\rightharpoonup \Det(\nabla u)\qquad\text{ in }\mathcal D'(\Om).
\end{align}
	Furthermore 
there exists $C>0$ such that 
	\begin{align}\label{bound_det}
		\|\Det(\nabla u)\|_{{\rm flat}}\leq C\|u\|_{W^{1,1}}.
	\end{align} 
\begin{Remark}\label{rem:strenghtened_convergence}\rm 
The convergence in \eqref{conv_det} can be 
strengthened into
\begin{align}\label{flat_strong}
\Vert \Det(\nabla v_k)- \Det(\nabla u)\Vert_{\text{{\rm Lip}}_0(\Om)'} \to 0.
\end{align}
Indeed, take a subsequence $(k_h)$; for 
any $\varphi\in  \text{{\rm Lip}}_0(\Om)$ write
\begin{equation}
\label{eq:estimate_for_improving_the_convergence}
\begin{aligned}
& \langle\Det(\nabla u)-\Det(\nabla v_{k_h}),
\varphi\rangle
\\
=& \int_\Om(\lambda_u-\lambda_{v_{k_h}})\cdot 
\nabla\varphi dx \leq\|\nabla \varphi\|_{L^\infty}\int_\Om|
\lambda_u-\lambda_{v_{k_h}}|dx
\\
\leq & \|\nabla \varphi\|_{L^\infty}\left(C_1 
\int_\Om |\nabla u-\nabla v_{k_h}|dx+C_2 \int_\Om|\nabla u \cdot (u-v_{k_h})|dx\right).
\end{aligned}
\end{equation}
Since 
$(u-v_{k_h})$ 
tends to zero in $L^1(\Om; \R^2)$ 
and since we can select a further subsequence $(k_{h_l})$ such that 
$(u-v_{k_{h_l}})$ tends to zero
weakly-star in $L^\infty$,
 we deduce that the limit of the right-hand side 
of \eqref{eq:estimate_for_improving_the_convergence} 
vanishes along the sub-subsequence, as $l\rightarrow +\infty$. 
In particular, taking the 
supremum of the left-hand side over $\varphi\in \text{{\rm Lip}}_0(\Om)$ with $\|\varphi\|_{ \text{{\rm Lip}}_0}\leq 1$, we infer
\begin{align}\label{flat_strong_sup}
\Vert \Det(\nabla v_{k_{h_l}})- \Det(\nabla u)\Vert_{\text{{\rm Lip}}_0(\Om)'} \to 0.
\end{align}
Thus, \eqref{flat_strong} follows from the Uryshon property.
\end{Remark}

\subsection{Maps with values in $\Suno$}
We collect here some useful tools and results, mostly
 on Sobolev maps taking values in $\Suno$. 
A large literature on this topic is available, e.g., following the results by Brezis and coauthors (see for instance \cite{BMbook} and references therein). 
Together with the Jacobian determinant, it is useful to introduce 
the notions of degree and lifting.
\begin{Definition}[\textbf{Degree}]\label{def:degree}
Let $B_r \subset \R^2$ be a disc of radius $r>0$,
and $\nu$ the outer unit normal vector to $\partial B_r$. 
The degree of a map $u\in W^{1,1}(\partial B_r;\Suno)$ 
 is defined as
\begin{align}
\deg(u;\partial B_r)=\frac{1}{2\pi}\int_{\partial B_r}
\Big(u_1\frac{\partial u_2}{\partial\tau}-u_2\frac{\partial u_1}{\partial\tau}\Big)
~d\mathcal H^1,
\end{align}
where $\tau:=\nu^\perp$.
\end{Definition}
Notice that $\deg(u;\partial B_r)\in \mathbb Z$.

\begin{Definition}[\textbf{Lifting}]\label{def:lifting}
Let $u = (u_1,u_2)\in BV(\Om;\Suno)$. We say that
$w\in BV(\Om)$ is a lifting of $u$ 
if $e^{iw}=(\cos w,\sin w)=(u_1,u_2)$ a.e. in $\Om$.
\end{Definition}

The following result holds
\cite{GMS2}, 
\cite[Th\'eor\`eme 0.1, Remarque 0.1]{DaIg:03}:
\begin{theorem}\label{teo_lifting}
Let $u\in
 BV(\Om;\Suno)$. Then there exists a 
lifting 
$w\in BV(\Om)$ of $u$ 
such that 
\begin{equation}\label{eq:norm_of_the_lifting}
\|w\|_{BV} \leq 2 \|u\|_{BV}.
\end{equation}
If furthermore 
$u \in  SBV(\Om;\Suno)$, then $w\in SBV(\Om)$.
\end{theorem}
Liftings $w$ provided by 
Theorem \ref{teo_lifting} 
satisfy the following important feature: if $u \in W^{1,1}(\Om; \Suno)$, then
\begin{align*}
	\Det(\nabla u)
=\frac12\curl(\nabla w) \qquad \text{ in }\mathcal D'(\Om),
\end{align*}
see \eqref{eq:distributional_curl}.
Indeed, for any $\varphi\in \testfunctions(\Om)$,
\begin{equation}\label{det=curl}
\begin{aligned}
	\langle\Det(\nabla u),\varphi\rangle&=\frac12\int_\Om \left(u_1
\nabla^\perp u_2-u_2\nabla^\perp u_1\right)\cdot\nabla \varphi dx
\\
	&=\frac12\int_\Om\nabla^\perp w\cdot\nabla \varphi dx=\frac12\langle \curl(\nabla w),\varphi\rangle.
\end{aligned}
\end{equation}

Let $B=B_R(0) \supset \supset \Om$ 
be an open disc,
 for some $R>0$ big enough, and 
$u\in W^{1,1}(\Om;\Suno)$.
We claim that there exists an extension $\overline u\in W^{1,1}(B;\Suno)$ of $u$.
Indeed, 
let $w \in BV(\Om, \Suno)$ be a lifting of $u$; 
since  $\Om$ has Lipschitz boundary, by \cite[pag. 162]{Giu:84} 
there exists $\widehat w\in W^{1,1}(B\setminus \overline\Om) \cap
BV(B)$,
with trace $\widehat w\res\partial\Om=w\res\partial\Om$.
 If we set
\begin{align}\label{extension_scalar}
\overline{w}:=\begin{cases}
	w&\text{ in }\Om,
\\
	\widehat{w}&\text{ in }B\setminus \overline\Om,
\end{cases}
\end{align} 
the map
\begin{align}\label{extension}
\overline u:=e^{i\overline{w}}
\end{align} 
is the map we are looking for. 
It is easy to see
that  by construction $\deg(\overline u,\partial B_r(0))=0$
 for a.e. $r>0$ 
with $\Om\subset\subset B_r(0)\subset B$. 


In what follows, we will need the following standard density result:
\begin{theorem}[\textbf{Density of $\mathcal C^\infty$ in $W^{1,1}(A; \Suno)$}]\label{teo:density_of_Cinfty_in_W11}
If $A\subset\R^2$ is a 
connected simply connected domain with smooth boundary, then the 
class 
\begin{align*}
\Big\{v\in W^{1,1}(A;\Suno):
\exists n\in \mathbb N, ~ \exists \{a_1,\dots,a_n\}\subset A,\;v\in C^\infty(A 
\setminus \{a_1,\dots,a_n\};\Suno)\Big\}
\end{align*}
is dense in $W^{1,1}(A;\Suno)$. Furthermore
\begin{align}\label{detS}
	\Det(\nabla v)=\pi\sum_{i=1}^nd_i\delta_{a_i} \qquad
\forall v\in W^{1,1}(A;\Suno)\cap C^\infty(A\setminus \{a_1,\dots,a_n\};\Suno),
\end{align}
where $d_i=\deg(v;\partial B_r(a_i))$ for any $r>0$ small enough.
\end{theorem}
\begin{proof}
See  \cite[Theorem 4 with $k=1$]{BZ}, 
and \cite[Lemma 2]{BBM} for the second part of the statement.
\end{proof}

The next theorem is an extension of \cite[Theorems 3 and 3']{BMP} 
to 
non-simply connected domains in $\R^2$. 
Even if 
it  can be directly obtained from  \cite{BMP} 
and \cite{BBM}, for convenience we give a quick proof; 
for a more detailed 
discussion we refer to \cite[Chapter 14]{BMbook}.
\begin{theorem}[\textbf{Distributional Jacobian of $\mathbb  S^1$-valued 
maps}]\label{teo:distributional_Jacobian_of_S1_valued_maps}
Let $u\in W^{1,1}(\Om;\Suno)$. 
Then 
$$
\frac1\pi\Det(\nabla u)\in X(\Om),
$$
i.e., there exists a
sequence $((x_i,y_i)) \subset \overline\Om \times \overline\Om \setminus
\diagonals$ 
such that 
$\sum_{i=1}^{+\infty} \dovOm(x_i,y_i)<+\infty$ and
 $\Det(\nabla u)=\pi\sum_{i=1}^{+\infty}(\delta_{x_i}-\delta_{y_i})$.
\end{theorem}

\begin{proof}
We use an argument similar to \cite[Lemma 12']{BBM}.
Let $\overline u \in W^{1,1}(B; \Suno)$ be an extension of $u$ as in \eqref{extension}.
Using Theorem 
\ref{teo:density_of_Cinfty_in_W11} we can select a sequence 
$$
(u_k)_k\subset 
\Big\{v\in W^{1,1}(B;\Suno):
\exists n\in \mathbb N,\; \exists \{a_1,\dots,a_n\}\subset  B,\;v\in C^\infty(B 
\setminus \{a_1,\dots,a_n\};\Suno)\Big\}
$$
converging to $\overline u$ 
in $W^{1,1}(B;\R^2)$. For all $k>0$ we can assume 
that  $\deg(u_k;\partial B_r(0))=0$ for some $r>0$ 
big enough with $\Om\subset\subset B_r(0)\subset B$, and so we can also rewrite 
\eqref{detS} in the following way:
\begin{equation}\label{det_points}
\frac1\pi\Det(\nabla u_k)=\sum_{i=1}^{n_k}(\delta_{x^k_i}-\delta_{y^k_i}),
\end{equation}
for suitable (not necessarily distinct) 
points $x_i^k,y_i^k\in B_r(0)$ 
and $n_k\in\mathbb N$. 
Furthermore, owing to \eqref{flat_strong} 
we may suppose\footnote{We use that $\|\cdot\|_{{\rm flat},B}$ and $\|\cdot\|_{{\rm Lip}_0(B)'}$ are equivalent.} 
\begin{align}\label{det_piccolo}
	\|\Det(\nabla u_{k+1})-\Det(\nabla u_{k})\|_{{\rm flat},B}\leq \frac{1}{2^{k}} \qquad \forall k >0.
\end{align}
As a consequence, 
we can write
$\Det(\nabla \overline u)=\Det(\nabla u_1)+\sum_{k=1}^{+\infty}(\Det(\nabla u_{k+1})-\Det(\nabla u_{k}))$,
the series being absolutely 
convergent in $ \text{{\rm Lip}}_0(B)'$. 
Up to relabelling the indices in \eqref{det_points}, we assume that for $k>0$ 
\begin{align}\label{rep1}
&	\Det(\nabla u_1)=\pi\sum_{i=1}^{m_{1}}(\delta_{x_i}-\delta_{y_i})\qquad \text{ in }\text{{\rm Lip}}_0(B)',\\
&\Det(\nabla u_{k+1})-\Det(\nabla u_{k})=\pi\sum_{i=m_k+1}^{m_{k+1}}(\delta_{x_i}-\delta_{y_i})\qquad \text{ in }\text{{\rm Lip}}_0(B)', \label{rep2}\end{align}
in such a way that 
$$\Det(\nabla \overline u)=\pi\sum_{i=1}^{+\infty}(\delta_{x_i}-\delta_{y_i})\qquad \text{ in }\text{{\rm Lip}}_0(B)'.$$
At the same time, restricting to $\text{{\rm Lip}}_0(\Om)$
the above linear functionals defined on $\text{{\rm Lip}}_0(B)$, 
we can replace the preceeding representations 
and assume that the points $x_i$ and $y_i$ in \eqref{rep1} and \eqref{rep2} belong to $\overline\Om$, are of the form \eqref{T_finite2}, 
and enjoy property $({\rm P})$. Up to a permutation of the points $(y_i)$, one can 
further suppose that 
\begin{align*}
	\|\Det(\nabla u_1)\|_{{\rm flat}}=\sum_{i=1}^{m_1}\dovOm(x_i,y_i),\qquad
	\|\Det(\nabla u_{k+1})-\Det(\nabla u_{k})\|_{{\rm flat}}=\sum_{i=m_k+1}^{m_{k+1}}\dovOm(x_i,y_i).
\end{align*}
This, together with \eqref{det_piccolo}, implies
$\|\Det(\nabla u)\|_{{\rm flat}}\leq \sum_{i=1}^{+\infty} \dovOm(x_i,y_i)<+\infty$, 
which concludes the proof.
\end{proof}

\section{Density results in $W^{1,1}(\Om; \Suno)$}
\label{sec:density_results_in_W11}
In this section we want to 
show the following density result, which is an immediate consequence
of
Lemmas \ref{lem:density_II} and \ref{lem:density_III} below. 

\begin{Proposition}[\textbf{Density in $W^{1,1}(\Om; \Suno)$}]
\label{prop:density_in_flat_norm}
Let 
$u\in W^{1,1}(\Om;\Suno)$. Then for all $\eps>0$ 
there exists a map $u_\eps\in W^{1,1}(\Om;\Suno)$ with the following properties:
\begin{itemize}
	\item[(i)] $\Det(\nabla u_\eps)=\pi\sum_{i=1}^{\indicefinalesomma_\eps}
(\delta_{x_i}-\delta_{y_i})$ in ${\rm Lip}_0(\Om)'$ for some
$\indicefinalesomma_\eps \in \mathbb N$,
with 
distinct and 
three by three not collinear points $x_i$, $y_i$ in $\Om$;
\item[(ii)] there exist positive numbers $\rho_{x_i}<\eps$ and 
$\rho_{y_i}<\eps$, $i=1,\dots,\indicefinalesomma$, such that 
the discs of the family $\{B_{\rho_{x_i}}(x_i),B_{\rho_{y_i}}(y_i):
x_i\in\Om, y_i\in \Om\}$ are contained in $\Om$ and 
pairwise disjoint, and\footnote{Here and in the
sequel, 
$\theta_x$ is the polar angular coordinate around $x$.}
 $u_\eps=e^{i\theta_{x_i}}$ in $  B_{\rho_{x_i}}(x_i)$, 
and $u_\eps=e^{-i\theta_{y_i}}$ in $  B_{\rho_{y_i}}(y_i)$;
\item[(iii)] 
$\|u-u_\eps\|_{W^{1,1}}+\|\Det(\nabla u)-\Det(\nabla u_\eps)\|_{{\rm flat}}
< \eps$.
\end{itemize}
\end{Proposition}

Recalling the Definition \ref{def:dipole_map} of the dipole map $w_{p,n}$,
we set  $v_{p,n}:=e^{iw_{p,n}}\in W^{1,1}(\Om;\Suno)$, which satisfies 
\begin{align*}
\Det(\nabla v_{p,n})=\pi(\delta_p-\delta_n).
\end{align*}

\begin{lemma}[\textbf{Density: finite number of singular points}]\label{lem:density_I}
	Let $u\in W^{1,1}(\Om;\Suno)$ 
and write, by property $({\rm P})$, $\Det(\nabla u)=\pi \sum_{i=1}^{+\infty}(\delta_{x_i}-\delta_{y_i})$, with $x_i,y_i\in \overline\Om$, $x_i \neq y_i$,
 and $\sum_{i=1}^{+\infty}|x_i-y_i|<+\infty$. Then, for all $\eps>0$ there exists a map $u_\eps\in W^{1,1}(\Om;\Suno)$  such that:  
	\begin{itemize}
		\item[(i)] 
$\Det(\nabla u_\eps)=\pi\sum_{i=1}^{N_\eps}(\delta_{x_i}-\delta_{y_i})$
for some 
$N_\eps\in \mathbb N$;
		\item[(ii)]
		$\|u-u_\eps\|_{W^{1,1}}+\|\Det(\nabla u)-\Det(\nabla u_\eps)\|_{{\rm flat}}<\eps$. 
	\end{itemize}
\end{lemma}

\begin{proof}
Let $\eta>0$ and choose $N_\eta\in \mathbb N$ 
so that $\sum_{i=N_\eta+1}^{+\infty} |x_i-y_i|<\eta/2$. 
Given 
$(x_i,y_i)$ with $i>N_\eta$,  
consider the dipole map $w_i:=w_{x_i,y_i}\in BV_{{\rm loc}}(\R^2)$ 
 in \eqref{dip_map} 
and
the cut-off function
$\psi^\eta_i:\R^2\rightarrow\R$ given by 
	\begin{equation}\label{eq:eta_i}
		\psi^\eta_i(x)=\varrho\Big(\frac{1}{\eta_i}d(x,\overline{ x_i y_i})\Big), 
\qquad \eta_i:=2^{-i}\eta,
	\end{equation}
	where 
$\varrho\in C^\infty([0,1])$ is non-increasing, $\rho \equiv 1$ 
in a right neighborhood of $0$, $\rho \equiv 0$ 
in a left neighborhood of $1$, and 
$|\varrho'|\leq {2}$. The support of $\psi^\eta_i$ satisfies
\begin{equation}\label{eq:measure_spt_psi_i}
|{\rm spt}(\psi^\eta_i)|\leq \pi\eta_i^2+2\eta_i|x_i-y_i|
\qquad \forall i>N_\eta.
\end{equation}
	By \eqref{gradw} 
one checks\footnote{This estimate can be obtained integrating the right-hand side of  \eqref{gradw} in the two discs $B_{\eta_i}(x_i)$ and $B_{\eta_i}(y_i)$, and estimating $|\nabla w_i|$ by $C/\eta_i$ in the remaining part of ${\rm spt}(\psi^\eta_i)$.},
for the approximate gradients, that there exists a constant
$C>0$ independent of $\eta$ 
such that
	\begin{align}\label{eq:gradw3}
		\int_{{\rm spt}(\psi^\eta_i)}|\nabla w_i|dx\leq C(\eta_i+|x_i-y_i|)\qquad \forall i>N_\eta.
	\end{align}
Let $w\in BV(\Om)$ be a lifting of $u$ given by Theorem \ref{teo_lifting} 
and consider its extension $\overline w$ 
in $B$ as in \eqref{extension_scalar}; 
substracting a phase contribution to $u$, we 
then define in $B$
	$$w_\eta:=\overline w-\sum_{i=N_\eta+1}^{+\infty} 
w_i\psi^\eta_i,\qquad\qquad u_\eta:=e^{iw_\eta}\in W^{1,1}(B;\Suno).$$
	Let also 
$\overline u:=e^{i\overline w}$, in particular $\overline u =u$ in $\Om$.
	 Setting $V_\eta:= \cup_{i>N_\eta}{\rm spt}(\psi^\eta_i) \subset \R^2$ we infer,
using \eqref{eq:gradw3} and \eqref{eq:measure_spt_psi_i},
\begin{equation}\label{gradw4}
	\begin{aligned}
		\int_{V_\eta}	|\nabla u_\eta|dx
=&
\int_{V_\eta}|\nabla w_\eta|dx
\leq 
\int_{V_\eta}|\nabla \overline w|dx
+\sum_{i=N_\eta+1}^{+\infty}\left(\frac{C}{\eta_i}|{\rm spt}(\psi^\eta_i)|
+C(\eta_i +|x_i-y_i|)\right)
\\
\leq& \int_{V_\eta}|\nabla \overline u|dx+C\sum_{i=N_\eta+1}^{+\infty} (\eta_i +|x_i-y_i|)
\\
\leq& \int_{V_\eta\cap \Om}|\nabla \overline u|dx+C\sum_{i=N_\eta+1}^{+\infty} (\eta_i +|x_i-y_i|)+o_\eta(1),
\end{aligned}
\end{equation}
where $o_\eta(1)\rightarrow 0$ as $\eta\rightarrow 0^+$.
The presence of $\overline u$ is due to the fact that in 
general $V_\eta\setminus \overline\Om$ might be nonempty. However, since $|V_\eta\setminus \overline\Om|\rightarrow 0$ as $\eta\rightarrow 0^+$,
the last estimate in \eqref{gradw4} holds.

From this and the definition of $\eta_i$ in \eqref{eq:eta_i}
we conclude
	\begin{align}\label{claimB}
	\|u-u_\eta\|_{W^{1,1}(\Om, \R^2)}\leq 2\int_{V_\eta}|\nabla u|dx+C\eta+o_\eta(1),
	\end{align}
	where we use that $u=u_\eta$ on $\Om\setminus V_\eta$. 

Now, we claim that \begin{align}\label{claimA}
	\Det(\nabla u_\eta)=\pi\sum_{i=1}^{N_\eta}(\delta_{x_i}-\delta_{y_i}),
	\end{align} which implies in turn that 
$$\|\Det(\nabla u)-\Det(\nabla u_\eta)\|_{{\rm flat}}
=
\pi \|\sum_{i=N_\eta+1}^{+\infty}(\delta_{x_i}-\delta_{y_i})\|_{{\rm flat}}\leq \pi\sum_{i=N_\eta+1}^{+\infty}|x_i-y_i|
< \frac{\pi \eta}{2}.
$$

To show the claim, for all $m>N_\eta$ define in $B$
$$
f_m:=\overline w-\sum_{i=N_\eta+1}^m w_i\psi^\eta_i, \qquad v_m:=e^{if_m}. 
$$
Using an estimate similar to \eqref{gradw4}, \eqref{claimB},
we see that  $v_m\rightarrow u_\eta$ 
in $W^{1,1}(\Om;\R^2)$ as 
$m\rightarrow+\infty$, and therefore,
owing to the same observation leading to \eqref{flat_strong},
	\begin{align}\label{det_conv}
\lim_{m \to +\infty}
\Vert \Det(\nabla v_m) - \Det(\nabla u_\eta)\Vert_{{\rm Lip}_0(\Om)'} \to 0.
	\end{align}
and also
	\begin{align}\label{det_conv_bis}
	 \Det(\nabla v_m) \rightharpoonup \Det(\nabla u_\eta)\qquad \text{ in }\mathcal D'(\Om).
\end{align}
	On the other hand $\Det(\nabla v_m)=\pi\sum_{i=1}^{N_\eta}(\delta_{x_i}-\delta_{y_i})+\pi\sum_{i=m+1}^{+\infty}(\delta_{x_i}-\delta_{y_i})$, and since the second term tends to zero in the flat distance, we 
conclude
	\begin{align}\label{eq:det_det}
	\Det(\nabla v_m)\rightarrow \pi\sum_{i=1}^{N_\eta}(\delta_{x_i}-\delta_{y_i}), \qquad \text{ in }\text{Lip}_0(\Om)'.
	\end{align}
In particular, from \eqref{det_conv_bis} and \eqref{eq:det_det}, \eqref{claimA} follows.
From this and \eqref{claimB} it suffices to choose $\eta=\eta(\eps)$ small enough to guarantee that (ii) holds. Hence setting $N_\eps:=N_{\eta}$ and $u_\eps:=u_\eta$ the thesis follows.
\end{proof}

Now we refine the approximation of Lemma \ref{lem:density_I}.

\begin{lemma}[\textbf{Density: not collinear points}]\label{lem:density_II}
Let $u\in W^{1,1}(\Om;\Suno)$ be such that $\Det(\nabla u)
=\pi \sum_{i=1}^\indicefinalesomma
(\delta_{x_i}-\delta_{y_i})$ 
is a representation as in \eqref{T_finite2} satisfying property
$({\rm P_f})$ in Section \ref{subsec:properties_P_and_P'}, 
with $x_i,y_i\in \overline\Om$, $x_i \neq y_i$,
$i=1,\dots,\indicefinalesomma$. Then, for all $\eps>0$ there exists a map $u_\eps\in W^{1,1}(\Om;\Suno)$ such that:
\begin{itemize}
	\item[(i)] $\Det(\nabla u_\eps)=\pi\sum_{i=1}^\indicefinalesomma
(\delta_{x^\eps_i}-\delta_{y^\eps_i})$ and the points $x_i^\eps$, $y_i^\eps$ 
in $\Om$ are distinct and 
three by three not collinear;
	\item[(ii)] 
	$\|u-u_\eps\|_{W^{1,1}}+\|\Det(\nabla u)-\Det(\nabla u_\eps)\|_{{\rm flat}}<  \eps$.
\end{itemize}
\end{lemma}

\begin{proof}
Define
$$
I^+:=\{i \in \{1,\dots,\indicefinalesomma\}:x_i\in \Om\},\qquad 
I^-:=\{i \in \{1,\dots,\indicefinalesomma\}:y_i\in \Om\}.
$$ 
Fix $\eta>0$.
For all $i\in I^+$ let us choose $\widehat x_i,\widehat y_i\in \Om$ with $\widehat y_i:=x_i$ and in such a way that the 
points $\widehat x_i$,  $i\in I^+$,
 are all distinct,  three by three not collinear, and satisfy
\begin{align}\label{lunghezze}
\sum_{i\in I^+}|\widehat x_i-\widehat y_i|
< \eta.
\end{align}
For all $i\in I^+$, let $\widehat w_i:=w_{\widehat x_i,\widehat y_i}$ be the dipole map defined in \eqref{dip_map},  and let $\psi^\eta_i:\R^2\rightarrow\R$ be the cut-off function given by
\begin{align*}
\psi^\eta_i(x)=\varrho\Big(
\frac1\eta d(x,\overline{\widehat x_i\widehat y_i})\Big),
\end{align*}
where $\varrho$ is as in the proof of Lemma \ref{lem:density_I}.
In particular $\psi_i^\eta$ is 
Lipschitz continuous with Lipschitz constant $\frac{2}{\eta}$ and is supported 
in $V^\eta_i:=\{x\in \R^2: d(x,\overline{\widehat x_i\widehat y_i})\leq \eta\}$. Supposing $\eta>0$ 
sufficiently small, we have 
$V_i^\eta \subset \Om$.
Now, using also \eqref{lunghezze}, 
notice that
\begin{equation}\label{eq:meas_V_eta}
|V^\eta_i|=\pi\eta^2+2\eta|{\widehat x_i-\widehat y_i}|\leq C\eta^2
\qquad\forall i\in I^+,
\end{equation}
where $C>0$ is a constant independent of $i$ and $\eta$.
Further, by \eqref{gradw} and \eqref{lunghezze}, we 
deduce that there is a constant, still denoted by $C>0$,
and  independent of $\eta$ 
and $i$,  such that 
\begin{align}\label{gradw2}
	\int_{V^\eta_i}|\nabla \widehat w_i|dx\leq C\eta\qquad \forall i\in I^+.
\end{align}
Similarly for all $i\in I^-$ we choose $\widetilde x_i,\widetilde y_i\in \Om$ with $\widetilde x_i:=y_i$ 
and in such a way that the points $\widehat x_i$, $i\in I^+$, 
and $\widetilde y_i$,
$i\in I^-$ are all distinct, three by three not collinear, 
and satisfy 
\begin{align}\label{lunghezze2}
	\sum_{i\in I^-}|\widetilde x_i-\widetilde y_i|
< \eta.
\end{align}
Also in this case we introduce, for $i\in I^-$, the maps $\widetilde w_i:=w_{\widetilde x_i,\widetilde y_i}$ and   $\phi_i^\eta:\R^2\rightarrow\R$, the latter  defined as $\phi^\eta_i(x):=\max\{0,1-\frac{1}{\eta}
d(x,\overline{\widetilde x_i\widetilde y_i})\},$ which 
enjoy the same features of $\psi_i^\eta$; in particular, the supports  $W_i^\eta$ of $\phi_i^\eta$, $i\in I^-$, are contained in $\Om$ and 
have Lebesgue measures 
bounded by $C\eta^2$.  The same estimate as in \eqref{gradw2} holds for $\widetilde w_i$.

Eventually, let us consider a
 lifting $w\in BV(\Om)$ 
of $u$ provided by Theorem \ref{teo_lifting}
 and \eqref{extension_scalar}; 
we define 
\begin{align*}
w_\eta:=w+\sum_{i\in I^+}\psi_i^\eta \widehat w_i+\sum_{i\in I^-}\phi_i^\eta\widetilde w_i,\qquad\qquad v_\eta:=e^{iw_\eta}.
\end{align*}
Owing to the fact that $w_\eta=w$ out of 
$\bigcup_{i\in I^+}V_i^\eta\cup\bigcup_{i\in I^-} W_i^\eta$,
it is immediate
 that $v_\eta\rightarrow u$ 
in $L^1(\Om;\R^2)$ as $\eta\rightarrow 0^+$. 
Since $|\nabla v_\eta|=|\nabla w_\eta|$ a.e. in $\Om$, 
we can estimate
\begin{align*}
|\nabla v_\eta|\leq |\nabla w|+\frac2\eta\left(\sum_{{i\in I^+}}\big\|\widehat w_i\|_{L^\infty}+\sum_{{i\in I^-}}\|\widetilde w_i\|_{L^\infty}\right)+\sum_{{i\in I^+}}|\nabla \widehat w_i|+\sum_{{i\in I^-}}|\nabla \widetilde w_i|.
\end{align*}
From this, in view of the fact that $v_\eta=u$ in $\Om\setminus 
(\bigcup_{i\in I^+}V_i^\eta\cup\bigcup_{i\in I^-} W_i^\eta)$, 
using \eqref{gradw2} we conclude
\begin{align*}
\|\nabla v_\eta-\nabla u\|_{L^1}\leq \|\nabla u\|_{L^1(A_\eta)}
+C\frac{|{A_\eta}|}{\eta}+ C\eta.
\end{align*}
The right-hand side being 
negligible as $\eta\rightarrow 0^+$
(see \eqref{eq:meas_V_eta}), we conclude 
\begin{equation}\label{vetaconv}
\lim_{\eta \to 0^+}
v_\eta = u \text{ in }W^{1,1}(\Om;\R^2).
\end{equation}
Furthermore, using \eqref{det=curl}, 
we readily see that
\begin{align*}
& \Det(\nabla v_\eta)=
\frac{1}{2} {\rm Curl}(\nabla w_\eta)
=\frac{1}{2}\left(\curl(\nabla w)+\sum_{{i\in I^+}}\curl(\psi_i^\eta \widehat w_i)+\sum_{{i\in I^-}}\curl(\phi_i^\eta \widetilde w_i)\right)
\\
=& \Det(\nabla u)+\sum_{{i\in I^+}}\Det(\nabla v_{\widehat p_i,\widehat n_i})+\sum_{{i\in I^-}}\Det(\nabla v_{\widetilde p_i,\widetilde n_i}),
\end{align*}
which implies
\begin{equation*}
\begin{aligned}
& 
\Det(\nabla u)
- \Det(\nabla v_\eta)
=
-\sum_{{i\in I^+}}\Det(\nabla v_{\widehat x_i,\widehat y_i})-
\sum_{{i\in I^-}}\Det(\nabla v_{\widetilde x_i,\widetilde y_i})
\\
=& 
- \sum_{i\in I^+}
(\delta_{\widehat x_i}-\delta_{\widehat y_i})
- \sum_{i \in I^-}
(\delta_{\widetilde x_i}-\delta_{\widetilde y_i})
\end{aligned}
\end{equation*}
in $\mathcal D'(\Om)$.
Thus,  using \eqref{lunghezze} and \eqref{lunghezze2}, 
we get
$$\|
\Det(\nabla u)
- \Det(\nabla v_\eta)
\|_{{\rm flat}}\leq 2\eta.
$$
 In particular, from this and \eqref{vetaconv}, 
setting $u_\eps:=v_\eta$ for $\eta>0$ small enough, the thesis follows.
\end{proof}

\begin{Remark}\rm
The noncollinearity condition will be used in the proof 
of Theorem \ref{teo:main1}, 
to guarantee the validity of condition \eqref{eq:reason_of_noncollinearity}. 
\end{Remark}

The approximating maps in Lemma \ref{lem:density_II}
can be suitably refined around the singular points as follows.

\begin{lemma}[\textbf{Density: behaviour near $x_i$, $y_i$}]
\label{lem:density_III}
Let  $u\in W^{1,1}(\Om;\Suno)$ be such 
that $\Det(\nabla u)=\pi\sum_{\indicegenerico
=1}^n(\delta_{ x_\indicegenerico}-\delta_{ y_\indicegenerico})$, 
with $x_\indicegenerico, y_\indicegenerico
\in \overline\Om$, $\indicegenerico=1,\dots,\indicefinalesommalemma$; let us assume that the points $ x_\indicegenerico, y_\indicegenerico$ which belong to $\Om$ are distinct and three by three not collinear, as in 
the thesis of  Lemma \ref{lem:density_II}.
Then, for all $\eps>0$ there exists a map $u_\eps\in W^{1,1}(\Om;\Suno)$ such that: 
\begin{itemize}
	\item[(i)] there exist positive numbers $\rho_{x_\indicegenerico}
<\eps$ and $\rho_{y_\indicegenerico}<\eps$ such that 
the discs of $\{B_{\rho_{x_\indicegenerico}}(x_i),B_{\rho_{y_\indicegenerico}}(y_i):
x_\indicegenerico\in\Om, y_i\in \Om\}$ are contained in $\Omega$ and 
pairwise disjoint, and 
$u_\eps=e^{i\theta_{x_\indicegenerico}}$ in $  B_{\rho_{x_\indicegenerico}}(x_\indicegenerico)$, $u_\eps=e^{-i\theta_{y_\indicegenerico}}$ in $  B_{\rho_{y_\indicegenerico}}(y_\indicegenerico)$;
	\item[(iii)]
	$\|u-u_\eps\|_{W^{1,1}}+\|\Det(\nabla u)-\Det(\nabla u_\eps)\|_{{\rm flat}}
< \eps$.
\end{itemize}
\end{lemma}
\begin{proof}
Let $\{z_\indicegenerico:\indicegenerico=1,\dots,N\}$ be the set, 
suitably relabelled, of those points
among the $x_\indicegenerico$'s and $y_\indicegenerico$'s which belong
to $\Om$. 
 Moreover, let $r>0$ be 
small enough so that 
the discs $B_{r}(z_\indicegenerico)$ are 
contained in $\Om$ and pairwise disjoint. 
We can choose $r>0$ arbitrarily 
small so that $u\res \partial B_r(z_\indicegenerico)\in W^{1,1}
(\partial B_r(z_\indicegenerico);\Suno)$
 for all $\indicegenerico=1,\dots,N$. We show how to modify $u$ 
in one of these discs, say $B_r(z_1)$, 
and then proceed similarly for the other discs. 

Let us assume without loss of generality that $z_1=x_\indicegenerico=0$ 
for some $\indicegenerico$ (i.e., $x$ is a positive pole 
at the origin), and write $B_r=B_r(0)$ in place of $B_r(x_\indicegenerico)$. 
Since in $B_r$ we have $\Det(\nabla u)=\pi\delta_{0}$, 
it is not difficult to see that 
\begin{align}\label{grado_u}
	\frac12\int_{\partial B_s}\left(
u_1\frac{\partial u_2}{\partial \tau} -u_2\frac{\partial u_1}{\partial \tau}\right)d\mathcal H^1=\pi\deg(u;\partial B_r)=\pi
\end{align}
for all $s\in(0,r)$ such that $u\res \partial B_s\in W^{1,1}(\partial B_s;\Suno)$. 
By the mean value theorem, we fix $d=d_\indicegenerico\in (r/2,r)$ so that 
\begin{align}\label{stima_sfera}
	\int_{\partial B_d}|\nabla u|d\mathcal H^1\leq \frac2 r\int_{r/2}^r\int_{\partial B_s}|\nabla u|d\mathcal H^1ds=\frac2r\int_{B_r\setminus B_{r/2}}|\nabla u|dx,
\end{align}
and $u\res\partial B_d\in W^{1,1}(\partial B_d;\Suno)$.
Let $\theta_u\in BV(\partial B_d)$ 
denote a lifting of $u\res \partial B_d$ such that, 
owing  to \eqref{grado_u}, $\theta_u$ has a unique jump point (say 
at $(d, 0)\in \partial B_d$) with $\jump{u}=2\pi$. 
Consider a polar coordinate system $(\rho,\theta)$ 
around $0$, and define
$H:B_d\setminus \overline B_{d/2}\rightarrow\R$  as
\begin{align*}
H(x):=2\theta_u\Big(d\frac{x}{|x|}\Big)\frac{|x|-d/2}{d}+\theta(x)\frac{d-|x|}{d}.
\end{align*}
The function $H$ has a jump of size $2\pi$ 
on the segment with endpoints $(d/2,0)$ and $(d,0)$. 
Also, $e^{iH}\in W^{1,1}(B_d\setminus 
\overline B_{d/2})$, it equals $u$ on $\partial B_d$ and $\frac{x}{|x|}$ on $\partial B_{d/2}$. 
We set
\begin{align*}
u_\eps(x):=\begin{cases}
	u(x)&\text{ if } x \in \Om\setminus B_d,\\
	e^{iH(x)}&\text{ if } x \in B_d\setminus \overline B_{d/2},
\\
	\frac{x}{|x|}&\text{ if } x \in B_{d/2},
\end{cases}
\end{align*}
in particular $u_\eps\in W^{1,1}(\Om;\Suno)$. 
Let us estimate the gradient of $H$; we have
\begin{align*}
	\nabla H(x)=&2
\dot \theta_u\left(\frac{d x}{|x|}\right)
\left(\frac{Id-\frac{x}{|x|}\otimes\frac{x}{|x|}}{|x|}\right)(|x|-d/2)
\\&+
2\theta_u\left(\frac{d x}{|x|}\right)\frac{x}{d|x|}
	+\nabla\theta(x)\frac{d-|x|}{d}-\theta(x)\frac{x}{d|x|},
\end{align*} 
where $\dot \theta_u$ denotes the (absolutely continuous part of the)
derivative of $\theta_u$.
Therefore, using that $|x|\in (d/2,d)$ for $x\in B_d\setminus \overline 
B_{d/2}$, there is a
constant $C>0$ independent of $d$ such that 
\begin{align*}
|	\nabla H(x)|\leq 2\left|\dot \theta_u\left(\frac{d x}{|x|}\right)\right|+\frac{C}{d}\qquad\qquad\text{for a.e. }x\in B_d\setminus 
\overline B_{d/2}.
\end{align*}
On the other hand, since $e^{i\theta_u}=u$ 
on $\partial B_d$, 
we have $\left|\dot \theta_u\left(\frac{d x}{|x|}\right)\right|=\left|\nabla u\left(\frac{d x}{|x|}\right)\right|$, and 
integrating on $B_d\setminus \overline B_{d/2}$ we get
\begin{align*}
	\int_{B_d\setminus \overline 
B_{d/2}}|\nabla u_\eps|dx&=\int_{B_d\setminus \overline B_{d/2}}|\nabla H|dx
\leq Cd+2\int_{d/2}^d\int_{\partial B_s}\left|\nabla u \left(\frac{d x}{|x|}\right)\right|d\mathcal H^1(x) ds\\
	&=Cd +C\int_{d/2}^d\int_{\partial B_d}|\nabla u|d\mathcal H^1ds\leq Cd+C\int_{B_r\setminus \overline B_{r/2}}|\nabla u|dx,
\end{align*}
where we have used \eqref{stima_sfera} and that $r/2<d<r$ in the last inequality.

Now, applying a similar modification of $u$ in the other discs centered at $z_i$,
we can finally estimate the distance between $u$ and $u_\eps$ in $W^{1,1}(\Om;\R^2)$, namely
\begin{align*}
&\|u-u_\eps\|_{L^1}\leq n\pi r^2,\\
&\|\nabla u-\nabla u_\eps\|_{L^1}\leq NCd+ C\sum_{\indicegenerico=1}^N
\int_{B_r(z_i)}|\nabla u|dx+\sum_{\indicegenerico=1}^N
\int_{B_{d_\indicegenerico/2}(z_\indicegenerico)}\left|\nabla \left(\frac{x}{|x|}\right)\right|dx.
\end{align*}
Since $r$ can be chosen 
arbitrarily small, 
the sum of the above right-hand sides can be bounded by $\eps$. Observing that $\Det(\nabla u)=\Det(\nabla u_\eps)$,  
the thesis follows by setting 
$\rho_{z_\indicegenerico}:=d_\indicegenerico/2$, 
for all $\indicegenerico =1,\dots,N$.
\end{proof}

\section{Proof of Theorem \ref{teo:upper_bound_W11Suno}}\label{sec:main_results}
In this section we prove Theorem 
\ref{teo:upper_bound_W11Suno}. 
Recalling the definition of $\Vert \cdot\Vert_{{\rm flat},\alpha}$ in \eqref{eq:def_M}, 
we start with:

\begin{theorem}\label{teo:main1}
Let $u\in W^{1,1}(\Om;\Suno)$. Suppose:
\begin{itemize}
	\item[(i)] $\frac1\pi\Det(\nabla u)=\sum_{i=1}^\indicefinalesommateorema
(\delta_{x_i}-\delta_{y_i})
=: T$ admits a representation in $\overline\Om$ satisfying ($P_f$) and such that the points $x_i,y_i$ belonging to $\Om$ are
distinct and three by three not collinear;
	\item[(ii)] 
there exists 
$\piccoloraggio >0$ such that the discs 
$B_{\piccoloraggio }(x_i)$, $B_{\piccoloraggio }(y_j)$ with $x_i,y_j\in \Om$, are 
contained in $\Om$ and pairwise disjoint, and 
$u=e^{i\theta_{x_i}}$ in $  B_{\piccoloraggio }(x_i)$, $u=e^{-i\theta_{y_j}}$ in $  B_{\piccoloraggio }(y_j)$.
\end{itemize}
Then 
\begin{equation}\label{eq:ril_area_upper_estimate}
	\rilarea(u,\Om)\leq \int_{\Om}\sqrt{1+|\nabla u|^2}~dx+
\Vert 
{\rm Det}(\grad u)
\Vert_{{\rm flat},\alpha}.
\end{equation}
\end{theorem}
\begin{proof}
We need to exhibit\footnote{See \eqref{eq:u_r} below.} a sequence $(u_r)\subset C^1(\Om; \R^2)$ converging to 
$u$ in $L^1(\Om; \R^2)$ such that $\liminf_{r \to 0^+}
\area(u_r, \Om)$ is less than or equal to the right-hand side
of \eqref{eq:ril_area_upper_estimate}.

For the measure 
$T$
we consider 
currents $R_{\textrm{min}}\in \mathcal D_0(\Om)$ and 
$S_{\textrm{min}}\in \mathcal D_1(\Om)$ 
given by Lemma \ref{lem:structure_of_minimizers_of_the_combinatorial_problem}. 
After relabelling, we write 
$$
R_{\textrm{min}}=\sum_{i=1}^k\sigma_i\delta_{z_i}, 
\ \ \sigma_i\in \{-1,+1\}, ~z_i\in \Om, \ k \leq \indicefinalesommateorema,\qquad 
\qquad S_{\textrm{min}}=\sum_{j\in J}\jump{\overline{y_jx_j}},
\ \ J\subset\{1,\dots,\indicefinalesommateorema\},
$$
with $T = R_{{\rm min}} + \partial S_{{\rm min}}$
(it may happen that $k=0$, in which
case we understand $R_{{\rm min}}=0$, 
or that $J = \emptyset$, in which case $S_{{\rm min}}=0$).
By Lemma \ref{lem:structure_of_minimizers_of_the_combinatorial_problem} (b),
the segment $\overline{x_jy_j}$ is contained in $\Om$, with the 
only (possible) exception 
of an endpoint (thanks to condition $({\rm P_f})$). We will 
work in a disc $B\supset \supset \Om$, that we fix from now on.

Take $r\in (0,\piccoloraggio /2)$, and consider the set 
$\{B_{2r}({z_i}): i=1,\dots,k\}$
(the $z_i$'s are among the $x_j$, $y_j$ and,
being contained in $\Om$, satisfy assumption (ii)); these discs are 
contained in $\Om$,  and the tubular neighborhoods 
$$
T_t(\overline{x_j y_j}):=\{x\in B: d(x,\overline{x_jy_j})<t\}, \qquad
j\in J,
$$
of $ \overline{x_jy_j}$ are 
disjoint from this family of discs. Moreover, by hypothesis (i), 
due to noncollinearity, 
the segments $\overline{x_jy_j}$ are 
pairwise disjoint (see Lemma \ref{lem:structure_of_minimizers_of_the_combinatorial_problem}), and so 
for all $t>0$ sufficiently small 
\begin{equation}\label{eq:reason_of_noncollinearity}
j_1, j_2 \in J, \ j_1 \neq j_2 \Rightarrow 
T_t(\overline{x_{j_1} y_{j_1}}) 
\cap 
T_t(\overline{x_{j_2} y_{j_2}}) 
= \emptyset.
\end{equation}
Set also
	\begin{align}\label{V_jt}
	V_t:=\bigcup_{j\in J}T_t(\overline{x_j y_j}).
	\end{align}	
If $k \geq 1$, for all $i=1,\dots,k$, we fix a 
simple polygonal\footnote{I.e., not self-intersecting and
obtained by a finite number of 
concatenations of segments.} curve $\connessionepoloi$ 
starting 
at $z_i$ and reaching the 
external boundary\footnote{$\dOm$ is, in general, not
connected, and consists of a finite number of loops. The external
boundary of $\dOm$ is the loop whose interior contains
all the others.} of $\partial \Om$. 
Curves $\connessionepoloi$'s can be chosen mutually 
disjoint, and disjoint from
$\overline V_t$. Further, it is convenient to extend $\connessionepoloi$ 
(keeping the same notation)
in order to reach $\partial B$ transversely. 
We set $
\ingrassatoconnessionepoloi
:=
\{x\in B: d(x,\connessionepoloi)<2r\}$, $i=1,\dots,k$,
and observe that $ B_{2r}(x_i)\subset \ingrassatoconnessionepoloi$ for all $i=1,\dots,k$.
If $r$ is small enough the elements of 
the family $\{\ingrassatoconnessionepoloi: i=1,\dots,k\}$ do not intersect each other, 
and moreover (choosing smaller $t$ and $r$ if necessary) $\overline{
\ingrassatoconnessionepoloi}\cap \overline V_t=\emptyset$ for all $i=1,\dots,k$. 
Consider also the connected curves $\gamma_{z_i}^{+,r}$
and $\gamma_{z_i}^{-,r}$
  which run parallel to 
$\connessionepoloi$ at distance $r$,
defined as
$$\gamma_{z_i}^{\pm,r}:=\{x\in \ingrassatoconnessionepoloi\setminus B_{r}(z_i):
\overline d(x,\connessionepoloi)
=\pm r\},$$
where $\overline 
d$ denotes a signed distance from $\connessionepoloi$ (defined in a suitable neighborhood of $\connessionepoloi$).
For every connected component $\partial_\ell \Om$
of $\partial \Om$ different from the external boundary ($\ell\in L$, 
$L$ some 
finite set of indices), 
we consider a simple polygonal curve $\connessionebordo_\ell
\subset B$ connecting $\partial_\ell \Om$ 
to $\partial B$, disjoint
from $\overline V_t$, from $\cup_{i=1}^k
\overline{\ingrassatoconnessionepoloi}$ and from 
$\partial_{\ell'}\Om$,
$\ell' \neq \ell$. Extending
slightly $\omega_\ell$ inside 
$\Om_\ell$, $\Om_\ell$ the region 
outside $\Om$ and enclosed by $\partial_\ell\Om$, we assume that $\omega_\ell$ starts at a point $\Delta_\ell\in \Om_l$ with 
$B_{2r}(\Delta_\ell)\subset\Om_\ell$, for all $\ell\in L$. 
Together with this we consider the connected curves $\connessionebordo_\ell^{+,r}$ 
and $\connessionebordo_\ell^{-,r}$ which run parallel 
to $\connessionebordo_\ell$ at distance $r$ from each side,  and join 
$\overline B_{r}(\Delta_\ell)$ with the external boundary,
$$\connessionebordo_\ell^{\pm,r}
 \subset \left\{x\in B\setminus B_{r}(\Delta_\ell):
\overline d(x,\connessionebordo_\ell)=\pm r\right\},
$$
where, again, $\overline 
d$ denotes a signed distance from $\omega_\ell$ 
(defined in a suitable neighborhood of $\omega_\ell$).
We may assume, choosing smaller $r$ and $t$ if necessary,  that 
 the curves $\connessionebordo_\ell,\connessionebordo_\ell^{\pm,r}$
are pairwise 
disjoint and do not intersect $\overline V_t\cup 
\cup_{i=1}^k \overline{\ingrassatoconnessionepoloi}
\cup \cup_{\ell' \neq \ell}
\partial_{\ell'}\Om$. 
Finally, we define 
\begin{equation}\label{eq:B_pm}
B^\pm:=B\setminus 
\Big[
\bigcup_{i=1}^k 
\overline B_{r}(z_i) ~\cup~ 
\bigcup_{j\in J}\overline{x_jy_j} ~\cup~
\bigcup_{i=1}^k \gamma_{z_i}^{\pm,r}~\cup~ 
\bigcup_{\ell \in L}\connessionebordo_\ell^{\pm,r}\cup \bigcup_{\ell \in L}B_r(\Delta_\ell)\Big].
\end{equation}
Using our assumptions it follows that $B^+$ and $B^-$ are connected; 
however, they are not necessarily simply connected. 
By 
construction, for any closed simple Lipschitz 
curve $\alpha:\Suno\rightarrow \Om\cap B^\pm$ such that $u^\alpha:=
u\circ \alpha\in W^{1,1}(\Suno;\Suno)$, we have
\begin{align*}
\frac{1}{2 \pi}	\int_{\Suno}(u^\alpha_1\nabla^\perp {u^\alpha_2} -
u_2^\alpha\nabla^\perp {u_1^\alpha})\cdot \dot\alpha\; ds=0,
\end{align*}
since the left-hand side is the 
degree of $u$ on the boundary of the domain 
enclosed by the support of $\alpha$,
and such curves cannot enclose any connected component of $B\setminus \Om$ due to the presence of 
$\bigcup_{\ell 
\in L}\connessionebordo_\ell$, 
and cannot enclose any single pole due to the presence of $\connessionepoloi$'s (note that they can enclose some 
segment $\overline{x_j y_j}$). In particular, there exist 
two liftings\footnote{if $J=\emptyset$, i.e., no dipoles,
and if $\Om$ is simply connected, then we can take 
$w_+ = w_-$.}
$w_\pm$ of $u$ with 
$$
w_+
\in W^{1,1}(\Om\cap B^+), \qquad w_-\in W^{1,1}(\Om\cap B^-).
$$
For $w_+$ and $w_-$ we consider (not-relabelled) 
extensions $w_+\in W^{1,1}(B^+)$ and $w_-\in W^{1,1}(B^-)$ as in \eqref{extension_scalar}.
We are now going to suitably smoothen these liftings through a function
$w_r$, that will allow us to eventually
define the map $u_r$ in \eqref{eq:u_r}.

We may assume that  
\begin{align}\label{wequals}
w_+=w_-\qquad \text{on }B^+\setminus \left(\bigcup_{i=1}^k\overline{T_r(\gamma_{z_i})}
\cup \bigcup_{\ell \in L}
\overline{\effeell}
\right)=B^-\setminus \left(\bigcup_{i=1}^k\overline{T_r(\gamma_{z_i})}\cup \bigcup_{\ell \in L}
\overline{\effeell} \right),
\end{align}
where $T_r(\connessionepoloi)
\subset \{x\in B:\overline d(x,
\connessionepoloi)\in(-r,r)\}$ 
is the region enclosed by $\gamma_{z_i}^{+,r}$ 
and $\gamma_{z_i}^{-,r}$ and  
$\effeell:=
\{x\in B\setminus \cup_{\ell\in L}B_{r}(\Delta_\ell):
d(x,\connessionebordo_\ell)\in(-r,r)\}$ is the tubular neighborhood of $\connessionebordo_\ell$ enclosed by 
$\connessionebordo^{+,r}_\ell$
and
$\connessionebordo^{-,r}_\ell$.
In addition, since the degree of $u$ around $z_i$ is $\sigma_i$, we 
see\footnote{The degree around a pole is computed
using counterclockwise turns, and this implicitely
determines an orientation of the jump of $w_\pm$.} that 
\begin{align}\label{jump_pm}
	\jump{w_+}=2\pi\sigma_i 
\ \ \mathcal H^1-\text{a.e. on }
\gamma_{z_i}^{+,r}, 
\quad \jump{w_-}=2\pi\sigma_i\ \ \mathcal H^1-\text{a.e. on }
\gamma_{z_i}^{-,r} 
\qquad
\forall i=1,\dots,k.
\end{align}
Furthermore 
\begin{align*}
	\jump{w_+}=\jump{w_-}=2\pi\qquad  \mathcal H^1-\text{a.e. on }
\overline{x_jy_j}
\qquad \forall j\in J,
\end{align*}
and 
$$\jump{w_\pm}\in 2\pi\mathbb Z \qquad  \mathcal H^1-\text{a.e. on }
\bigcup_{\ell \in L}\connessionebordo_\ell^{\pm,r}.$$
From \eqref{jump_pm} it follows that
\begin{align}\label{wpm}
w_+=w_-+2\pi\sigma_i\quad \text{ a.e. in }T_r(\connessionepoloi)
 \qquad
\forall i=1,\dots,k.
\end{align}
Similarly, given
$\ell\in L$, there exists 
$h_\ell\in \mathbb Z$ such that 
\begin{align}\label{wpmbis}
w_+=w_-+2\pi h_\ell\quad \text{ a.e. in }\effeell.
\end{align}
Finally, we extend $w_\pm$ to $0$ on $\bigcup_{i=1}^k\overline{
	B_{r}(z_i)}\cup \bigcup_{\ell\in L}\overline{B_r(\Delta_\ell)}$, and mollify $w_\pm$ 
using a kernel $\varrho_r$ supported in $B_{r/4}(0)$.
In particular, using also \eqref{wpm} and \eqref{wpmbis}, we infer that the
traces of the mollifications on $\connessionepoloi$ and $\omega_\ell$ satisfy
\begin{align*}
&w_+*\varrho_r
=w_-*\varrho_r
+2\pi\sigma_i\qquad\mathcal H^1\text{-a.e. on }
\connessionepoloi,\\
&w_+*\varrho_r
=w_-*\varrho_r
+2\pi h_\ell\qquad\mathcal H^1\text{-a.e. on }
\omega_\ell,
\end{align*}
and therefore, setting $B_r^- := \{x \in B : d(x, \partial B) > r\}$ and
defining $w_r: B_r^-\setminus (\cup_{i=1}^k{\connessionepoloi}
\cup\cup_{i=1}^k \overline{
B_{r}(z_i)
}\cup \cup_{\ell 
\in L}\connessionebordo_\ell\cup\cup_{\ell 
\in L}\overline\Om_\ell)
\to \R^2$ 
as
\begin{align*}
w_r:=\begin{cases}
	w_+*\varrho_r
&\text{ in }B_r^-\setminus (\cup_{i=1}^k{T_r(\gamma_{z_i})})\setminus (\cup_{i=1}^k\overline{
	B_{r}(z_i)
})\setminus (\cup_\ell T_{r}(\omega_\ell))\setminus (\cup_\ell\Om_\ell),\\
	w_-*\varrho_r
&\text{ in }\cup_{i=1}^k\{x\in 
T_r(\gamma_{z_i})
\setminus  \overline{B_{r}(z_i)}:\overline d(x,\connessionepoloi)\in(0,r)\},\\
	w_+*\varrho_r
&\text{ in }\cup_{i=1}^k\{x\in 
T_r(\gamma_{z_i})
\setminus  \overline{B_{r}(z_i)}:\overline d(x,\connessionepoloi)\in(-r,0)\},\\
	w_-*\varrho_r
&\text{ in }\cup_{\ell \in L}\{x\in \effeell \setminus 
\Om_\ell
:\overline d(x,\connessionebordo_\ell)\in(0,r)\}\\
		w_+*\varrho_r
&\text{ in }\cup_{\ell \in L}\{x\in \effeell \setminus 
\Om_\ell
:\overline d(x,\connessionebordo_\ell)\in(-r,0)\},
\end{cases}
\end{align*}
we see that $w_r\in C^\infty( 
B_r^-\setminus (\cup_{i=1}^k{\connessionepoloi}
\cup \cup_{i=1}^k \overline{
	B_{r}(z_i)
}\cup \cup_{\ell 
\in L}\connessionebordo_\ell)
\cup \cup_{\ell 
	\in L}\overline\Om_\ell),$ and 
\begin{align}\label{jumpu_r}
&\jump{w_r}=2\pi \sigma_i\qquad \mathcal H^1\text{-a.e. on }\connessionepoloi,\; i=1,\dots, k,\nonumber\\
&\jump{w_r}=2\pi h_\ell\qquad \mathcal H^1\text{-a.e. on }\connessionebordo_\ell,\; \ell\in L.
\end{align}
Eventually, for all $i=1,\dots,k$, by the 
assumptions on $u$ and the choice of $r \in (0,\piccoloraggio /2)$, we have
$u(x)=e^{i\sigma_i\theta_{z_i}}$ 
for $x \in B_{2r}(z_i) \setminus \{z_i\}$
for suitable $\sigma_i \in \{\pm 1\}$.

Thus
\begin{align*}
w_\pm-\sigma_i\theta_{z_i}\in 2\pi\mathbb Z\qquad \text{ in }B_{2r}(z_i)\setminus B_{r}(z_i).
\end{align*} 
Assuming without loss of generality that $\theta_{z_i}$ jumps on $\connessionepoloi$ in $B_{2r}(z_i)\setminus B_{r}(z_i)$, by \eqref{wequals}, for all $i=1,\dots,k$ we 
find an integer $\zeta_i$ such that 
\begin{align*}
	w_+=w_-=\sigma_i\theta_{z_i}+2\pi \zeta_i \quad \text{ in }B_{2r}(z_i)\setminus B_{r}(z_i)\setminus T_{r}(\connessionepoloi),
\end{align*}
whereas in $B_{2r}(z_i)
\cap T_{r}(\connessionepoloi)
\setminus 
B_{r}(z_i)$ 
we have 
\begin{align*}
	&w_+=\sigma_i\theta_{z_i}+2\pi \zeta_i-2\pi\sigma_i 
 & \text{in }\{x:\overline d(x,\connessionepoloi)\in(0,r)\},
\\
	&w_+=\sigma_i\theta_{z_i}+2\pi \zeta_i-2\pi(\sigma_i+1) & \text{in }\{x:\overline d(x,\connessionepoloi)\in(-r,0)\},\\
	&w_-=\sigma_i\theta_{z_i}+2\pi \zeta_i &  \text{in }\{x:\overline d(x,
\connessionepoloi)\in(0,r)\},
\\
	&w_-=\sigma_i\theta_{z_i}+2\pi (\zeta_i-1) &  \text{in }\{x:\overline d(x,\connessionepoloi)\in(-r,0)\}.
\end{align*}
Therefore
\begin{align*}
	&w_r-\sigma_i\theta_{z_i}*\varrho_r\in2\pi 
\mathbb Z&\text{ in } (B_{\frac53r}(z_i)\setminus B_{\frac43r}(z_i))\setminus 
T_r(\connessionepoloi)\\
	&w_r-\sigma_i\widehat\theta_{z_i}*\varrho_r\in2\pi \mathbb Z&\text{ 
in }(B_{\frac53r}(z_i)\setminus B_{\frac43r}(z_i))\cap T_r(\connessionepoloi),
\end{align*}
where $\sigma_i\widehat \theta_{z_i}$ is any 
lifting of $\frac{x-z_i}{|x-z_i|}$ 
which is continuous in $T_r(\connessionepoloi)$. 

We introduce a nondecreasing cut-off function 
$\psi:[0,2r]\rightarrow [0,1]$ of class $C^1$ such that $\psi=0$ 
on $[0,\frac43r]$, $\psi=1$ on $[\frac53r,2r]$, with $\psi'\leq \frac{12}{r}$.
Finally, we define 
\begin{align}\label{eq:u_r}
u_r(x):=\begin{cases}
e^{iw_r(x)}&\text{if }x\in B_r^-\setminus \cup_{j=1}^k B_{2r}(z_j),
\\
e^{iw_r(x)}\psi(|x-z_j|)&\text{if }x\in B_{2r}(z_j)\text{ for some }j=1,\dots,k,
\end{cases}
\end{align}
where we extend $w_r$ in $B$ to $0$ outside its domain.
In particular
$u_r(x)=e^{i\sigma_i\widehat\theta_{z_i}*\varrho_r}\psi(|x-z_i|)$ for
$x \in (B_{\frac52r}(z_i)\setminus B_{\frac43r}(z_i))\cap T_r(\connessionepoloi)$,
for any $i=1,\dots,k$.
We also observe that if we suppose that the kernel $\varrho_r$ 
is radial, a direct computation shows that $\widehat\theta_{z_i}*\varrho_r=\widehat\theta_{z_i}$ in $B_{\frac53r}(z_i)\setminus B_{\frac43r}(z_i)$. 
So  
$$
u_r(x)=e^{i\sigma_i\widehat\theta_{z_i}}\psi(|x-z_i|)\qquad
\forall x \in (B_{\frac53r}(z_i)\setminus B_{\frac43r}(z_i))\cap 
T_r(\connessionepoloi).
$$
Now, $u_r$ is of class $C^1$, $|u_r|\leq 1$ and it is straightforward that $u_r\rightarrow u$ 
pointwise almost everywhere in $\Om$
as $r \to 0^+$. In particular, $\lim_{r \to 0^+}
u_r = u$
in $L^1(\Om; \R^2)$. 

We are now in a position to 
estimate the graph area of the map $u_r$. In order to 
estimate it in $\Om\setminus \cup_i \overline B_{2r}(x_i)$ 
it is convenient to consider a lifting $$w\in W^{1,1}\big(B\setminus (\bigcup_{i=1}^k
(\connessionepoloi\cup \overline B_{r}(z_i))\cup \bigcup_{j\in J}\overline{x_jy_j})\cup \bigcup_{\ell\in L}(\omega_\ell\cup \overline B_r(\Delta_\ell))\big),$$ which coincides with $w_\pm$ in the set in \eqref{wequals}. Such a lifting $w\in BV(\Om\setminus \cup_{i=1}^k \overline B_{r}(z_i) )$ 
satisfies
\begin{align*}
&\jump{w}=2\pi\sigma_i\qquad\mathcal H^1\text{-a.e. on }\connessionepoloi,\;i=1,\dots,k,\nonumber\\
&\jump{w}=2\pi\qquad\mathcal H^1\text{-a.e. on }\overline{x_jy_j},\;j\in J,\nonumber\\
&\jump{w}=2\pi h_\ell\qquad\mathcal H^1\text{-a.e. on }\connessionepoloi,\;\ell\in L.
\end{align*}
Notice that 
$\lim_{r \to 0^+}w_r = w$ strictly in $BV(V_t)$ 
(for $t$ small enough as in \eqref{V_jt}), 
since $w_r=w*\varrho_r$ on these sets. 
In particular, 
by classical results (see for instance
\cite[Theorem 2.39]{AmFuPa:00})
one has
\begin{equation}\label{99}
\begin{aligned}
\int_{V_t}\sqrt{1+|\nabla u_r|^2}dx
=&
\int_{V_t}\sqrt{1+|\nabla w_r|^2}dx\rightarrow 
\int_{V_t}\sqrt{1+|\nabla w|^2}dx+|Dw|(V_t)
\\
=&\int_{V_t}\sqrt{1+|\nabla u|^2}dx+2\pi\sum_{j\in J}|x_j-y_j|,
\end{aligned}
\end{equation} 
as $r\rightarrow0^+$. 
Concerning the integral over $\Om\setminus V_t$, using that $u_r$ takes values in $\Suno$ in $\Om\setminus B_{2r}(z_i)$, 
we can estimate
\begin{align}\label{100}
&\int_{\Om\setminus V_t}\sqrt{1+|\nabla u_r|^2+|\det(\nabla u_r)|^2}dx\leq \int_{\Om\setminus V_t}\sqrt{1+|\nabla u_r|^2}dx+\int_{B_{2r}(z_i)}|\det(\nabla u_r)|dx\nonumber\\
=&\int_{\Om\setminus V_t\setminus B_{2r}(z_i)}\sqrt{1+|\nabla u_r|^2}dx+\int_{B_{2r}(z_i)}\sqrt{1+|\nabla u_r|^2}dx+\int_{B_{2r}(z_i)}|\det(\nabla u_r)|dx.
\end{align}
 Let us estimate the last term; so fix $i\in\{1,\dots,k\}$ and assume without loss of generality that $\sigma_i=1$. In $B_{2r}(z_i)$ we then have 
$$u_r=(\cos(\theta_{z_i}),\sin(\theta_{z_i}))\psi(\rho_{z_i}),$$
where $(\rho_{z_i},\theta_{z_i})$ is a polar coordinate system around $z_i$. 
Thus
\begin{align}\label{101}
\nabla u_r(\rho_{z_i},\theta_{z_i})=\begin{pmatrix}
	\psi'(\rho_{z_i})\cos(\theta_{z_i})&-\frac{\psi(\rho_{z_i})}{\rho_{z_i}}\sin(\theta_{z_i})\\
	\psi'(\rho_{z_i})\sin(\theta_{z_i})&\frac{\psi(\rho_{z_i})}{\rho_{z_i}}\cos(\theta_{z_i})
\end{pmatrix}
\end{align}
and therefore $\det(\nabla u_r(\rho_{z_i},\theta_{z_i}))=\psi'(\rho_{z_i})\psi(\rho_{z_i})\rho_{z_i}^{-1}$, which gives
\begin{equation}\label{eq:contributo_polo}
	\int_{B_{2r}(z_i)}|\det(\nabla u_r)|dx=2\pi\int_{\frac43r}^{\frac53r}\psi(\rho)\psi'(\rho)d\rho=\pi.
	\end{equation}
Moreover \eqref{101} 
implies $|\nabla u_r|\leq \frac{C}{|x-z_i|}$ in $B_{2r}(z_i)$. 
In particular,
\begin{align*}
\int_{\cup_i B_{2r}(z_i)}\sqrt{1+|\nabla u_r|^2}dx\rightarrow 0\qquad\text{ as }r\rightarrow0^+.
\end{align*}
Finally, 
due to the fact that $\nabla u_r\rightarrow \nabla u$ 
in $L^1(\Om\setminus V_t; \R^{2 \times 2})$ we infer
\begin{equation}\label{eq:we_infer}
\int_{\Om\setminus V_t\setminus \cup_i B_{2r}(z_i)}
\sqrt{1+|\nabla u_r|^2}dx\rightarrow \int_{\Om\setminus V_t}\sqrt{1+|\nabla u|^2}dx.
\end{equation}
All in all, we have proved that the right-hand side in formula \eqref{100} 
tends
to
$\int_{\Om\setminus V_t}\sqrt{1+|\nabla u|^2}dx+k\pi$
 as $r\rightarrow0^+$. Thus, from \eqref{99} and  \eqref{100}
we get
\begin{align*}
\liminf_{r\rightarrow 0^+}\int_\Om\sqrt{1+|\nabla u_r|^2+|\det(\nabla u_r)|^2}dx\leq \int_{\Om}\sqrt{1+|\nabla u|^2}dx+k\pi+2\pi\sum_{j\in J}|x_j-y_j|,
\end{align*} 
which, in view 
of the results in Section \ref{sec:a_minimization_problem_for_atomic_distributions} 
concerning $\Vert \cdot\Vert_{{\rm flat},\alpha}$,
 gives \eqref{eq:ril_area_upper_estimate}. 
\end{proof}

We are now in a position to conclude the proof of Theorem \ref{teo:upper_bound_W11Suno} .

\begin{proof}[Proof of Theorem \ref{teo:upper_bound_W11Suno} ]
Let $u\in W^{1,1}(\Om;\Suno)$. 
In 
view of Proposition \ref{prop:density_in_flat_norm}
we can pick a sequence $(u_k)_k\subset W^{1,1}(\Om;\Suno)$ 
such that: 
\begin{itemize}
\item[(a)]
$\Det(\nabla u_k)=\pi\sum_{i=1}^{\indicefinalesommateorema_k}(\delta_{x^k_i}-\delta_{y^k_i})$ for some $N_k \in \mathbb N$, 
with
each $u_k$ satisfying (i) and (ii) of Proposition \ref{prop:density_in_flat_norm};
\item[(b)]
$\|u-u_k\|_{W^{1,1}}+\|\Det(\nabla u)-\Det(\nabla u_k)\|_{{\rm flat}}< \frac1k$
for all $k \in \NN$.
\end{itemize}
Hence, owing to (a), we are in a position to apply
 Theorem \ref{teo:main1} to each $u_k$, so that
\begin{align}\label{122_intro}
	\rilarea(u_k,\Om)\leq \int_\Om\sqrt{1+
|\nabla u_k|^2}dx+
\Vert \Det(\nabla u_k)\Vert_{{\rm flat}, 
\alpha} \qquad \forall k >0, 
\end{align}
and therefore
\begin{equation}\label{1234}
\begin{aligned}
\rilarea(u,\Om)&\leq 
\liminf_{k \to +\infty} \rilarea(u_k, \Om) 
\leq
\lim_{k\rightarrow +\infty}
\left(\int_\Om\sqrt{1+|\nabla u_k|^2}dx
+
\Vert \Det(\nabla u_k))\Vert_{{\rm flat}, \alpha}
\right)
\\
&=\int_\Om\sqrt{1+|\nabla u|^2}dx+
\Vert \Det(\nabla u)\Vert_{{\rm flat}, \alpha}.
\end{aligned}
\end{equation}
\end{proof}

\section{On the countably subadditive interior envelope of $\rilarea$}
\label{sec:on_the_countably_subadditive_interior_envelope}
As we have seen, the nonlocality of $\rilarea(u, \cdot)$ is
unavoidable. Therefore, it seems interesting to consider the 
largest countably subadditive set function non larger than
$\rilarea(u, \cdot)$, as defined in \eqref{eq:subadditive_interior_envelope}. 
 We have the following integral representation result:

\begin{Proposition}[\textbf{``Double'' relaxation}]\label{prop_2rel}
	Let $u\in W^{1,1}(\Om;\Suno)$. Then 
	\begin{align}\label{statement1}
		\rilrilarea(u,\Om)=\int_\Om\sqrt{1+|\nabla u|^2}dx.
	\end{align}
\end{Proposition}
\begin{proof}
	Since $\rilarea(u,A)\geq \int_A\sqrt{1+|\nabla u|^2}dx$ 
for any open set $A \subseteq \Om$, we only 
need to show the $\leq$ inequality in \eqref{statement1}.
	 
We know from Theorem \ref{teo:distributional_Jacobian_of_S1_valued_maps}
that $\Lambda := 
\frac1\pi\Det(\nabla u)=\sum_{i=1}^{+\infty}(\delta_{x_i}-\delta_{y_i})$ with $\sum_{i=1}^{+\infty}|x_i-y_i|<+\infty$. 
Fix $\eps>0$ and $\indicefinalesommaproposition_\eps \in \mathbb N$ so that $\sum_{i=\indicefinalesommaproposition_\eps+1}^{+\infty}|x_i-y_i|<\eps$. 
Set
$$\Lambda_\eps:=\sum_{i=\indicefinalesommaproposition_\eps+1}^{+\infty}(\delta_{x_i}-\delta_{y_i}).$$ 
Then, as $	\|\Lambda
_\eps\|_{{\rm flat}}<  \eps$, we infer
\begin{equation}\label{eq:flat_eps_flat_2eps}
\Vert \Lambda_\eps\Vert_{{\rm flat},\alpha}< 2\eps.
\end{equation}
Let $\{z_k: k=1,\dots,m_\eps\}$ be the set 
of points in $\{x_k,y_k:k\leq \indicefinalesommaproposition_\eps\}$ which are 
contained in $\Om$. 
Choose mutually 
disjoint discs $\overline B_{2r}(z_k)
\subset\Om$ for $k=1,\dots,m_\eps$, and  set
$$
\outofballs
:=\Om\setminus (\cup_{i=1}^{m_\eps}\overline B_{r}(z_k));
$$
notice that $\outofballs$ and $\overline B_{2r}(z_k)$ overlap on annuli
of radii $r$ and $2 r$.
 By definition of 
$\rilrilarea(u, \cdot)$ and using 
Theorem \ref{teo:upper_bound_W11Suno},  we have
\begin{equation}\label{main_est}
\begin{aligned}
\rilrilarea(u,\Om)\leq &
\rilarea(u,\outofballs)+
\sum_{k=1}^{m_\eps}\rilarea(u,\overline B_{2r}(z_k))
\leq 
\int_{\outofballs}\sqrt{1+|\nabla u|^2}dx
+\Vert \Det(\nabla u)\Vert_{{\rm flat},\alpha, \outofballs}
\\
& + \sum_{k=1}^{m_\eps}\left(
\int_{\overline B_{2r}(z_k)}\sqrt{1+|\nabla u|^2}dx
+\Vert \Det(\nabla u)\Vert_{{\rm flat},\alpha, 
\overline B_{2r}(z_k)}
\right).
\end{aligned}
\end{equation}
We claim that 
\begin{equation}\label{eq:wecla}
\Vert \Det(\nabla u)\Vert_{{\rm flat},\alpha,\outofballs}+
\sum_{k=1}^{m_\eps}
\Vert \Det(\nabla u)\Vert_{{\rm flat},\alpha,\overline B_{2r}(z_k)}
\leq 5\pi\eps,
\end{equation}
for $r>0$ sufficiently small.
Since the discs $\overline B_{2r}(z_k)$ 
are mutually disjoint for $k\geq1$, we 
see that 
\begin{align}\label{Teps}
\Vert \Lambda_\eps\Vert_{{\rm flat},\alpha,\cup_{k=1}^{m_\eps}
\overline B_{2r}(z_k)}=
\sum_{k=1}^{m_\eps}
\Vert \Lambda_\eps\res 
\overline B_{2r}(z_k)\Vert_{{\rm flat},\alpha,\overline B_{2r}(z_k)},
\end{align}
whereas, recalling 
\eqref{eq:flat_eps_flat_2eps},
\begin{align}\label{T_0}
\Vert \Lambda_\eps\Vert_{{\rm flat},\alpha,\outofballs}
\leq \Vert \Lambda_\eps\Vert_{{\rm flat},\alpha}<2\eps.
\end{align}
In $\overline B_{2r}(z_k)$, connecting $z_k$ 
to $\partial 
\overline B_{2r}(z_k)$ with a segment, we see that 
\begin{align*}
\Vert h_k\delta_{z_k}\Vert_{{\rm flat},\alpha}\leq 2|h_k|r,\qquad k=1,\dots,m_\eps,
\end{align*}
where $h_k\in\mathbb Z$ denotes the multiplicity of $z_k$ 
in the distribution $\Lambda-\Lambda_\eps=\sum_{k=1}^{\indicefinalesommaproposition_\eps}(\delta_{x_k}-\delta_{y_k})$.
On the other hand, by construction $(\Lambda-\Lambda_\eps)\res \overline B_{2r}(z_k)=h_k\delta_{z_k}$ and therefore
$$\Lambda\res \overline B_{2r}(z_k)=\Lambda_\eps\res \overline B_{2r}(z_k)+h_k\delta_{z_k},$$
which implies 
$$\Vert \Lambda\Vert_{{\rm flat},\alpha,\overline B_{2r}(z_k)}
\leq \Vert \Lambda_\eps\res 
\overline B_{2r}(z_k)\Vert_{{\rm flat},\alpha}+2|h_k|r_k.$$
Summing over $k=1,\dots, m_\eps$, by \eqref{Teps} one gets
\begin{align}
\sum_{k=1}^{m_\eps}
\Vert \Lambda\Vert_{{\rm flat},\alpha,\overline B_{2r}(z_k)}
&\leq \sum_{k=1}^{m_\eps}
\Vert \Lambda_\eps\res \overline B_{2r}(z_k)\Vert_{{\rm flat},
\alpha,\overline B_{2r}(z_k)}+\sum_{k=1}^{m_\eps}2|h_k|r\nonumber\\
&=\Vert \Lambda_\eps\Vert_{{\rm flat}, \alpha,\cup_{k=1}^{m_\eps}
\overline B_{2r}(z_k)}
+\sum_{k=1}^{m_\eps}2|h_k|r
\leq 2\eps+\sum_{k=1}^{m_\eps}2|h_k|r,\label{Mestimate}
\end{align}
where the last inequality follows from 
\eqref{eq:flat_eps_flat_2eps}, 
since 
$\Vert \Lambda_\eps\Vert_{{\rm flat},\alpha,
\cup_{k=1}^{m_\eps}\overline B_{2r}(z_k)}\leq \Vert \Lambda_\eps\Vert_{{\rm flat},\alpha}$.
{}From \eqref{T_0} 
we conclude 
$$
\Vert \Det(\nabla u)\Vert_{{\rm flat},\alpha,\outofballs}+
\sum_{k=1}^{m_\eps}
\Vert \Det(\nabla u)\Vert_{{\rm flat},
\alpha,\overline B_{2r}(z_k)}
\leq 4\pi\eps+\pi\sum_{k=1}^{m_\eps}2|h_k|r,$$
for any $r>0$ small enough, and \eqref{eq:wecla} follows.
 
Now, from \eqref{main_est}, 
we conclude that for every $r>0$ sufficiently small we have
\begin{align*}
	\rilrilarea(u,\Om)&\leq 
\int_{\outofballs}\sqrt{1+|\nabla u|^2}dx+
\sum_{k=1}^{m_\eps}\left(
\int_{\overline B_{2r}(z_k)}\sqrt{1+|\nabla u|^2}dx
\right)+5\pi\eps\nonumber\\
	&=\int_{\Om}\sqrt{1+|\nabla u|^2}dx+\sum_{k=1}^{m_\eps}\int_{B_{2r}(z_k)\setminus \overline B_{r}(z_k)}\sqrt{1+|\nabla u|^2}dx+5\pi\eps,
\end{align*}
which in turn (since $\grad u\in L^1(\Om; \R^{2\times 2})$) implies, 
letting $r\rightarrow0^+$,
$$\rilrilarea(u,\Om)\leq\int_{\Om}\sqrt{1+|\nabla u|^2}dx+5\pi\eps.$$
By the arbitrariness of $\eps$ we 
have $\rilrilarea(u,\Om)\leq\int_{\Om}\sqrt{1+|\nabla u|^2}dx$,
and this concludes the proof.
\end{proof}

As a direct consequence of
Proposition \ref{prop_2rel}, we readily get the following:
\begin{Corollary}[\textbf{Integral representation}]\label{cor:integral_representation}
Let $u\in W^{1,1}(\Om;\Suno)$. Then the set function $E\mapsto \rilrilarea(u,E)$ defines a Borel measure
absolutely continuous with respect to the Lebesgue measure $\mathcal L^2$ and coincides with
\begin{align*}
	\rilrilarea(u,E)=\int_E\sqrt{1+|\nabla u|^2}dx\qquad{\text{ for 
all Borel sets } E\subseteq\Om}.
\end{align*}
\end{Corollary}

In order to prove Theorem \ref{teo:doppio_rilassato_introduction} 
we now extend Theorem \ref{teo:main1} 
to the case of open sets 
obtained from $\Om$ by removing a finite set of points:

\begin{theorem}\label{teo:main2}
Let $C:=\{\punto_1,\dots,\punto_N\}$ be a finite set of distinct points of 
$\Om$. Let $u\in W^{1,1}(\Om;\Suno)$ satisfy (i) and (ii) of Theorem \ref{teo:main1}, and suppose that 
\begin{equation}\label{eq:p_i_x_k}
	\begin{cases}
		{\rm dist}(x_i, \partial \Om)
		\neq \vert x_i - \punto_k\vert 
		\\
		{\rm dist}(y_i, \partial \Om)
		\neq \vert y_i - \punto_k\vert 
	\end{cases}
	\qquad \forall k = 1,\dots, N, \ 
	\forall i=1,\dots,n.
\end{equation}
and that the points 
$x_i,y_i,\punto_k$, $i=1,\dots,n,\;k=1,\dots,N$, are  
three by three not collinear. Then
\begin{align*}
	\rilarea(u,\Om\setminus C)\leq \int_{\Om}\sqrt{1+|\nabla u|^2}dx+
\Vert {\rm Det}(\grad u)\Vert_{{\rm flat},\alpha,\Om\setminus C}.
\end{align*}
\end{theorem} 
\begin{proof}
Fix $\overline\rho>0$ small enough so that the discs $\overline B_{2\overline\rho}(\punto_k)$, 
$k=1,\dots,N$, are contained in $\Om$, are 
mutually disjoint and, for each $x_i$ or $y_i \in \Om$, we have 
$x_i, y_i \in \Om\setminus (\cup_{k=1}^N\overline B_{2\overline\rho}(\punto_k))$, respectively.
For all $\rho\in (0,2\overline\rho)$ let us denote $\Om_\rho:=\Om\setminus(\cup_{k=1}^N
\overline B_{\rho}(\punto_k))$, and let $R^\rho_{{\rm min}}\in \zerocurrfin(\Om_\rho)$ and $S^\rho_{{\rm min}}\in \mathcal S(\Om_\rho)$ be 
minimal currents 
as in Lemma \ref{lem:structure_of_minimizers_of_the_combinatorial_problem},
with $T = \frac{1}{\pi} {\rm Det}(\grad u)$. 
In particular, 
$R^\rho_{{\rm min}} = \sum_{i=1}^h \sigma_i \delta_{z_i}$ (notice that,
possibly reducing $\overline \rho>0$,  
$
R^\rho_{{\rm min}} = R_{{\rm min}}$ becomes independent of $\rho$), 
$S^\rho_{{\rm min}}=\sum_{j\in J}\jump{\overline{p_jq_j}}$ 
and there might be
points $p_j = p^\rho_j$ or $q_j = q^\rho_j$ on $\partial B_\rho(\punto_k)$  
for some $j$ and $k$. 
 However, since by assumption the points  
$x_i,y_i,\punto_k$ 
are three by three not collinear, it is easy to see that the points 
in $\{p_j,q_j\in \partial B_\rho(\punto_k),\text{for some }j,k\}$, if any, are
distinct. In particular, the segments  $\overline{p_jq_j}$, $j\in J$ are 
pairwise disjoint. Finally, as a consequence of \eqref{eq:p_i_x_k}, 
we may assume that if $\eta^\rho\in \partial B_\rho(\punto_k)$ is one of 
the points in the set $\{p_j,q_j\in \partial B_\rho(\punto_k) ~ \text{for some }j,k\}$, then 
 \begin{align}\label{slim}
 \eta^\rho\rightarrow \punto_k\qquad \text{ as }\;\rho\rightarrow 0^+.
 \end{align}
Using $R_{{\rm min}}$ and $S_{{\rm min}}^\rho$, we can now consider the 
sequence $(u^\rho_r)\subset 
C^1(\Om_\rho;\R^2)$ found in the proof of Theorem 
\ref{teo:main1} (see \eqref{eq:u_r}), with
$\Om$ replaced by $\Om_\rho$; 
so, for any $m\in \mathbb N$, we find $r_m \in (0, 2 \overline \rho)$ 
small enough so that $v^\rho_m:=u^\rho_{r_m}$ satisfies
\begin{equation}\label{eq:estimates_and_convergence}
\begin{aligned}
& \|v^\rho_m-u\|_{L^1(\Om_\rho;\R^2)}\leq \frac{1}{m},
\\
& \area(v^\rho_m,\Om_\rho)\leq \int_{\Om_\rho}\sqrt{1+|\nabla u|^2}dx+
\Vert
{\rm Det}(\grad u)\Vert_{{\rm flat},\alpha,\Om_\rho}+\frac1m. 
\end{aligned}
\end{equation}
%
Furthermore, as $\lim_{m\rightarrow+\infty} v^\rho_m = u$  in $W^{1,1}(\Om_\rho;\R^2)$, we may also suppose 
\begin{align}\label{est1}
\int_{\Om_{\rho}\setminus 
\overline \Om_{2\rho}}|\nabla v^\rho_m|dx\leq \int_{\Om_{\rho}\setminus \overline 
\Om_{2\rho}}|\nabla u|dx +\frac1m.
\end{align}
Now we suitably 
modify $v^\rho_m$ and extend it to $\Om\setminus C$: For any $k=1,\dots,N$, 
we pick a radius $s_k\in (\rho,2\rho)$ so that 
\begin{equation}\label{eq:sottomedia_unod}
\int_{\partial B_{s_k}(\punto_k)}|\nabla v^\rho_m|~d
\mathcal H^1\leq \frac{1}{\rho}\int_{B_{2\rho}(\punto_k)\setminus\overline B_{\rho}(\punto_k)}|\nabla v^\rho_m|dx.
\end{equation}
Then we define 
\begin{align}\label{est2}
\overline v^\rho_m(x):=\begin{cases}
	v_m^\rho(x)&\text{if }x\in \Om_{s_k}\\
	v_m^\rho(\punto_k+s_k\frac{x-\punto_k}{|x-\punto_k|})&\text{if }x\in B_{s_k}(\punto_k)\setminus \{\punto_k\}, \ \ k=1,\dots,N,
\end{cases}
\end{align}
where, with a little abuse of notation, we denoted $\Om_{s_k}:=\Om\setminus(\cup_{k=1}^N B_{s_k}(\punto_k))$.
From \eqref{est1} we get
$$
\int_{\Om_{s_k} \setminus \overline\Om_{2\rho}} \vert \grad \overline v_m^\rho\vert~dx =
\int_{\Om_{s_k} \setminus \overline\Om_{2\rho}} \vert \grad v_m^\rho\vert~dx \leq 
\int_{\Om_{\rho}\setminus \overline\Om_{2\rho}}|\nabla v_m^\rho|dx\leq 
\int_{\Om_{\rho}\setminus \overline\Om_{2\rho}}|\nabla u|dx + \frac{1}{m}.
$$
On the other hand, for any $k=1,\dots,N$,
$$
\begin{aligned}
 \int_{B_{s_k}(\punto_k)} \vert \grad \overline v^\rho_m\vert~dx =
\int_0^{s_k} 
\int_{\partial B_s(\punto_k)} \vert \grad \overline v_m^\rho\vert~d\mathcal
H^1~ds 
= s_k 
\int_{\partial B_{s_k}(\punto_k)} \vert \grad v_m^\rho\vert~d\mathcal
H^1, 
\end{aligned}
$$
where the last equality follows since, by definition, $\overline v^\rho_m$ is $0$-homogeneous in $B_{s_k}(c_k)$, and the integral $\int_{\partial B_s(\punto_k)} \vert \grad \overline v_m^\rho\vert~d\mathcal
H^1$ does not depend on $s\in (0,s_k)$.
Using \eqref{eq:sottomedia_unod}, and since $s_k\leq 2\rho$, it then follows
\begin{equation}\label{eq:estimate_grad_v_m}
\begin{aligned}	
\int_{\Om\setminus \overline\Om_{2\rho}}|\nabla \overline v^\rho_m|dx=\int_{\Om_{s_k} \setminus \overline\Om_{2\rho}} \vert \grad \overline v_m^\rho\vert~dx+\sum_{k=1}^N\int_{B_{s_k}(c_k)}|\nabla \overline v^\rho_m|dx
\leq 2\int_{\Om_{\rho}\setminus \overline\Om_{2\rho}}|\nabla u|dx +\frac{2}{m}.
\end{aligned}
\end{equation}
Furthermore, using that $|v_m^\rho|\leq1$, also $|\overline v_m^\rho|\leq1$ and it easily follows
\begin{align}\label{lim_uL1}
	\limsup_{m\rightarrow +\infty}\|\overline v_m^\rho-u\|_{L^1(\Om;\R^2)}\leq	\limsup_{m\rightarrow +\infty}\|\overline v_m^\rho-u\|_{L^1(\Om\setminus \Om_{2\rho};\R^2)}\leq 8\pi N\rho^2.
\end{align}
Since in $\Om_\rho\setminus \overline \Om_{2\rho}$ the map $v^\rho_m$ 
takes values in $\Suno$, using \eqref{eq:estimates_and_convergence}
and \eqref{eq:estimate_grad_v_m}
we can estimate 
\begin{equation}\label{final_est}
\begin{aligned}
\area(\overline v^\rho_m,\Om\setminus C)&\leq\area(\overline v^\rho_m;\Om\setminus (\cup_{k=1}^NB_{s_k}(\punto_k)))+\sum_{k=1}^N\pi s_k^2+2\int_{\Om_{\rho}\setminus \Om_{2\rho}}|\nabla u|dx+\frac2m
\\
&\leq \int_{\Om_\rho}\sqrt{1+|\nabla u|^2}dx+
\Vert {\rm Det}(\grad u)\Vert_{{\rm flat},\alpha,\Om_\rho}
+2\int_{\Om_{\rho}\setminus \Om_{2\rho}}|\nabla u|dx+\frac3m.
\end{aligned}
\end{equation}
We have proved that for every $\rho>0$ small enough 
we can find a Lipschitz map $\overline v^\rho_m :\Om\setminus C \to \R^2$ 
satisfying \eqref{final_est}. 
In particular, choosing a sequence $\rho_h\searrow0$, by a diagonal 
argument we find a sequence $(v_h)$ of Lipschitz maps\footnote{Even if 
$v_h$ is not $C^1$, by a density argument finding such
a sequence is sufficient (see \cite[Proposition 3.5]{BCS}).}
on $\Om\setminus C$ to $\R^2$ satisfying $v_h\rightarrow u$ in $L^1(\Om;\R^2)$ (by \eqref{lim_uL1}), and  
\begin{align}\label{final_est_bis}
	\area(v_h,\Om\setminus C)\leq 
\int_{\Om_{\rho_h}}\sqrt{1+|\nabla u|^2}dx+
\Vert {\rm Det}(\grad u)\Vert_{{\rm flat},\alpha,\Om_{\rho_h}}
+
2|Du|({\Om_{\rho_h}\setminus \overline \Om_{2\rho_h}})+\frac3h.
\end{align}
Letting $h\rightarrow +\infty$, we use that $ {\rm Det}(\grad u)$ is a measure so that it is easy to see, owing to \eqref{Mdef_loc},  
$$\Vert {\rm Det}(\grad u)\Vert_{{\rm flat},\alpha,\Om_{\rho_h}} \to 
\Vert {\rm Det}(\grad u)\Vert_{{\rm flat},\alpha,\Om\setminus C};$$ 
 then the thesis easily follows.
\end{proof}

Using Theorem \ref{teo:main2} and a density 
argument, we can prove the following:
\begin{theorem}\label{teo:doppio_rilassato}
Let $C$ be a finite set of distinct points of $\Om$,
and  $u\in W^{1,1}(\Om;\Suno)$. Then
\begin{align*}
	\rilarea(u,\Om\setminus C)\leq \int_{\Om}\sqrt{1+|\nabla u|^2}dx+
\Vert {\rm Det}(\grad u)\Vert_{{\rm flat},\alpha,\Om\setminus C}.
\end{align*}
\end{theorem}
\begin{proof}
It is sufficient to argue along the lines
 of the proof of Theorem \ref{teo:upper_bound_W11Suno},
 replacing, in \eqref{122_intro}, 
$\Vert \cdot\Vert_{{\rm flat},\alpha}$ by $\Vert \cdot\Vert_{{\rm 
flat},\alpha, \Om\setminus C}$. More specifically, 
in 
view of Proposition \ref{prop:density_in_flat_norm}
we can pick a sequence $(u_k)_k\subset W^{1,1}(\Om;\Suno)$ 
satisfying (a) and (b) of the proof
of Theorem \ref{teo:upper_bound_W11Suno}, and \eqref{eq:p_i_x_k}.
Applying
 Theorem \ref{teo:main2} to each $u_k$, we obtain
\begin{align}\label{122}
	\rilarea(u_k,\Om\setminus C)\leq \int_\Om\sqrt{1+
|\nabla u_k|^2}dx+
\Vert \Det(\nabla u_k)\Vert_{{\rm flat}, 
\alpha, \Om \setminus C} \qquad \forall k >0.
\end{align}
By lower-semicontinuity we conclude
\begin{equation}\label{123}
\begin{aligned}
\rilarea(u,\Om)&\leq 
\liminf_{k \to +\infty} \rilarea(u_k, \Om) 
\leq
\lim_{k\rightarrow +\infty}
\left(\int_\Om\sqrt{1+|\nabla u_k|^2}dx
+
\Vert \Det(\nabla u_k))\Vert_{{\rm flat}, \alpha, \Om \setminus C}
\right)
\\
&=\int_\Om\sqrt{1+|\nabla u|^2}dx+
\Vert \Det(\nabla u))\Vert_{{\rm flat}, \alpha},
\end{aligned}
\end{equation}
 where the equality is obtained since by 
(b), 
$u_k\rightarrow u$ in $W^{1,1}(\Om;\R^2)$ ,
$\Vert {\rm Det}(\grad u) - {\rm Det}(\grad u_k)\Vert_{{\rm flat}} \to 0$, 
and $\Vert \cdot\Vert_{{\rm flat},\alpha}$ is continuous in the flat metric. 
\end{proof}

\begin{theorem}\label{teo_chiuso}
Let $u\in W^{1,1}(\Om;\Suno)$. 
Then for every $\eps>0$ there exists a finite
set $C_\eps$ of points of 
$\Om$ 
such that 
\begin{equation}\label{eq:A_bar_C_eps}
	\rilarea(u,\Om\setminus C_\eps)\leq \int_\Om\sqrt{1+|\nabla u|^2}dx+\eps.
\end{equation} 
\end{theorem}

\begin{proof}
We know that $\frac1\pi\Det(\nabla u)=\sum_{i=1}^{+\infty}
(\delta_{x_i}-\delta_{y_i})$ 
with $\sum_{i=1}^{+\infty}
|{x_i}-{y_i}|<+\infty$. Take $N_\eps\in \mathbb N$ 
so that  $\sum_{i=N_\eps+1}^{+\infty}|x_i-y_i|<\frac{\eps}{2\pi}$, and 
let  $C_\eps
:=
\{x_k \in \Om: 1 \leq k\leq N_\eps\}\cup
\{y_k \in \Om: 1\leq k\leq N_\eps\}$.
Set $T:=\Det(\nabla u)\res(\Om\setminus C_\eps)$,
so that $\Vert T\Vert_{{\rm flat},\alpha,\Om\setminus C_\eps}
\leq \eps$. 
Then Theorem \ref{teo:doppio_rilassato} implies \eqref{eq:A_bar_C_eps}.
\end{proof}

Using Theorem \ref{teo_chiuso} we can positively answer 
to a modification of a
conjecture by De Giorgi \cite{DG}, adapted to the context of 
$\Suno$-valued Sobolev maps.
\begin{Corollary}\label{cor:De_Giorgi}
	Let $u\in W^{1,1}(\Om;\Suno)$. 
Then
	$$
	\rilrilarea (u, \Om)=
\inf_{\substack{C\subset\Om, \mathcal H^0(C)< +\infty}}
	\rilarea(u,\Om\setminus C).
$$
\end{Corollary}
\begin{proof}
{}From Theorem \ref{teo_chiuso} we get
$$
\inf_{\substack{C\subset\Om, \mathcal H^0(C)<+\infty}}
	\rilarea(u,\Om\setminus C) 
\leq \int_\Om \sqrt{1+\vert \grad u\vert^2}~dx=\rilrilarea(u,\Om),
$$
where the last equality follows from Corollary \ref{cor:integral_representation}.
On the other hand, for any finite set $C \subseteq \Om$, 
we know \cite{ADM} that 
$$
	\rilarea(u,\Om\setminus C) \geq  
\int_{\Om\setminus C} \sqrt{1+\vert \grad u\vert^2}~dx
=\int_{\Om} \sqrt{1+\vert \grad u\vert^2}~dx
=\rilrilarea(u,\Om).
$$
\end{proof} 

\section{Appendix}\label{sec:appendix}
In this appendix we collect a first 
standard result, and a proposition with an
independent interest.

\begin{lemma}\label{fed}
Let $\Lambda\in {\rm Lip}_0(\Om)'$. Then 
\begin{equation}\label{eq:inf_equal_inf}
\begin{aligned}
 \sup_{{\substack{\varphi\in C_c^1(\Om)\\\|\varphi\|_{{\rm Lip}_0, \alpha}\leq 1}}}\langle\Lambda,\varphi\rangle
 =&\inf
\Big\{|R|_\Om+ \alpha^{-1}|S|_\Om:
(R,S)\in \mathcal D_0(\Om) \times  \mathcal D_1(\Om),\;\Lambda=
R+\partial S\Big\}
\\
=& 
 \sup_{{\substack{\varphi\in {\rm Lip}_0(\Om)\\\|\varphi\|_{{\rm Lip}_0, \alpha}\leq 1}}}\langle\Lambda,\varphi\rangle.
\end{aligned}
\end{equation}
\end{lemma}

\begin{proof}
We adapt the arguments 
of \cite[page 367]{Federer}. 
Let $R\in \mathcal D_0(\Om)$ and $S\in \mathcal D_1(\Om)$ 
be such that $\Lambda=R+\partial S$ in $\mathcal D_0(\Om)$. 
Then, as $\alpha=\frac12$, 
$$
\left|\langle\Lambda,\varphi\rangle\right|\leq |R(\varphi)|+
|S(d\varphi)|\leq (|R|_\Om+2|S|_\Om)\|\varphi\|_{{\rm Lip}_0, \alpha}
\qquad
\forall \varphi\in C^1_c(\Om),
$$
so the inequality
$\le$ holds in the first line of 
\eqref{eq:inf_equal_inf}.
To prove the converse inequality, set
$$Y:=\left
\{(\varphi,\psi)\in C^1_c(\Om)\times C^0_c(\Om;\R^2)\right\}$$ 
endowed with the norm $\|(\varphi,\psi)\|_Y:=\max\{\|\varphi\|_{L^\infty},
\frac{1}{2}\|\psi\|_{L^\infty}\}$, 
and define the linear injective operator
$$
Q:C^1_c(\Om)\rightarrow Y, \qquad
Q(\varphi):=(\varphi,\nabla\varphi)\qquad \forall \varphi\in C^1_c(\Om).$$
Since
$Q(C^1_c(\Om)) \subset Y$, we have 
$$
\langle \Lambda,Q^{-1}(\varphi, \nabla \varphi)\rangle
\leq \Vert \Lambda\Vert_{{\rm flat},\alpha} \|(\varphi, \nabla \varphi)\|_Y
\qquad \forall \varphi\in C^1_c(\Om),$$
and therefore
we can extend the linear functional $\Lambda\circ Q^{-1}:
Q(C^1_c(\Om))
\rightarrow \R$ to some linear functional $L:Y\rightarrow \R$ with 
\begin{equation}\label{eq:L_below_norm}
L(\varphi,\psi)\leq \Vert \Lambda\Vert_{{\rm flat},\alpha}\|(\varphi, \psi)\|_Y
\qquad \forall (\varphi,\psi)\in Y.
\end{equation} 
Now we define 
\begin{align*}
&R(\varphi):=L(\varphi,0)\qquad \forall \varphi\in C^1_c(\Om),\\
&S(\psi):=L(0,\psi)\qquad \forall \psi\in C^0_c(\Om;\R^2),
\end{align*}
so that, from \eqref{eq:L_below_norm},
 $\|(\varphi,\psi)\|_Y\leq 1$ implies $R(\varphi)+
S(\psi)=L(\varphi,\psi)\leq \Vert \Lambda\Vert_{{\rm flat},\alpha}$. 
In particular, $R\in 
\mathcal D_0(\Om)$, $S\in \mathcal D_1(\Om)$, and passing to the supremum,
$$|R|_\Om+2|S|_\Om\leq\Vert \Lambda\Vert_{{\rm flat},\alpha}.
 $$
Since $R(\varphi)+S(d\varphi)=L(\varphi,\nabla \varphi)=\langle\Lambda,\varphi\rangle$ for all $ \varphi\in C^1_c(\Om)$,
it follows $\Lambda=R+\partial S$, and 
the first equality in \eqref{eq:inf_equal_inf} follows.

 To show the second equality 
in \eqref{eq:inf_equal_inf}, we first observe that, from
the first equality, 
$$\sup_{{\substack{\varphi\in C_c^1(\Om)\\\|\varphi\|_{{\rm Lip_0}, \alpha}
\leq 1}}}\langle\Lambda,\varphi\rangle\leq 
\Vert \Lambda\Vert_{{\rm flat},\alpha},$$
and so $\leq$ holds. 
On the other hand, if $R\in \mathcal \mathcal R_f$ and $S\in \mathcal S$ 
are such that $\Lambda=R+\partial S$ in $\mathcal D_0(\Om)$, then 
$$\left|\langle\Lambda,\varphi\rangle\right|=|R(\varphi)|+
|S(d\varphi)|\leq (|R|_\Om+2|S|_\Om)\|\varphi\|_{{\rm Lip}_0, \alpha}
\qquad
\forall \varphi\in \text{{\rm Lip}}_0(\Om),
$$
so also the inequality
$\le$ holds, thanks to Corollary \ref{cor3.3}. 
\end{proof}

The next result has been used in the proof of Proposition
\ref{lem_min_finita}, and is based on Lemma 
\ref{lem:structure_of_minimizers_of_the_combinatorial_problem}.

\begin{Proposition}\label{minimal_lip}
Let $T=\sum_{i=1}^n(\delta_{x_i}-\delta_{y_i})\in \mathcal D'(\Om)$ be as in \eqref{T_finite2} and 
satisfying $({\rm P_f})$, 
where $x_i,y_i\in \overline\Om$, $x_i \neq y_i$. Let $I_P$, $I_D$, $\tau$,
 $R_{{\rm min}}$ and $S_{{\rm min}}$ be 
as in \eqref{defRSprimo_bis}. Then 
	\begin{align*}
\exists \varphi\in \text{{\rm Lip}}_0(\Om) ~~ {\rm with~}  
~ \| \varphi\|_{{\rm Lip}_0,\alpha}\leq 1 ~~{\rm such~that}~~ 
		\langle T, \varphi\rangle
=|R_{\rm min}|_\Om+ \alpha^{-1}|S_{\rm min}|_\Om.
	\end{align*}
As a consequence, 
for all $k\in I_P$ and $j\in I\setminus \tau(I_D)$ with $x_k\in \Om$ 
and $y_j\in \Om$, we have $\varphi(x_k)=-\varphi(y_j)=1$, 
and for all $k\in I_D$ we have $\varphi(x_k)-\varphi(y_{\tau(k)})=\alpha^{-1}
|x_k-y_{\tau(k)}|$.
In particular, 
$$\min\{|R|_\Om+\alpha^{-1}|S|_\Om:R\in \zerocurrfin,\;S\in {\mathcal S},\;T=R+\partial S\}=\max_{\substack{\phi\in {\rm Lip}_0(\Om)\\
		\|\phi\|_{{\rm Lip}_0,\alpha}\leq 1}}\langle T,\phi\rangle.$$
\end{Proposition}

\begin{proof}
Define 
$P^+:=\{k\in I_P:x_k\in \Om\}$ and
$P^-:=\{k\in I\setminus \tau(I_D):y_k\in \Om\}.$
	The function $\varphi$ in the statement must satisfy $\varphi(x_k)=1$ for all $k\in P^+$ and $\varphi(y_k)=-1$ for all $k\in P^-$, and, recalling that $\alpha^{-1}=2$, also $\varphi(x_k)-\varphi(y_{\tau(k)})=2|x_k-y_{\tau(k)}|$ for all $k\in I_D$.
	
	For any $k\in P^+$, we define $\phi_k(x)
:=1-2\vert x-x_k\vert$ for all $x\in \overline\Om$, and set 
$$
\phi
:=\begin{cases}
\max_{k\in P^+}\{\phi_k\} & {\rm if}~ P^+\neq \emptyset,
\\
-1 & {\rm if}~  P^+ = \emptyset.
\end{cases}
$$
Define also
	$\Psi_0
:=\max\{\phi_0,\phi\},$
where $\phi_0(x) :=\max\{-1,-2d(x,\partial\Om)\}$
for all $x\in \overline\Om$, i.e., 
\begin{align}\label{defpsi0}
	\Psi_0(x)=\max\{-1,-2d(x,\partial\Om),1-2|x-x_k|,k\in P^+\}\qquad \forall x\in \overline\Om,
\end{align}
and observe that 
$\Psi_0
\in \text{{\rm Lip}}_0(\Om)$
with 
$\| \Psi_0\|_{{\rm Lip}_0,\alpha}\leq 1$. 

	Using the minimality, and in particular \eqref{min_2}, one
verifies that $\Psi_0(x_k)=1$ for all $k\in P^+$.
Let us check that $\Psi_0(y_k)=-1$ for all $k\in P^-$. If not, either $\phi_0(y_k)>-1$ or $\phi(y_k)>-1$,
and both the two cases are excluded again by \eqref{min_2}.
	
Now, we have to take into account the dipoles and the boundary values of $\Psi_0$.
We divide the proof into three steps.
	
	\textit{Step 1:} For all $k\in I_D$ with $y_{\tau(k)},x_k\in \Om$, 
	\begin{align}\label{prop_Phi0}
		\Psi_0(x_k)\leq 1, \qquad \qquad \Psi_0(y_{\tau(k)})\leq 1-2|x_k-y_{\tau(k)}|.
	\end{align}
	Moreover, if 
either  $y_{\tau(k)}\in \partial\Om$ or $x_k\in \partial\Om$
for some $k\in I_D$ , then
	\begin{align}\label{prop_Phi0_bis}
		\text{either }\Psi_0(x_k)\leq 2 d(x_k,\partial\Om) \quad 
\text{or } \quad \Psi_0(y_{\tau(k)})=-2 d(y_{\tau(k)},\partial\Om),
	\end{align}
	respectively. 

Let us check \eqref{prop_Phi0}.  
The first inequality follows since $\Psi_0\leq 1$ on $\overline\Om$.
The second inequality is deduced as follows: By 
\eqref{defpsi0}, if $\Psi_0(y_{\tau(k)})=-2 d(x_k,\partial\Om)$, then we conclude by \eqref{min_3}. If instead $\Psi_0(y_{\tau(k)})=1-2|x_h-y_{\tau(k)}|$ for some $h\in P^+$, then we conclude 
since by minimality
$|x_h-y_{\tau(k)}|\geq |x_k-y_{\tau(k)}|$ (where we have used that $h\in P^+$).
	
	Let us now check \eqref{prop_Phi0_bis}. Assume $y_{\tau(k)}\in \partial\Om$ and, by \eqref{defpsi0}, that $\Psi_0(x_k)=1-2|x_h-x_k|$ for some $h\in P^+$. By 
the triangle inequality and \eqref{min_2}, we have
	$$\Psi_0(x_k)\leq 1+2|x_k-y_{\tau(k)}|-2|x_h-y_{\tau(k)}|\leq2|x_k-y_{\tau(k)}|=2 d(x_k,\partial\Om).$$ 
	Assume instead $x_k\in \partial\Om$. If, by contradiction, $\Psi_0(y_{\tau(k)})>-2\dist(y_{\tau(k)},\partial\Om)=-2|y_{\tau(k)}-x_k|$, for some $h\in P^+$ we necessarily have $\Psi_0(y_{\tau(k)})=1-2|x_h-y_{\tau(k)}|>-2|y_{\tau(k)}-x_k|$. This contradicts the minimality of $R_{{\rm min}}$ and  $S_{{\rm min}}$, because 
a
direct check shows that 
	$$|R'|_\Om+2|S'|_\Om\leq |R_{{\rm min}}|_\Om+2|S_{{\rm min}}|_\Om-1+2|x_h-y_{\tau(k)}|-2|x_k-y_{\tau(k)}|<|R_{{\rm min}}|_\Om+2|S'_{{\rm min}}|_\Om,$$
where $R':=R_{{\rm min}}-\delta_{x_h}$
and $S':=S_{{\rm min}}-\jump{\overline{y_{\tau(k)}x_k} }+\jump{\overline{y_{\tau(k)}x_h}}$.

	Eventually we check that 
	\begin{align}\label{prop_Phi0_tris}
		\Psi_0=0\qquad {\rm on}~ \partial\Om.
	\end{align}
	Indeed, if $\Psi_0(x)=1-2|x_k-x|>0$ for some $k\in P^+$, then arguing as before we can define $R'=R_{{\rm min}}-\delta_{x_k}$ and $S'=S_{{\rm min}}+\jump{\overline{xx_k}}$, and an easy check shows that $|R'|_\Om+2|S'|_\Om<|R_{{\rm min}}|_\Om+2|S_{{\rm min}}|_\Om$, against the minimality.
\medskip

Before proceeding to the next step, for the sake of simplicity and without loss of generality, we 
relabel the indices and assume that $\tau:I_D\rightarrow I_D$ is the identity map, so that  $I_P=I\setminus I_D=I\setminus \tau(I_D)$. 
The function $\varphi$, that will be constructed 
starting from $\Psi_0$,
 should satisfy $\varphi(x_k)=\Psi_0(x_k)=1$ for all $k\in P^+$, $\varphi(y_k)=\Psi_0(y_k)=-1$ for all $k\in P^-$, 
and 
$\varphi(x_k)-\varphi(y_k)=2|x_k-y_k|$
for every $k\in I_D$. 
To build such a $\varphi$, we apply 
a recursive procedure.  We define, for all $m\geq1$, the function $\Psi_m$ as follows:
	\begin{align}
		\Psi_m(x) :=\max\{\Psi_{m-1}(x),\Phi_{m}(x)\}
\qquad \forall x\in \overline\Om,\label{defpsim}
	\end{align}
	where $\Phi_m$ is given by
	$$\Phi_m(x):=\max_{k\in I_D}\{\phi_k^m(x)\},\qquad \phi_k^m(x):=\Psi_{m-1}(n_{k})+2|x_k-n_{k}|-2|x_k-x|,\qquad \forall x\in \overline\Om.$$
	Trivially, $\Phi_m$ is Lipschitz continuous with Lipschitz constant $2$, for all $m\geq0$.

\textit{Step 2:} 
We claim that 
	\begin{itemize}
		\item[(a)]  $\Psi_m(y_k)\leq 1-2|x_k-y_k|$ and $\Psi_m(x_k)\leq 1$ for all $k\in  I_D$ with $x_k,y_k\in \Om$;
		\item[(b)] $\Psi_m(x_k)=1$ and $\Psi_m(y_h)=-1$ for all $k\in P^+$, $h\in P^-$;
		\item[(c)]  $\Psi_m(x)=0$ for $x\in \partial\Om$;
		\item[(d)] if, for some $k\in I_D$, either  $y_{k}\in \partial\Om$ or $x_k\in \partial\Om$, then
		\begin{align}\label{prop_m}
			\text{either }\Psi_m(x_k)\leq 2\dist(x_k,\partial\Om),\qquad\qquad \text{or }\Psi_m(y_k)=-2\dist(y_k,\partial\Om),
		\end{align}
		respectively.
	\end{itemize}
	Proof of (a): First we notice that, if $\{h_1,\dots,h_j\}\subseteq I_D$ with $h_i\neq h_{i'}$ for $i\neq i'$, then the minimality of $S_{{\rm min}}$ implies that 
	\begin{align}\label{min_4}
		\sum_{i=1}^j|x_{h_i}-y_{h_i}|\leq |x_{h_1}-y_{h_2}|+|x_{h_2}-y_{h_3}|+\dots+|x_{h_j}-y_{h_1}|.
	\end{align}
	Now, by construction and definition of $\Phi_m$, we can find a set of $r\geq0$ indices (possibly $r=0$) $0= m_1<\dots<m_r\leq m$, and indices $k_1,\dots,k_r\in I_D$ such that 
	\begin{align}\label{sequence(c)}
		&\Psi_m(y_k)=\phi^{m}_{k}(y_k)=\Psi_{m_r}(y_{k_r})+2|x_{k_r}-y_{k_r}|-2|x_{k_r}-y_{k}|,\nonumber\\	&\Psi_{m_r}(y_{k_r})=\phi^{m_{r}}_{k_{r}}(y_{k_r})=\Psi_{m_{r-1}}(y_{k_{r-1}})+2|x_{k_{r-1}}-y_{k_{r-1}}|-2|x_{k_{r-1}}-y_{k_r}|,\nonumber\\
		&\qquad\dots\qquad\\\nonumber
		&\Psi_{m_3}(y_{k_3})=\phi^{m_{3}}_{k_{3}}(n_{k_{3}})=\Psi_{m_{2}}(y_{k_{2}})+2|x_{k_{2}}-y_{k_{2}}|-2|x_{k_{2}}-y_{k_3}|,\\
		&\Psi_{m_2}(y_{k_2})=\phi^{m_{2}}_{k_{2}}(y_{k_{2}})=\Psi_{0}(y_{k_{1}})+2|x_{k_{1}}-y_{k_{1}}|-2|x_{k_{1}}-y_{k_2}|.\nonumber
	\end{align}
	Notice that, if $r=0$ we simply have $\Phi_m(y_k)=\Phi_0(y_k)$ and (a) follows from \eqref{prop_Phi0} and thanks to the fact that 
$\Psi_0\in \text{{\rm Lip}}_0(\Om)$,
$\| \Psi_0\|_{{\rm Lip}_0,\alpha}\leq 1$.
Assume then that $r>0$, and so from \eqref{sequence(c)} it follows  that 
	$$\Psi_m(y_k)=\Psi_0(y_{k_1})+2\sum_{i=1}^r|x_{k_i}-y_{k_i}|-2\sum_{i=1}^r|x_{k_i}-y_{k_{i+1}}|,$$
	where we have set $y_{k_{r+1}}=y_k$. 
Now we have two cases:
\begin{itemize}
\item[(a1)]$\Psi_0(y_{k_1})=-2\dist(\partial\Om,y_{k_1})$;
\item[(a2)] $\Psi_0(y_{k_1})=1-2|x_h-y_{k_1}|$ for some $h\in P^+$.
\end{itemize}
In case (a1) we will show that 
	\begin{equation}\label{caso1(c)}
		2\sum_{i=1}^r|x_{k_i}-y_{k_i}|-2\sum_{i=1}^r|x_{k_i}-y_{k_{i+1}}|\leq 1+2\dist(\partial\Om,y_{k_1})-2|x_k-y_k|,
	\end{equation}
whereas in case (a2)  we will show that 
	\begin{equation}\label{caso2(c)}
		2\sum_{i=1}^r|x_{k_i}-y_{k_i}|-2\sum_{i=1}^r|x_{k_i}-y_{k_{i+1}}|\leq 2|x_h-y_{k_1}|-2|x_k-y_k|,
	\end{equation}
	and this will conclude (a). First, in view of \eqref{min_3}, we can assume that the indices $k_1,k_2,\dots,k_r$ are all distinct. Indeed, if  $k_{i}=k_{i'}$ for some $i\neq i'$, then we can  erase the indices $k_i,k_{i+1},\dots,k_{i'-1}$, since $\sum_{j=i}^{i'-1}|x_{k_j}-y_{k_j}|-\sum_{j=i}^{i'-1}|x_{k_j}-y_{k_{j+1}}|\leq0$ by \eqref{min_3}.
	
	Therefore, assuming (a1), let us prove \eqref{caso1(c)}. Consider the point $p$ on $\partial\Om$ so that $|p-y_{k_1}|=\dist(\partial\Om,y_{k_1})$; then \eqref{caso1(c)} is a consequence of the minimality of $R_{{\rm min}}$ and $S_{{\rm min}}$. Indeed, setting $$R'=R_{{\rm min}}+\delta_{x_k}\qquad S'=S_{{\rm min}}-\sum_{i=1}^{r+1}\jump{\overline{x_{k_i}y_{k_i}}}+\sum_{i=1}^{r}\jump{\overline{x_{k_i}y_{k_{i+1}}}}+\jump{\overline{py_{k_1}}},$$
	the inequality $|R'|_\Om+2|S'|_\Om\geq |R_{{\rm min}}|_\Om+2|S_{{\rm min}}|_\Om$ is equivalent to \eqref{caso1(c)}.
	
	In case (a2), instead, to get \eqref{caso2(c)}, arguing as before, it suffices to write $|R'|_\Om+2|S'|_\Om\geq |R_{{\rm min}}|_\Om+2|S_{{\rm min}}|_\Om$, with
	$$R'=R_{{\rm min}}+\delta_{y_{k_1}}+\delta_{x_h}\qquad S'=S_{{\rm min}}-\sum_{i=1}^{r+1}\jump{\overline{x_{k_i}y_{k_i}}}+\sum_{i=1}^{r}\jump{\overline{x_{k_i}y_{k_{i+1}}}}+\jump{\overline{x_{h}y_{k_1}}}.$$
	So far we have proved that $\Psi_m(y_k)\leq 1-2|x_k-y_k|$, for all $k\in I_D$. The fact that $\Psi_m(x_k)\leq 1$ for all $k\in I_D$ 
follows thanks to the fact that
$\Psi_m
\in \text{{\rm Lip}}_0(\Om)$,  
with 
$\| \Psi_m\|_{{\rm Lip}_0,\alpha}\leq 1$.

	Proof of (b): Let us check that $\Psi(p_h)\leq 1$ for all $h\in P^+$. To this aim it is 
sufficient to observe that $\Phi_m(x)\leq 1 $ for all $x\in \overline\Om$, since $\phi_k^m(x)\leq \phi_k^m(x_k)=\Psi_{m-1}(y_k)+2|x_k-y_k|\leq 1$, for all $k\in I_D$, thanks to point (a).
	
	Let us check that $\Psi(x_h)=-1$ for all $h\in P^-$. Arguing as in \eqref{sequence(c)}, we can find a set of $r\geq0$ indices $0= m_1<\dots<m_r\leq m$, and indices $k_1,\dots,k_r\in I_D$ such that
	\begin{align*}
		\Psi_m(y_h)=\Psi_0(y_{k_1})+2\sum_{i=1}^r|x_{k_i}-y_{k_i}|-2\sum_{i=1}^r|x_{k_i}-y_{k_{i+1}}|,
	\end{align*} 
	where $k_{r+1}:=h$. If $r=0$, we readily conclude that $\Psi_m(y_h)=\Psi_0(y_{h})$ and the thesis follows from the fact that $\Psi_0(y_{h})=-1$ for all $h\in P^-$. Hence assume $r>0$. Also here we have to cases: 
	\begin{itemize}
		\item[(b1)] $\Psi_0(y_{k_1})=1-2|x_j-y_{k_1}|$ for some $j\in P^+$;
		\item[(b2)]$\Psi_0(y_{k_1})=-2\dist(y_{k_1},\partial\Om)$.
		\end{itemize}
	In the first case the thesis is equivalent to
	\begin{equation}\label{caso1(b)}
		1-2|x_j-y_{k_1}|+2\sum_{i=1}^r|x_{k_i}-y_{k_i}|-2\sum_{i=1}^r|x_{k_i}-y_{k_{i+1}}|\leq -1.
	\end{equation}
	As before we might assume that the indices $k_i$ 
are distinct.
	Writing $R':=R_{{\rm min}}+\delta_{x_h}-\delta_{y_j}$ and $S':=S_{{\rm min}}-\sum_{i=1}^r\jump{\overline{x_{k_i}y_{k_i}}}+\sum_{i=1}^r\jump{\overline{x_{k_i}y_{k_{i+1}}}}+\jump{\overline{x_{j}y_{k_1}}}$, \eqref{caso1(b)} readily follows by the inequality $|R'|_\Om+2|S'|_\Om\geq |R_{{\rm min}}|_\Om+2|S_{{\rm min}}|_\Om$.
	
	Now, if (b2) holds, we will conclude by showing that 
	\begin{equation}
		-2\dist(y_{k_1},\partial\Om)+2\sum_{i=1}^r|x_{k_i}-y_{k_i}|-2\sum_{i=1}^r|x_{k_i}-y_{k_{i+1}}|\leq -1.
	\end{equation}
	Also this is obtained using the minimality of $R_{{\rm min}}$ and $S_{{\rm min}}$, by setting $R':=R_{{\rm min}}+\delta_{y_h}$ and $S':=S_{{\rm min}}-\sum_{i=1}^r\jump{\overline{x_{k_i}y_{k_i}}}+\sum_{i=1}^r\jump{\overline{x_{k_i}y_{k_{i+1}}}}+\jump{\overline{x_jy_{k_1}}}$ where $x_j\in \partial\Om$ is such that $|x_j-y_{k_1}|=\dist(y_{k_1},\partial\Om)$.
	
	Proof of (c): To show this, fix $x\in \partial\Om$; if $\Psi_m(x)=\Psi_0(x)=0$ there is nothing to prove. If not, we can find a set of $r>0$ indices $0= m_1<\dots<m_r\leq m$, and indices $k_1,\dots,k_r\in I_D$ such that
	\begin{align*}
		\Psi_m(x)=\Psi_0(y_{k_1})+2\sum_{i=1}^r|x_{k_i}-y_{k_i}|-2\sum_{i=1}^r|x_{k_i}-y_{k_{i+1}}|,
	\end{align*} 
	where $y_{k_{r+1}}:=x$. If $\Psi_0(y_{k_1})=-2\dist(y_{k_1},\partial\Om)$, we show that 
	$$\sum_{i=1}^r|x_{k_i}-y_{k_i}|\leq \sum_{i=1}^r|x_{k_i}-y_{k_{i+1}}|+\dist(y_{k_1},\partial\Om).$$
	As usual we might assume that the indices $k_i$ are distinct; being $x\in \partial\Om$, we have $|x_{k_r}-x|\geq \dist(x_{r_k},\partial\Om)$, and so the previous inequality is obtained by minimality of $R_{{\rm min}}$ and $S_{{\rm min}}$, arguing similarly as in the preceding cases.
	
	If instead $\Psi_0(y_{k_1})=1-2|x_h-y_{k_1}|$ for some $h\in P^+$, we reduce ourselves to prove that 
	$$1+\sum_{i=1}^r|x_{k_i}-y_{k_i}|\leq \sum_{i=1}^r|x_{k_i}-y_{k_{i+1}}|+|x_h-y_{k_1}|,$$
	which is again implied by the minimality of $R_{{\rm min}}$ and $S_{{\rm min}}$.
	
	Proof of (d): 
	The first condition in \eqref{prop_m} is a consequence of point (c) 
and the fact that $\Psi_m 
\in \text{{\rm Lip}}_0(\Om)$,  
with 
$\| \Psi_m\|_{{\rm Lip}_0,\alpha}\leq 1$. Let us prove the second condition. If $\Psi_m(y_k)=\Psi_0(y_k)$ then the thesis follows from \eqref{prop_Phi0_bis}; if not, we can find a sequence of $r>0$ indices $0= m_1<\dots<m_r\leq m$, and indices $k_1,\dots,k_r\in I_D$ such that
	\begin{align*}
		\Psi_m(y_k)=\Psi_0(y_{k_1})+2\sum_{i=1}^r|x_{k_i}-y_{k_i}|-2\sum_{i=1}^r|x_{k_i}-y_{k_{i+1}}|,
	\end{align*} 
	where $k_{r+1}:=k$. Now, either $\Psi_0(y_{k_1})=-2\dist(\partial\Om,y_{k_1})$ or $\Psi_0(y_{k_1})=1-2|x_h-y_{k_1}|$ for some $h\in P^+$. 
Again assuming that the indices $k_i$ are distinct, in the first case it is 
sufficient to observe that 
	$$\dist(\partial\Om,y_{k})+\sum_{i=1}^r|x_{k_i}-y_{k_i}|\leq \dist(\partial\Om,n_{k_1})+\sum_{i=1}^r|x_{k_i}-y_{k_{i+1}}|,$$
	which follows once more by the minimality of $R_{{\rm min}}$ and $S_{{\rm min}}$, since $\dist(\partial\Om,y_{k})=|x_k-y_k|$, and $x_k\in \partial\Om$.
	
	Assuming instead that $\Psi_0(y_{k_1})=1-2|x_h-y_{k_1}|$ for some $h\in P^+$, we can show that 
	$$1+\dist(\partial\Om,y_{k})+\sum_{i=1}^r|x_{k_i}-y_{k_i}|\leq \sum_{i=1}^r|x_{k_i}-y_{k_{i+1}}|+|x_h-y_{k_1}|,$$
	using once again the minimality of $R_{{\rm min}}$ and $S_{{\rm min}}$.
		
Before passing to the next step, observe that 
since $\Psi_{m}\leq \Psi_{m+1}$ for all $m\geq0$, we can take the limit 
	$$\varphi:=\lim_{m\rightarrow +\infty}\Psi_m,$$
	and thanks to the properties of $\Psi_m$ we easily infer that conditions (a), (b), (c), and (d) are still valid for $\varphi$. Furthermore, since $\Psi_m$ is Lipschitz continuous with Lipschitz constant $2$, we also have that $\Psi_m\rightarrow\varphi$ uniformly on $\overline\Om$.
\medskip

	The next step concludes the proof of the lemma. 

	\textit{Step 3.} The function $\varphi$ satisfies
	\begin{align}\label{claim}
		\varphi(x_k)=\varphi(n_k)+2|x_k-y_k|\qquad \forall k\in I_D.
	\end{align}
To see this we define
	\begin{align*}
		\overline\varphi
:=\max_{k\in I_D}\{\phi_k\},\qquad \phi_k(x):=\varphi(y_{k})+2|x_k-y_{k}|-2|x_k-x|,\qquad \forall k\in I_D,\;\forall x\in \overline\Om.
	\end{align*}
	In order to prove 
\eqref{claim} we show that $\varphi\geq \overline\varphi$; 
	this implies
 that $\varphi(x_k)\geq \phi_k(x_k)=\varphi(y_{k})+2|x_k-y_{k}|$ 
for all $k\in I_D$, and since the opposite inequality is guaranteed 
due to the fact that  $\varphi$ is $2$-Lipschitz, \eqref{claim} follows.
	Now, for all $\eps>0$ we can find $m_\eps$ so that $\Psi_{m}(x)+\eps\geq \varphi(x)$ for all $x\in \overline\Om$ and  $m\geq m_\eps$.
	We compute	\begin{align*}
		\phi_k(x)&=\varphi(y_{k})+2|x_k-y_{k}|-2|x_k-x|\leq \eps+\Psi_m(x)+2|x_k-y_{k}|-2|x_k-x|\\
		&\leq \eps+\Psi_{m+1}(x)\leq \eps+\varphi(x),
	\end{align*}
	where the last but one inequality follows from the definition of $\Psi_{m+1}$. This implies $\overline\varphi(x)\leq \eps+\varphi(x)$ which, by the arbitrariness of $\eps>0$, implies the claim.
\end{proof}

\begin{Remark}\label{finalremark}
{\rm If one knows in advance the regularity
result
$$
\Vert T\Vert_{{\rm flat},\alpha} = 
\min
\big\{
|R|_\Om + 
\alpha^{-1} |S|_\Om: 
(R,S) \in \zerocurrfin \times  {\mathcal S},
T=
R+\partial S\},
$$
since $\Vert T\Vert_{{\rm flat},\alpha} = \max_{\substack{\varphi\in {\rm Lip}_0(\Om)
		\\
		\| \varphi\|_{{\rm Lip}_0,\alpha}\leq 1}}\langle T
,\varphi\rangle = 
\langle T, \overline \varphi\rangle$, it is not 
difficult to check
that a maximizing $\overline\varphi$ satisfies the properties
of the function $\varphi$ in the proof of Proposition \ref{minimal_lip}. }
\end{Remark}

\textbf{Acknowledgements:}
We are grateful to Andrea Marchese for stimulating discussions and advices.
The authors are members of the Gruppo Nazionale per l'Analisi Matematica, la Probabilit\`a e le loro Applicazioni (GNAMPA) of the Istituto Nazionale di Alta Matematica (INdAM) (RS joins the project CUP\_E53C22001930001). RS also acknowledges the partial financial support of the F-CUR project number $\textrm{2262-2022-SR-CONRICMIUR\_PC-FCUR2022\_002}$  of the University of Siena. 


\begin{thebibliography}{99}
%
%
\bibitem{ADM} E. Acerbi and G. Dal Maso, {\it New lower semicontinuity results for polyconvex integrals}, Calc.
Var. Partial Differential Equations 2 (1994), 329--371.

\bibitem{AmFuPa:00}
L. Ambrosio, N. Fusco and D. Pallara, ``Functions of Bounded Variation and
Free Discontinuity Problems'', Oxford Mathematical
Monographs, Oxford University Press, New York, 2000.

\bibitem{BCS}{ G.~Bellettini, S. Carano and R.~Scala},
{\it The relaxed area of $\Suno$-valued singular maps
	in the strict $BV$-convergence 
},
ESAIM: Control, Optimization and Calculus of Variations
{\bf 28} (2022), 1-38.

\bibitem{BCS2}{ G.~Bellettini, S. Carano and R.~Scala},
{\it Relaxed area of graphs of piecewise Lipschitz maps in the strict
	$BV$-convergence},
Preprint (2023).


\bibitem{BEPS}
G. Bellettini, A. Elshorbagy, M. Paolini and  R. Scala,
{\it 
	On the relaxed area of the graph of discontinuous maps from the plane to the 
	plane taking three values with no symmetry assumptions}, 
Ann. Mat. Pura Appl. {\bf 199}, 445–477 (2020).

\bibitem{BES}{ G.~Bellettini, A.~Elshorbagy and R.~Scala}
{\it The $L^1$-relaxed area of the graph of the vortex map}, 
submitted. Preprint arXiv 2107.07236, https://arxiv.org/abs/2107.07236 (2021).

\bibitem{BMS}{ G.~Bellettini, R. Marziani and R.~Scala},
{\it 
	A non-parametric Plateau problem with partial free boundary},  submitted. Preprint arXiv 2201.06145, https://arxiv.org/abs/2201.06145 (2022).


\bibitem{BP}
G. Bellettini and M. Paolini,
{\it  On the area of the graph of a singular
	map from the plane to the plane taking three values},
Adv. Calc.
Var. {\bf 3} (2010), 371-386.

\bibitem{BePaTe:15}
G. Bellettini, M. Paolini and L. Tealdi,
{\it On the area of the graph of a piecewise smooth map from the
	plane to the plane with a curve discontinuity},
ESAIM: Control Optim. Calc. Var.
{\bf 22} (2015), 29--63.

\bibitem{BePaTe:16}
G. Bellettini, M. Paolini and L. Tealdi,
{\it Semicartesian surfaces and  the relaxed area of
	maps from the plane to the plane with a line discontinuity},
Ann. Mat. Pura Appl.
{\bf 195} (2016), 2131--2170.




\bibitem{BZ}
F. Bethuel, X. Zheng, \emph{
Density of Smooth Functions
between Two Manifolds in Sobolev Spaces}, J. Funct. Analysis {\bf 80} (1988),
60--75. 

 \bibitem{BBM}
 J. Bourgain, H. Brezis, P. Mironescu, \emph{$H^{1/2}$ maps with values into the circle: minimal connections, lifting and Ginzburg–Landau equation,}
 Publ. Math. Inst. Hautes Etudes Sci. {\bf 99} (2004), 1--115.
 
 
 \bibitem{BMbook}
 H. Brezis, P. Mironescu, 
 ``Sobolev Maps to the Circle'', Birkh\"auser New York, 2021.
 
\bibitem{BMP}
H. Brezis, P. Mironescu, A. C. Ponce, 
\emph{$W^{1,1}$ maps with values into $\Suno$}, 
Contemp. Math. {\bf 368} (2005), 69--100.

\bibitem{BMP2}
H. Brezis, P. Mironescu, A. C. Ponce, 
\emph{Complements to the paper ``$W^{1,1}$ maps with values into $\Suno$''},
preprint.

\bibitem{DaIg:03}
J. D\`avila, R. Ignat, \emph{Lifting of $BV$ functions with values in 
$\Suno$}, C. R. Acad. Sci. Paris, Ser. I {\bf 337} (2003) 159--164.

\bibitem{DG}
E. De Giorgi, \emph{On the relaxation of functionals defined on cartesian manifolds}, in “Developments in Partial Differential Equations
and Applications in Mathematical Physics” (Ferrara 1992), Plenum Press, New York (1992).

\bibitem{DLSVG}
L. De Luca, R. Scala, N. Van Goethem, \emph{A new approach to topological singularities via a weak notion of Jacobian for functions of bounded variation}, To appear on Indiana Univ. Math. J. (2022).

\bibitem{Federer}
	{H.~Federer}, {``Geometric Measure Theory''}, 
Die Grundlehren der mathematischen Wissenschaften, Vol. 153, Springer-Verlag, New York Inc., New York, (1969).




\bibitem{GMS2}
M. Giaquinta, G. Modica and J. Sou\v{c}ek, ``Cartesian Currents in the
Calculus of Variations I'', vol. 37, Springer-Verlag, Berlin, 1998.

\bibitem{GMS3} { M. Giaquinta, G. Modica and J. Sou\u{c}ek}, 
``Cartesian Currents in the Calculus of Variations II. Variational Integrals'',
Ergebnisse der Mathematik und ihrer Grenzgebiete, Vol. 38,
Springer-Verlag, Berlin-Heidelberg, 1998.

\bibitem{Giu:84} { E. Giusti},
``Minimal Surfaces and Functions of Bounded Variation'',
Birkh\"auser, Boston
1984.

\bibitem{I} R. Ignat,
\emph{The space $BV(\mathbb S^2,\Suno)$: minimal connection and 
optimal lifting}, Ann. I. H. Poincar\`e, Anal. Nonlineaire {\bf 22}
(2005), 283--302.
%
\bibitem{Krantz_Parks} S.G. Krantz, H.R. Parks,
``Geometric Integration Theory'',
Cornerstones, Birk\"auser Boston, Inc., Boston, MA, 2008

\bibitem{Mu}
{D. Mucci}, \textit{Strict convergence with equibounded area and minimal completely vertical liftings},{ 
	Nonlinear Anal.} { \bf 221} {(2022)}.




\bibitem{Ponce} A. C. Ponce, \emph{On the distributions of the form $\sum_i(\delta_{p_i}-\delta_{n_i})$}, J. Func. Anal. {\bf 210} (2004),
391-435.




\bibitem{S}
R. Scala, {\it Optimal estimates for the triple junction function and other surprising aspects of the
	area functional}, Ann. Sc. Norm. Super. Pisa Cl. Sci. (5) {\bf 20}(2), 491-564. (2020).
\bibitem{SS}
R. Scala, G. Scianna, {\it On the $L^1$-relaxed area of graphs of
	BV piecewise constant maps}, In preparation.

\end{thebibliography}
\end{document}